\documentclass[a4paper,12pt,twoside]{article}
\usepackage{amsmath,amssymb,ifthen}
\usepackage[amsmath,thmmarks]{ntheorem}
\usepackage{paper}
\usepackage[all]{xy}
\usepackage{mathrsfs}
\usepackage{url}
\newcommand{\red}{\mathrm{red}}
\DeclareMathOperator{\Spf}{Spf}
\DeclareMathOperator{\Spa}{Spa}
\DeclareMathOperator{\GSp}{GSp}
\DeclareMathOperator{\Sp}{Sp}
\DeclareMathOperator{\Res}{Res}
\DeclareMathOperator{\Nilp}{\mathbf{Nilp}}
\DeclareMathOperator{\Ind}{Ind}
\DeclareMathOperator{\cInd}{c-Ind}
\DeclareMathOperator{\supp}{supp}
\DeclareMathOperator{\Zel}{Zel}
\DeclareMathOperator{\Rep}{\mathbf{Rep}}
\DeclareMathOperator{\Mod}{\mathbf{Mod}}
\DeclareMathOperator{\spp}{sp}
\DeclareMathOperator{\Tor}{Tor}
\DeclareMathOperator{\Tot}{Tot}
\DeclareMathOperator{\D}{D}
\newcommand{\Lotimes}{\overset{\mathbb{L}}{\otimes}}
\newcommand{\M}{\mathcal{M}}
\newcommand{\X}{\mathbb{X}}
\newcommand{\qcet}{\mathrm{qc\acute{e}t}}
\newcommand{\cohet}{\mathrm{coh\acute{e}t}}

\SelectTips{cm}{11}
\title{Zelevinsky involution and $\ell$-adic cohomology of the Rapoport-Zink tower}
\author{Yoichi Mieda}
\def\headertitle{Zelevinsky involution and $\ell$-adic cohomology of the Rapoport-Zink tower}

\begin{document}
\maketitle

\begin{firstfootnote}
 The Hakubi Center for Advanced Research / Department of Mathematics, Kyoto University, Kyoto, 606--8502, Japan

 E-mail address: \texttt{mieda@math.kyoto-u.ac.jp}

 2010 \textit{Mathematics Subject Classification}.
 Primary: 11F70;
 Secondary: 14G35, 22E50.
\end{firstfootnote}

\begin{abstract}
 In this paper, we investigate how the Zelevinsky involution appears in the $\ell$-adic cohomology
 of the Rapoport-Zink tower. We generalize the result of Fargues on the Drinfeld tower to
 the Rapoport-Zink towers for symplectic similitude groups.
\end{abstract}

\section{Introduction}

The non-abelian Lubin-Tate theory says that the local Langlands correspondence for $\GL(n)$ is
geometrically realized in the $\ell$-adic cohomology of the Lubin-Tate tower and the Drinfeld tower
(\cf \cite{MR1044827}, \cite{MR1464867}, \cite{MR1876802}).
It urges us to consider how representation-theoretic operations are translated into geometry.
In \cite{Fargues-Zel}, Fargues found a relation between the Zelevinsky involution and
the Poincar\'e duality of the $\ell$-adic cohomology of the Drinfeld tower.
Furthermore, by using Faltings' isomorphism between the Lubin-Tate tower and the Drinfeld tower
(\cf \cite{MR1936369}, \cite{MR2441311}), he obtained a similar result for the Lubin-Tate tower.
This result is useful for study of the $\ell$-adic cohomology itself.
For example, it played a crucial role in Boyer's work \cite{MR2511742},
which completely determined the $\ell$-adic cohomology of the Lubin-Tate tower.

A Rapoport-Zink tower (\cf \cite{MR1393439}) is a natural generalization of the Lubin-Tate tower
and the Drinfeld tower.
It is a projective system of \'etale coverings of rigid spaces $\{M_K\}$ lying over the rigid generic
fiber $M$ of a formal scheme $\M$.
The formal scheme $\M$, called a Rapoport-Zink space, 
is defined as a moduli space of deformations by quasi-isogenies
of a $p$-divisible group over $\overline{\F}_p$ with additional structures.
For a prime number $\ell\neq p$, consider the compactly supported
$\ell$-adic cohomology $H^i_c(M_\infty)=\varinjlim_K H^i_c(M_K,\overline{\Q}_\ell)$.
It is naturally endowed with an action of $G\times J\times W$, where $G$ and $J$ are
$p$-adic reductive groups and $W$ is the Weil group of some $p$-adic field
(a local analogue of a reflex field).
The cohomology $H^i_c(M_\infty)$ is expected to be described by the local Langlands correspondence
for $G$ and $J$ (\cf \cite{MR1403942}), but only few results are known.

In this paper, we will give a generalization of Fargues' result mentioned above to Rapoport-Zink
towers other than the Lubin-Tate tower and the Drinfeld tower.
Although our method should be valid for many Rapoport-Zink towers
(see Remark \ref{rem:generalizations}),
here we restrict ourselves to the case of $\GSp(2n)$ for the sake of simplicity.
In this case, the Rapoport-Zink tower is a local analogue of the Siegel modular variety,
and is also treated in \cite{RZ-LTF}. The group $G$ is equal to $\GSp_{2n}(\Q_p)$, and
$J$ is an inner form $\mathrm{GU}(n,D)$ of $G$, where $D$ is the quaternion division algebra
over $\Q_p$.

The main difference between Fargues' case and ours is that
the Rapoport-Zink space $\M$ is a $p$-adic formal scheme in the former, while not in the latter.
Owing to this difference, we need to introduce a new kind of cohomology $H^i_{\mathcal{C}_\M}(M_\infty)$.
Contrary to the compactly supported cohomology $H^i_c(M_\infty)$, this cohomology depends on
the formal model $\M$ of the base $M$ of the Rapoport-Zink tower. Roughly speaking,
it is a cohomology with compact support in the direction of the formal model $\M$
(for a precise definition, see Section \ref{subsec:formal-adic} and Section \ref{subsec:RZ-tower}).
If $\M$ is a $p$-adic formal scheme, it coincides with the compactly supported cohomology.
By using these two cohomology, our main theorem is stated as follows:

\begin{thm}[Theorem \ref{thm:RZ-Zel}]\label{thm:main-intro}
 Let $\widetilde{K}$ be an open compact-mod-center subgroup of $G$
 and $\tau$ an irreducible smooth representation of $\widetilde{K}$.
 Denote by $\chi$ the central character of $\tau^\vee$.
 For a smooth $G$-representation $V$, 
 put $V_\tau=\Hom_{\widetilde{K}}(\tau,V\otimes_{\mathcal{H}(Z_G)}\chi^{-1})$,
 where $\mathcal{H}(Z_G)$ is the Hecke algebra of the center $Z_G$ of $G$.
 Let $\mathfrak{s}\in I_\chi$ be a Bernstein component of the category of smooth representations
 of $J$ with central character $\chi$.
 An integer $\iota(\mathfrak{s})$ is naturally attached to $\mathfrak{s}$; $\mathfrak{s}$ is
 supercuspidal if and only if $\iota(\mathfrak{s})=0$ (\cf Section \ref{sec:Zelevinsky}).

 Assume that the $\mathfrak{s}$-component $H^q_c(M_\infty)_{\tau,\mathfrak{s}}$ of
 $H^q_c(M_\infty)_\tau$ is a finite length $J$-representation for every integer $q$.
 Then, for each integer $i$, we have an isomorphism of $J\times W_{\Q_p}$-representations
 \[
  H^{2d+\iota(\mathfrak{s})-i}_{\mathcal{C}_\M}(M_\infty)_{\tau^\vee,\mathfrak{s}^\vee}(d)\cong
 \Zel\bigl(H^i_c(M_\infty)_{\tau,\mathfrak{s}}^\vee\bigr).
 \]
 Here $W_{\Q_p}$ denotes the Weil group of $\Q_p$, $d=n(n+1)/2$ the dimension of $M$,
 and $\Zel$ the Zelevinsky involution with respect to $J$ (see Section \ref{sec:Zelevinsky}).
\end{thm}

By applying this theorem to the case where $\cInd_{\widetilde{K}}^G\tau$ becomes supercuspidal,
we obtain the following consequence on the supercuspidal part of $H^i_c(M_\infty/p^\Z)$
(here $p^\Z$ is regarded as a discrete subgroup of the center of $J$).

\begin{cor}[Corollary \ref{cor:vanishing-noncusp-part}, Corollary \ref{cor:A-packet}]\label{cor:application-intro}
 Let $\pi$ be an irreducible supercuspidal representation of $G$, $\rho$ an irreducible
 non-supercuspidal representation of $J$, and $\sigma$ an irreducible $\ell$-adic representation
 of $W_{\Q_p}$. 
 Under some technical assumption (Assumption \ref{assump:finiteness}), the following hold.
 \begin{enumerate}
  \item The representation $\pi\otimes\rho$ does not appear as a subquotient
	of $H^{d-\dim\M^\red}_c(M_\infty/p^\Z)$. In particular, if $n=2$, 
	$\pi\otimes\rho$ does not appear as a subquotient of $H^2_c(M_\infty/p^\Z)$.
  \item If $n=2$, $\pi\otimes\rho\otimes\sigma$ appears as a subquotient of $H^3_c(M_\infty/p^\Z)$
	if and only if $\pi^\vee\otimes\Zel(\rho^\vee)\otimes\sigma^\vee(-3)$
	appears as a subquotient of $H^4_c(M_\infty/p^\Z)$.
 \end{enumerate}
\end{cor}
The proof of ii) also requires a result of \cite{RZ-GSp4}, which measures the difference between 
the two cohomology groups $H^i_c(M_\infty/p^\Z)$ and $H^i_{\mathcal{C}_\M}(M_\infty/p^\Z)$.
This corollary will be very useful to investigate how non-supercuspidal representations of $J$
contribute to $H^i_c(M_\infty/p^\Z)$ for each $i$.

The outline of this paper is as follows.
In Section \ref{sec:Zelevinsky}, we briefly recall the definition and properties of
the Zelevinsky involution. The main purpose of this section is to fix notation,
and most of the proofs are referred to \cite{MR1471867} and \cite{Fargues-Zel}.
Section \ref{sec:preliminaries} is devoted to give some preliminaries on algebraic and rigid geometry
used in this article. The cohomology appearing in Theorem \ref{thm:main-intro} is defined in
this section. Because we prefer to use the theory of adic spaces
(\cf \cite{MR1306024}, \cite{MR1734903}) as a framework of rigid geometry,
we need to adopt the theory of smooth equivariant sheaves for Berkovich spaces
developed in \cite[\S IV.9]{MR2441311} to adic spaces. In fact proofs become simpler;
especially the compactly supported cohomology and the Godement resolution can be easily treated.
In Section \ref{sec:duality}, we prove a duality theorem under a general setting.
Finally in Section \ref{sec:application}, after introducing the Rapoport-Zink tower for $\GSp(2n)$,
we deduce Theorem \ref{thm:main-intro} from the duality theorem proved in Section \ref{sec:duality}.
We also give the applications announced in Corollary \ref{cor:application-intro}.

\medbreak
\noindent{\bfseries Acknowledgment}\quad
This work was supported by JSPS KAKENHI Grant Number 24740019.

\medbreak
\noindent{\bfseries Notation}\quad
For a field $k$, we denote its algebraic closure by $\overline{k}$.

\section{Zelevinsky involution}\label{sec:Zelevinsky}
In this section, we recall briefly about the Zelevinsky involution.
See \cite[\S 1]{Fargues-Zel} for details.

Let $p$ be a prime number and $G$ be a locally pro-$p$ group; namely, $G$ is a locally compact group
which has an open pro-$p$ subgroup.
Fix a $0$-dimensional Gorenstein local ring $\Lambda$ in which $p$ is invertible.
We write $\Rep_\Lambda(G)$ for the category of smooth representations of $G$ over $\Lambda$.
Denote by $G^{\mathrm{disc}}$ the group $G$ with the discrete topology.
We have a natural functor $i_G\colon \Rep_{\Lambda}(G)\longrightarrow \Rep_{\Lambda}(G^\mathrm{disc})$,
which has a right adjoint functor $\infty_G\colon \Rep_\Lambda(G^{\mathrm{disc}})\longrightarrow \Rep_\Lambda(G)$;
$V\longmapsto \varinjlim_{K}V^K$.
Here $K$ runs through compact open subgroups of $G$. The functor $\infty_G$ is not exact in general.

Let $\mathcal{D}_c(G)$ be the convolution algebra of compactly supported $\Lambda$-valued
distributions on $G$. It contains the Hecke algebra $\mathcal{H}(G)$ of $G$ consisting of
compactly supported distributions invariant under some compact open subgroup of $G$.
For each open pro-$p$ subgroup $K$ of $G$, an idempotent $e_K$ of $\mathcal{H}(G)$
is naturally attached.
We denote by $\Mod(\mathcal{D}_c(G))$ the category of $\mathcal{D}_c(G)$-modules.
We have a natural functor $i_{\mathcal{D}}\colon \Rep_\Lambda(G)\longrightarrow \Mod(\mathcal{D}_c(G))$. 
The right adjoint functor $\infty_{\mathcal{D}}\colon \Mod(\mathcal{D}_c(G))\longrightarrow \Rep_\Lambda(G)$ of
$i_{\mathcal{D}}$ is given by $M\longmapsto \varinjlim_{K}e_KM$, where $K$ runs through compact open
pro-$p$ subgroups of $G$. Note that $\infty_{\mathcal{D}}$ is an exact functor.

For a compact open subgroup $K$ of $G$, we have a functor
\[
 \cInd_{\mathcal{D}_c(K)}^{\mathcal{D}_c(G)}\colon \Mod(\mathcal{D}_c(K))\longrightarrow \Mod(\mathcal{D}_c(G));\
 M\longmapsto \mathcal{D}_c(G)\otimes_{\mathcal{D}_c(K)}M.
\] 
This functor is exact and the following diagrams are commutative:
\[
 \xymatrix{%
 \Rep_{\Lambda}(K)\ar[rr]^-{\cInd_K^G}\ar[d]^-{i_{\mathcal{D}}}&& \Rep_{\Lambda}(G)\ar[d]^-{i_{\mathcal{D}}}\\
 \Mod\bigl(\mathcal{D}_c(K)\bigr)\ar[rr]^-{\cInd_{\mathcal{D}_c(K)}^{\mathcal{D}_c(G)}}&&
 \Mod\bigl(\mathcal{D}_c(G)\bigr)\lefteqn{,}
 }\
 \xymatrix{%
 \Mod\bigl(\mathcal{D}_c(K)\bigr)\ar[rr]^-{\cInd_{\mathcal{D}_c(K)}^{\mathcal{D}_c(G)}}\ar[d]^-{\infty_{\mathcal{D}}}&&
 \Mod\bigl(\mathcal{D}_c(G)\bigr)\ar[d]^-{\infty_{\mathcal{D}}}\\
 \Rep_{\Lambda}(K)\ar[rr]^-{\cInd_K^G}&&
 \Rep_{\Lambda}(G)\lefteqn{.}
 }
\] 
Let us observe the commutativity of the right diagram, as it is not included in \cite{Fargues-Zel}.
Take a system of representatives $\{g_i\}_{i\in I}$ of $G/K$. Then, as in \cite[\S 1.4]{Fargues-Zel},
we have $\mathcal{D}_c(G)=\bigoplus_{i\in I}\delta_{g_i^{-1}}*\mathcal{D}_c(K)$,
where $\delta_{g_i^{-1}}$ denotes the Dirac distribution at $g_i^{-1}\in G$.
Therefore, for a $\mathcal{D}_c(K)$-module $M$, 
an element $x$ of $\cInd_{\mathcal{D}_c(K)}^{\mathcal{D}_c(G)}M$ can be written uniquely
in the form $\sum_{i\in I}\delta_{g_i^{-1}}\otimes x_i$ with $x_i\in M$.
Put $M_\infty=i_{\mathcal{D}}\infty_{\mathcal{D}}M$. It is the $\mathcal{D}_c(K)$-submodule of $M$
consisting of $x\in M$ such that $e_{K'}x=x$ for some open pro-$p$ subgroup $K'$ of $G$.
By the left diagram and the fact that $i_{\mathcal{D}}$ is fully faithful,
it suffices to prove that $\cInd_{\mathcal{D}_c(K)}^{\mathcal{D}_c(G)}M_\infty=(\cInd_{\mathcal{D}_c(K)}^{\mathcal{D}_c(G)}M)_\infty$.
First take $x=\sum_{i\in I}\delta_{g_i^{-1}}\otimes x_i$ in
$\cInd_{\mathcal{D}_c(K)}^{\mathcal{D}_c(G)}M_\infty$. For each $i\in I$ with $x_i\neq 0$,
take an open pro-$p$ subgroup $K_i$ of $K$ such that $e_{K_i}x_i=x_i$.
We can find an open pro-$p$ subgroup $K'$ of $G$ such that $g_iK'g_i^{-1}\subset K_i$ for every
$i\in I$ with $x_i\neq 0$.
Then, we have $e_{K'}x=\sum_{i\in I,x_i\neq 0}\delta_{g_i^{-1}}\otimes e_{g_iK'g_i^{-1}}x_i=x$,
and thus $x\in (\cInd_{\mathcal{D}_c(K)}^{\mathcal{D}_c(G)}M)_\infty$.
Next, take $x=\sum_{i\in I}\delta_{g_i^{-1}}\otimes x_i$ in 
$(\cInd_{\mathcal{D}_c(K)}^{\mathcal{D}_c(G)}M)_\infty$.
Then, there exists an open pro-$p$ subgroup $K'$ of $G$ such that $e_{K'}x=x$.
We may shrink $K'$ so that $g_iK'g_i'^{-1}\subset K$ for every $i\in I$ with $x_i\neq 0$.
Then, we have $x=e_{K'}x=\sum_{i\in I,x_i\neq 0}\delta_{g_i^{-1}}\otimes e_{g_iK'g_i^{-1}}x_i$.
Hence $e_{g_iK'g_i^{-1}}x_i$ is equal to $x_i$ for every $i\in I$ with $x_i\neq 0$,
which implies $x_i\in M_\infty$. Thus $x\in \cInd_{\mathcal{D}_c(K)}^{\mathcal{D}_c(G)}M_\infty$.
Now we conclude that $\cInd_{\mathcal{D}_c(K)}^{\mathcal{D}_c(G)}M_\infty=(\cInd_{\mathcal{D}_c(K)}^{\mathcal{D}_c(G)}M)_\infty$.

\begin{defn}
 For $\pi\in\Rep_\Lambda(G)$, consider $\Hom_G(\pi,\mathcal{H}(G))$, where $\mathcal{H}(G)$ is regarded
 as a smooth representation of $G$ by the left translation. As $\mathcal{H}(G)$ has another smooth 
 $G$-action by the right translation, $\Hom_G(\pi,\mathcal{H}(G))$ has a structure of 
 a $\mathcal{D}_c(G)$-module. Therefore we get a contravariant functor from $\Rep_\Lambda(G)$ to
 $\Mod(\mathcal{D}_c(G))$, for which we write $\D^m$.
 Composing with $\infty_{\mathcal{D}}$, we obtain a contravariant functor
 $\D=\infty_{\mathcal{D}}\circ \D^m\colon \Rep_\Lambda(G)\longrightarrow \Rep_\Lambda(G)$.
 We denote by $R\!\D^m$ (resp.\ $R\!\D$) the right derived functor of $\D^m$ (resp.\ $\D$).
 As $\infty_{\mathcal{D}}$ is exact, we have $R\!\D=\infty_{\mathcal{D}}\circ R\!\D^m$.
\end{defn}

\begin{prop}\label{prop:dual-cInd}
 Let $K$ be an open pro-$p$ subgroup of $G$
 and $\rho$ a smooth representation of $K$ over $\Lambda$. Then there is
 a natural $\mathcal{D}_c(G)$-linear injection
 \[
 \cInd_{\mathcal{D}_c(K)}^{\mathcal{D}_c(G)}(\rho^*)\hooklongrightarrow \D^m(\cInd_K^G\rho).
 \]
 Here $\rho^*=\Hom_K(\rho,\Lambda)$ denotes the algebraic dual, which is naturally equipped with
 a structure of a $\mathcal{D}_c(K)$-module.
 Applying $\infty_\mathcal{D}$ to the injection above, we obtain a $G$-equivariant injection
 \[
  \cInd_K^G(\rho^\vee)\hooklongrightarrow \D(\cInd_K^G\rho),
 \]
 where $\rho^\vee=\infty_{\mathcal{D}}(\rho^*)=\infty_K(\Hom_K(\rho,\Lambda))$ denotes the contragredient
 representation of $\rho$.

 If moreover $\rho$ is a finitely generated $K$-representation
 (in other words, $\rho$ has finite length as a $\Lambda$-module), then we have
 \begin{gather*}
  R\!\D^m\bigl(\cInd_K^G(\rho)\bigr)=\D^m\bigl(\cInd_K^G(\rho)\bigr)=\cInd_{\mathcal{D}_c(K)}^{\mathcal{D}_c(G)}(\rho^*),\\
  R\!\D\bigl(\cInd_K^G(\rho)\bigr)=\D\bigl(\cInd_K^G(\rho)\bigr)=\cInd_K^G(\rho^\vee).
 \end{gather*}
\end{prop}

\begin{prf}
 See \cite[Lemme 1.10, Lemme 1.12, Lemme 1.13]{Fargues-Zel}.
\end{prf}

\begin{cor}\label{cor:derived-dual}
 Assume that $\Lambda$ is a field, $\Rep_\Lambda(G)$ is noetherian and has finite projective dimension.
 \begin{enumerate}
  \item Let $D^b_{\mathrm{fg}}(\Rep_\Lambda(G))$ be the full subcategory of $D^b(\Rep_\Lambda(G))$ consisting of
	complexes whose cohomology are finitely generated $G$-representations.
	Then the contravariant functor $R\!\D$ maps $D^b_{\mathrm{fg}}(\Rep_\Lambda(G))$ into itself.
  \item The contravariant functor
	\[
	 R\!\D\colon D^b_{\mathrm{fg}}\bigl(\Rep_\Lambda(G)\bigr)\longrightarrow D^b_{\mathrm{fg}}\bigl(\Rep_\Lambda(G)\bigr)
	\]
	satisfies $R\!\D\circ R\!\D\cong \id$.
  \item For a field extension $\Lambda'$ of $\Lambda$, the following diagram is $2$-commutative:
	\[
	 \xymatrix{%
	D^b_{\mathrm{fg}}\bigl(\Rep_\Lambda(G)\bigr)\ar[r]^-{R\!\D}\ar[d]&
	D^b_{\mathrm{fg}}\bigl(\Rep_\Lambda(G)\bigr)\ar[d]\\
	D^b_{\mathrm{fg}}\bigl(\Rep_{\Lambda'}(G)\bigr)\ar[r]^-{R\!\D}&
	D^b_{\mathrm{fg}}\bigl(\Rep_{\Lambda'}(G)\bigr)\lefteqn{.}
	}
	\]
	Here the vertical arrows denote the base change functor.
 \end{enumerate}
\end{cor}

\begin{prf}
 For i) and ii), see \cite[Proposition 1.18]{Fargues-Zel}. iii) follows immediately from 
 Proposition \ref{prop:dual-cInd}.
\end{prf}

In the remaining part of this section, we assume that $\Lambda=\C$.
Let $F$ be a $p$-adic field and $\mathbf{G}$ a connected reductive group over $F$.
Write $G$ for $\mathbf{G}(F)$. 
Fix a discrete cocompact subgroup $\Gamma$ of the center $Z_G=\mathbf{Z}_{\mathbf{G}}(F)$
and put $G'=G/\Gamma$. We simply write $\Rep(G')$ for $\Rep_{\C}(G')$.
Note that we have a natural decomposition of a category
\[
 \Rep(G')=\prod_{\chi\colon Z_G/\Gamma\to\C^\times}\Rep_{\chi}(G),
\]
where $\chi$ runs through smooth characters of the compact group $Z_G/\Gamma$,
and $\Rep_{\chi}(G)$ denotes the category of smooth representations of $G$ with central character $\chi$.
We will apply the theory above to the group $G'$.
In this case all the assumptions in Corollary \ref{cor:derived-dual} are satisfied
(\cf \cite[Remarque 3.12]{MR771671}, \cite[Proposition 37]{MR1159104}).

Let $I_\Gamma$ be the set of inertially equivalence classes of cuspidal data $(\mathbf{M},\sigma)$ such that
$\sigma\vert_{\Gamma}$ is trivial. We have the Bernstein decomposition (\cf \cite[Th\'eor\`eme VI.7.2]{MR2567785})
\[
 \Rep(G')=\prod_{\mathfrak{s}\in I_\Gamma}\Rep(G')_{\mathfrak{s}}.
\]
For $V\in \Rep(G')$, we denote the corresponding decomposition by $V=\bigoplus_{\mathfrak{s}\in I_\Gamma}V_{\mathfrak{s}}$.

\begin{prop}
 For $\mathfrak{s}=[(\mathbf{M},\sigma)]\in I_{\Gamma}$, 
 put $\mathfrak{s}^\vee=[(\mathbf{M},\sigma^\vee)]\in I_{\Gamma}$.
 Then, $R\!\D$ induces a contravariant functor 
 $D^b_{\mathrm{fg}}(\Rep(G')_{\mathfrak{s}})\longrightarrow D^b_{\mathrm{fg}}(\Rep(G')_{\mathfrak{s}^\vee})$.
\end{prop}

\begin{prf}
 See \cite[Remarque 1.5]{Fargues-Zel}.
\end{prf}

For $\mathfrak{s}=[(\mathbf{M},\sigma)]\in I_{\Gamma}$, put
$\iota(\mathfrak{s})=r_{\mathbf{G}}-r_{\mathbf{M}}$, where $r_{\mathbf{G}}$ (resp.\ $r_{\mathbf{M}}$)
denotes the split semisimple rank of $\mathbf{G}$ (resp.\ $\mathbf{M}$).
The number $\iota(\mathfrak{s})$ is $0$ if and only if $\mathbf{M}=\mathbf{G}$.

\begin{thm}[{{\cite[Theorem III.3.1]{MR1471867}}}]\label{thm:Schneider-Stuhler}
 Fix $\mathfrak{s}\in I_{\Gamma}$.
 Let $\Rep(G')_{\mathfrak{s}}^{\mathrm{fl}}$ be the full subcategory of $\Rep(G')_{\mathfrak{s}}$
 consisting of representations of finite length.
 For $\pi\in \Rep(G')_{\mathfrak{s}}^{\mathrm{fl}}$, we have
 $R^i\!\D(\pi)=0$ if $i\neq \iota(\mathfrak{s})$. Moreover, $R^{\iota(\mathfrak{s})}\!\D(\pi)$ has finite length.
\end{thm}

\begin{defn}
 For $\mathfrak{s}\in I_{\Gamma}$ and $\pi\in \Rep(G')_{\mathfrak{s}}^{\mathrm{fl}}$,
 put $\Zel(\pi)=R^{\iota(\mathfrak{s})}\!\D(\pi^\vee)$. 
 The (covariant) functor $\Zel\colon \Rep(G')_{\mathfrak{s}}^{\mathrm{fl}}\longrightarrow \Rep(G')_{\mathfrak{s}}^{\mathrm{fl}}$ is called the Zelevinsky involution.
 It is an exact categorical equivalence. In particular, it preserves irreducibility.
\end{defn}

\begin{prop}
 For an irreducible smooth representation $\pi$ of $G'$, we have 
 an isomorphism $\Zel(\pi^\vee)\cong \Zel(\pi)^\vee$. In particular, $\Zel(\Zel(\pi))\cong \pi$.
\end{prop}

\begin{prf}
 It is an immediate consequence of \cite[Proposition IV.5.4]{MR1471867}.
\end{prf}

 Let $\chi\colon Z_G\longrightarrow \C^\times$ be a smooth character, which is not necessarily unitary.
 We can consider a variant of $\Zel$ on $\Rep_\chi(G)$ as follows.
 Let $\mathcal{H}_\chi(G)$ be the set of locally constant $\C$-valued functions $f$ such that
 \begin{itemize}
  \item $f(zg)=\chi(z)^{-1}f(g)$ for every $z\in Z_G$ and $g\in G$,
  \item and $\supp f$ is compact modulo $Z_G$.
\end{itemize}
Let $\D_\chi\colon \Rep_\chi(G)\longrightarrow \Rep_{\chi^{-1}}(G)$ be the contravariant functor defined by
\[
 \D_\chi(\pi)=\Hom_{\Rep_{\chi}(G)}\bigl(\pi,\mathcal{H}_\chi(G)\bigr),
\]
and $R\!\D_{\chi}$ be the derived functor of $\D_\chi$. As in Corollary \ref{cor:derived-dual} i), $R\!\D_{\chi}$
induces a contravariant functor
\[
 R\!\D_{\chi}\colon D^b_{\mathrm{fg}}\bigl(\Rep_\chi(G)\bigr)\longrightarrow
 D^b_{\mathrm{fg}}\bigl(\Rep_{\chi^{-1}}(G)\bigr).
\]
Let $I_\chi$ be the set of inertially equivalence classes of cuspidal data $(\mathbf{M},\sigma)$
such that $\sigma\vert_{Z_G}=\chi$. We have the Bernstein decomposition
\[
 \Rep_\chi(G)=\prod_{\mathbf{s}\in I_\chi}\Rep_\chi(G)_{\mathfrak{s}}.
\]
Let $\Rep_\chi(G)_{\mathfrak{s}}^{\mathrm{fl}}$ be the full subcategory of $\Rep_\chi(G)_{\mathfrak{s}}$
consisting of representations of finite length.
By \cite[Theorem III.3.1]{MR1471867}, for $\pi\in\Rep_\chi(G)_{\mathfrak{s}}^{\mathrm{fl}}$ we have
$R^i\!\D_\chi(\pi)=0$ if $i\neq \iota(\mathfrak{s})$ and $R^{\iota(\mathfrak{s})}\!\D_\chi(\pi)$
has finite length. Hence we can give the following definition.

\begin{defn}\label{defn:Zel-variant}
 For $\mathfrak{s}\in I_\chi$ and $\pi\in \Rep_\chi(G)_{\mathfrak{s}}^{\mathrm{fl}}$,
 put $\Zel_\chi(\pi)=R^{\iota(\mathfrak{s})}\!\D_{\chi^{-1}}(\pi^\vee)$. 
 It induces an exact categorical equivalence 
 $\Zel_\chi\colon \Rep_\chi(G)_{\mathfrak{s}}^{\mathrm{fl}}\longrightarrow \Rep_\chi(G)_{\mathfrak{s}}^{\mathrm{fl}}$
 satisfying $\Zel_\chi^2\cong \id$.
\end{defn}

\begin{lem}\label{lem:Zel-twist}
 \begin{enumerate}
  \item For $\pi\in \Rep_\chi(G)_{\mathfrak{s}}^{\mathrm{fl}}$ and a smooth character $\omega$ of $G$,
	we have $\Zel_{\chi\otimes\omega}(\pi\otimes\omega)\cong \Zel_\chi(\pi)\otimes\omega$.
  \item Let $\Gamma\subset Z_G$ and $G'=G/\Gamma$ be as above. 
	If $\chi$ is trivial on $\Gamma$, then for every
	$\pi\in\Rep_\chi(G)_{\mathfrak{s}}^{\mathrm{fl}}$ we have $\Zel_\chi(\pi)\cong \Zel(\pi)$.
	In the right hand side, $\pi$ is regarded as an object of $\Rep(G')_{\mathfrak{s}}^{\mathrm{fl}}$.
 \end{enumerate}
\end{lem}

\begin{prf}
 i) is clear from definition. For ii), note that $\mathcal{H}(G')=\bigoplus_{\chi'}\mathcal{H}_{\chi'}(G)$,
 where $\chi'$ runs through smooth characters of $Z_G$ which are trivial on $\Gamma$.
 By the decomposition $\Rep(G')=\prod_{\chi'}\Rep_{\chi'}(G)$, we have
 \[
 R\!\D(\pi^\vee)=R\Hom_{\Rep(G')}\bigl(\pi^\vee,\mathcal{H}(G)\bigr)=R\Hom_{\Rep_{\chi^{-1}}(G)}\bigl(\pi^\vee,\mathcal{H}_{\chi^{-1}}(G)\bigr)=R\!\D_{\chi^{-1}}(\pi^\vee).
 \]
 Hence we have $\Zel_\chi(\pi)\cong \Zel(\pi)$, as desired.
\end{prf}

By this lemma, we can simply write $\Zel$ for $\Zel_\chi$ without any confusion.

\section{Preliminaries}\label{sec:preliminaries}
\subsection{Compactly supported cohomology for partially proper schemes}
In \cite[\S 5]{MR1734903}, Huber defined the compactly supported cohomology for adic spaces
which are partially proper over a field as the derived functor of $\Gamma_c$.
This construction is also applicable to schemes over a field.

\begin{defn}
 Let $f\colon X\longrightarrow Y$ be a morphism between schemes.
 \begin{enumerate}
  \item The morphism $f$ is said to be specializing if for every $x\in X$ and every specialization $y'$
	of $y=f(x)$, there exists a specialization $x'$ of $x$ such that $y'=f(x')$.
	If an arbitrary base change of $f$ is specializing, $f$ is said to be universally specializing.
  \item The morphism $f$ is said to be partially proper
	if it is separated, locally of finite type and universally specializing.
 \end{enumerate}
\end{defn}

\begin{prop}\label{prop:partially-proper}
 \begin{enumerate}
  \item A morphism of schemes is proper if and only if partially proper and quasi-compact.
  \item Partially properness can be checked by the valuative criterion.
  \item Let $f\colon X\longrightarrow Y$ be a partially proper morphism between schemes.
	Assume that $Y$ is noetherian. Then, for every quasi-compact subset $T$ of $X$,
	the closure $\overline{T}$ of $T$ is quasi-compact.
 \end{enumerate}
\end{prop}

\begin{prf}
 i) can be proved in the same way as \cite[Lemma 1.3.4]{MR1734903}.
 ii) is straightforward and left to the reader.
 Let us prove iii). We may assume that $T$ is open in $X$. Moreover we may assume that
 $T$ and $Y$ are affine. Put $T=\Spec A$ and $Y=\Spec B$.
 Let $\widetilde{B}$ be the integral closure of the image of $B\longrightarrow A$ in $A$,
 and consider the topological space $T_v=\Spa(A,\widetilde{B})$. Here we endow $A$ with the discrete topology.
 As a set, $T_v$ can be identified with the set of pairs $(x,V_x)$ where
 \begin{itemize}
  \item $x\in T$, and
  \item $V_x$ is a valuation ring of the residue field $\kappa_x$ at $x$ such that the composite
	$\Spec\kappa_x\longrightarrow T\longrightarrow Y$ can be extended to
	$\Spec V_x\longrightarrow Y$.
 \end{itemize}
 Therefore, by ii), we can construct a map $\phi\colon T_v\longrightarrow X$ as follows.
 For $(x,V_x)\in T_v$, the $Y$-morphism $\Spec\kappa_x\longrightarrow T$ uniquely extends to
 a $Y$-morphism $\Spec V_x\longrightarrow X$. We let $\phi(x,V_x)$ be the image of the closed point
 in $\Spec V_x$ under this morphism. Since $X$ is quasi-separated locally spectral and $T$ is quasi-compact open,
 each point in $\overline{T}$ is a specialization of some point in $T$ (\cite[Corollary of Theorem 1]{MR0251026}).
 Thus $\overline{T}$ coincides with $\phi(T_v)$.

 We will prove that $\phi$ is continuous. 
 Fix $(x,V_x)\in T_v$ and take an affine neighborhood $U=\Spec C$ of $y=\phi(x,V_x)$.
 We can find $u\in A$ such that $T'=\Spec A[1/u]$ is an open neighborhood of $x$ contained in $U$. 
 On the other hand, as $f$ is locally
 of finite type, $C$ is a finitely generated $B$-algebra. Take a system of generators $c_1,\ldots,c_n\in C$
 $(n\ge 1)$ and consider the images of them under the ring homomorphism $C\longrightarrow A[1/u]$ that comes from
 the inclusion $T'\hooklongrightarrow U$. There exist integers $k_1,\ldots,k_n$ and
 $a_1,\ldots,a_n\in A$ such that the image of $a_i$ in $A[1/u]$ coincides with the image of
 $u^{k_i}c_i$ under $C\longrightarrow A[1/u]$.
 Let $W$ be the open subset of $T_v$ defined by the condition
 $v(a_i)\le v(u^{k_i})\neq 0$ $(i=1,\ldots,n)$. Then it is easy to observe that
 $(x,V_x)$ belongs to $W$ and $\phi(W)$ is contained in $U$.
 This proves the continuity of $\phi$.

 By \cite[Theorem 3.5 (i)]{MR1207303}, $T_v$ is quasi-compact. Therefore we conclude that 
 $\overline{T}=\phi(T_v)$ is quasi-compact, as desired.
\end{prf}

\begin{cor}
 Let $X$ be a scheme which is partially proper over a noetherian scheme.
 Then $X$ is a locally finite union of quasi-compact closed subsets.
\end{cor}

\begin{prf}
 Let $S$ be the set of minimal points in $X$. For $\eta\in S$, we denote by $Z_\eta$ the closure of $\{\eta\}$.
 By Proposition \ref{prop:partially-proper} iii), $Z_\eta$ is quasi-compact.
 It is easy to observe that $\{Z_\eta\}$ cover $X$. Take a quasi-compact open subset $U$ of $X$.
 As $U$ is noetherian, it contains finitely many minimal points, thus $U$ intersects finitely many $Z_\eta$.
 This concludes that the closed covering $\{Z_\eta\}$ is locally finite.
\end{prf}

In the rest of this subsection, let $k$ be a field and $X$ a scheme which is partially proper over $k$.

\begin{defn}
 \begin{enumerate}
  \item For an (abelian \'etale) sheaf $\mathcal{F}$ on $X$, let $\Gamma_c(X,\mathcal{F})$ be
	the subset of $\Gamma(X,\mathcal{F})$ consisting of $s\in \Gamma(X,\mathcal{F})$ such that
	$\supp s$ is proper over $k$.
	As $\supp s$ is closed in $X$, this condition is equivalent to saying that
	$\supp s$ is quasi-compact (\cf Proposition \ref{prop:partially-proper} ii)).
  \item Let $H^i_c(X,-)$ be the $i$th derived functor of the left exact functor $\Gamma_c(X,-)$.
 \end{enumerate}
\end{defn}

\begin{prop}\label{prop:H_c-limit}
 Let $\mathcal{F}$ be a sheaf on $X$.
 \begin{enumerate}
  \item We have $H_c^i(X,\mathcal{F})\cong \varinjlim_Y H^i_Y(X,\mathcal{F})$, where $Y$ runs through
	quasi-compact closed subsets of $X$.
  \item We have $H_c^i(X,\mathcal{F})\cong \varinjlim_U H^i_c(U,\mathcal{F}\vert_U)$, where $U$ runs through
	quasi-compact open subsets of $X$.
 \end{enumerate}
\end{prop}

\begin{prf}
 i) If $i=0$, then the claim follows immediately from the definition of $\Gamma_c$.
 On the other hand, if $\mathcal{F}$ is injective, then $\varinjlim_YH^i_Y(X,\mathcal{F})=0$ for $i>0$.
 Therefore we have the desired isomorphism.

 ii) By Proposition \ref{prop:partially-proper} iii), for each quasi-compact open subset $U$ of $X$,
 we can find a quasi-compact closed subset $Y$ of $X$ containing $U$. On the other hand, for such a $Y$,
 we can find a quasi-compact open subset $U'$ of $X$ containing $Y$. Under this situation, we have
 push-forward maps
 \[
  H^i_c(U,\mathcal{F}\vert_U)\longrightarrow H^i_Y(X,\mathcal{F})\longrightarrow H^i_c(U',\mathcal{F}\vert_{U'}).
 \]
 These induce an isomorphism $\varinjlim_U H^i_c(U,\mathcal{F}\vert_U)\cong \varinjlim_Y H^i_Y(X,\mathcal{F})$.
 Hence the claim follows from i).
\end{prf}

\begin{cor}\label{cor:H_c-filtered-indlim}
 The functor $H^i_c(X,-)$ commutes with filtered inductive limits.
\end{cor}

\begin{prf}
 For a quasi-compact open subset $U$ of $X$, $H^i_c(U,-)$ commutes with filtered inductive limits;
 indeed, for a compactification $j\colon U\hooklongrightarrow \overline{U}$, we have
 $H^i_c(U,-)=H^i(\overline{U},j_!(-))$, and both $j_!$ and $H^i(\overline{U},-)$ commute with
 filtered inductive limits (for the later, note that $\overline{U}$ is quasi-compact and quasi-separated).
 On the other hand, the restriction functor to $U$ also commutes with filtered inductive limits.
 Hence Proposition \ref{prop:H_c-limit} ii) tells us that $H^i_c(X,-)$ commutes with filtered inductive limits.
\end{prf}

\begin{cor}\label{cor:H_c-direct-sum}
 For $i\ge 0$, the functor $H^i_c(X,-)$ commutes with arbitrary direct sums.
 Namely, for a set $\Lambda$ and sheaves $\{\mathcal{F}_\lambda\}_{\lambda\in\Lambda}$ on $X$
 indexed by $\Lambda$, we have an isomorphism 
 $H^i_c(X,\bigoplus_{\lambda\in\Lambda}\mathcal{F}_\lambda)\cong \bigoplus_{\lambda\in\Lambda}H^i_c(X,\mathcal{F}_\lambda)$.
\end{cor}

\begin{prf}
 For a finite subset $\Lambda_0$ of $\Lambda$, put $\mathcal{F}_{\Lambda_0}=\bigoplus_{\lambda\in\Lambda_0}\mathcal{F}_\lambda$. Then $\bigoplus_{\lambda\in\Lambda}\mathcal{F}_\lambda$ can be written as
 the filtered inductive limit $\varinjlim_{\Lambda_0\subset\Lambda}\mathcal{F}_{\Lambda_0}$.
 As $H^i_c(X,-)$ commutes with filtered inductive limits and finite direct sums,
 we obtain the desired result.
\end{prf}

\begin{rem}
 By exactly the same method as in \cite{MR1626021}, we can extend the definitions and properties above
 to $\ell$-adic coefficients.
\end{rem}

\subsection{Formal schemes and adic spaces}\label{subsec:formal-adic}
Let $R$ be a complete discrete valuation ring with separably closed residue field $\kappa$,
$F$ the fraction field of $R$,
and $k$ a separable closure of $F$. We denote by $k^+$ the valuation ring of $k$.
For a locally noetherian formal scheme $\mathcal{X}$ over $\Spf R$, we can associate an adic space
$t(\mathcal{X})$ over $\Spa(R,R)$ (\cf \cite[Proposition 4.1]{MR1306024}).
Its open subset $t(\mathcal{X})_\eta=t(\mathcal{X})\times_{\Spa(R,R)}\Spa(F,R)$ is
called the rigid generic fiber of $\mathcal{X}$.
In the following, we assume that $\mathcal{X}$ is special in the sense of Berkovich
\cite[\S 1]{MR1262943}. Then $t(\mathcal{X})_\eta$ is locally of finite type over $\Spa(F,R)$.
Therefore, we can make the fiber product
$t(\mathcal{X})_{\overline{\eta}}=t(\mathcal{X})_\eta\times_{\Spa(F,R)}\Spa(k,k^+)$,
which we call the rigid geometric generic fiber of $\mathcal{X}$.
The morphism $t(\mathcal{X})\longrightarrow \mathcal{X}$ of locally ringed spaces
induces a continuous map $t(\mathcal{X})_\eta\longrightarrow \mathcal{X}=\mathcal{X}^\red$.
We also have morphisms of sites 
$t(\mathcal{X})_{\overline{\eta},\et}\longrightarrow t(\mathcal{X})_{\eta,\et}\longrightarrow \mathcal{X}_\et\cong (\mathcal{X}^\red)_\et$.

\begin{prop}\label{prop:generic-fiber-partially-proper}
 Assume that $\mathcal{X}^\red$ is partially proper over $\kappa$.
 Then $t(\mathcal{X})_\eta$ is partially proper over $\Spa(F,R)$.
\end{prop}

\begin{prf}
 In \cite[Proposition 4.23]{formalnearby}, we have obtained the same result under the assumption
 that $\mathcal{X}^{\mathrm{red}}$ is proper.
 In fact, the proof therein only uses the partially properness of $\mathcal{X}^{\mathrm{red}}$.
\end{prf}

For simplicity, we put $X=t(\mathcal{X})_{\overline{\eta}}$.
We denote the composite 
$t(\mathcal{X})_{\overline{\eta}}\longrightarrow t(\mathcal{X})_\eta\longrightarrow \mathcal{X}^\red$ by $\spp$.
We also write $\spp$ for the morphism of \'etale sites $X_\et\longrightarrow (\mathcal{X}^\red)_\et$.
For a closed subset $Z$ of $\mathcal{X}$, consider the interior $\spp^{-1}(Z)^\circ$ of
$\spp^{-1}(Z)$ in $X$. It is called the tube of $Z$.

\begin{prop}\label{prop:tube-property}
 \begin{enumerate}
  \item Let $\mathcal{Z}$ be the formal completion of $\mathcal{X}$ along $Z$.
	Then the natural morphism $t(\mathcal{Z})_{\overline{\eta}}\longrightarrow t(\mathcal{X})_{\overline{\eta}}$ induces an isomorphism $t(\mathcal{Z})_{\overline{\eta}}\cong \spp^{-1}(Z)^\circ$.
  \item If $Z$ is quasi-compact, $\spp^{-1}(Z)^\circ$ is the union of countably many quasi-compact open subsets
	of $X$.
  \item Assume that $\mathcal{X}^\red$ is partially proper over $\kappa$.
	Then $\spp^{-1}(Z)^\circ$ is partially proper over $\Spa(k,k^+)$.
  \item Assume that $\mathcal{X}$ is locally algebraizable (\cf \cite[Definition 3.19]{formalnearby})
	and $Z$ is quasi-compact. Then, for a noetherian torsion ring $\Lambda$
	whose characteristic is invertible in $R$,
	$H^i(\spp^{-1}(Z)^\circ,\Lambda)$ and $H^i_c(\spp^{-1}(Z)^\circ,\Lambda)$ are 
	finitely generated $\Lambda$-modules.
 \end{enumerate}
\end{prop}

\begin{prf}
 i) is proved in \cite[Lemma 3.13 i)]{MR1620118}.
 For ii), we may assume that $\mathcal{X}$ is affine.
 Then, the claim has been obtained in the proof of \cite[Lemma 3.13 i)]{MR1620118}.
 By i), to prove iii) and iv), we may assume that $Z=\mathcal{X}^\red$.
 Then iii) follows from Proposition \ref{prop:generic-fiber-partially-proper}.
 As for iv), we have $H^i(X,\Lambda)=H^i(\mathcal{X}^\red,R\spp_*\Lambda)$.
 By \cite[Theorem 3.1]{MR1395723}, $R\spp_*\Lambda$ is a constructible complex on $\mathcal{X}^\red$.
 Thus $H^i(X,\Lambda)$ is a finitely generated $\Lambda$-module.
 On the other hand, by \cite[Proposition 3.21, Theorem 4.35]{formalnearby},
 $H^i_c(X,\Lambda)$ is a finitely generated $\Lambda$-module.
\end{prf}

\begin{defn}\label{defn:vertically-compact-supported}
 Assume that $\mathcal{X}^\red$ is partially proper over $\kappa$.
 We write $\Gamma_{c,\mathcal{X}}(X,-)$ for the composite functor $\Gamma_c(\mathcal{X}^\red,\spp_*(-))$.
 Denote the derived functor of $\Gamma_{c,\mathcal{X}}(X,-)$ by $R\Gamma_{c,\mathcal{X}}(X,-)$, and
 the $i$th cohomology of $R\Gamma_{c,\mathcal{X}}(X,-)$ by $H^i_{c,\mathcal{X}}(X,-)$.
\end{defn}

\begin{rem}\label{rem:p-adic-case}
 If $\mathcal{X}$ is a $p$-adic formal scheme, then for a quasi-compact closed subset $Z$
 of $\mathcal{X}^\red$, $\spp^{-1}(Z)$ is quasi-compact. Hence we have
 $\Gamma_{c,\mathcal{X}}(X,-)=\Gamma_c(X,-)$, $R\Gamma_{c,\mathcal{X}}(X,-)=R\Gamma_c(X,-)$ and
 $H^i_{c,\mathcal{X}}(X,-)=H^i_c(X,-)$ in this case. 
\end{rem}

The cohomology $H^i_{c,\mathcal{X}}(X,-)$ will appear in our main result.
We would like to discuss functoriality of this cohomology with respect to $X$.
For this purpose, we give another interpretation of $H^i_{c,\mathcal{X}}(X,-)$.
Here, more generally, let $X$ denote an adic space which is locally of finite type and partially proper over $\Spa(k,k^+)$.

\begin{defn}\label{defn:support-set}
 A support set $\mathcal{C}$ of $X$ is a set consisting of closed subsets of $X$
 satisfying the following conditions:
 \begin{itemize}
  \item For $Z,Z'\in \mathcal{C}$, we have $Z\cup Z'\in \mathcal{C}$.
  \item For $Z\in\mathcal{C}$ and a closed subset $Z'$ of $Z$, we have $Z'\in \mathcal{C}$.
 \end{itemize}
 For a support set $\mathcal{C}$, we define $\Gamma_{\mathcal{C}}(X,-)=\varinjlim_{Z\in\mathcal{C}}\Gamma_Z(X,-)$.
 Let $R\Gamma_{\mathcal{C}}(X,-)$ be the derived functor of $\Gamma_{\mathcal{C}}(X,-)$, and
 $H^i_{\mathcal{C}}(X,-)$ the $i$th cohomology of $R\Gamma_{\mathcal{C}}(X,-)$.
 It is easy to see that $H^i_{\mathcal{C}}(X,\mathcal{F})=\varinjlim_{Z\in\mathcal{C}}H^i_Z(X,\mathcal{F})$
 for a sheaf $\mathcal{F}$ on $X$.
\end{defn}

\begin{defn}\label{defn:support-set-mor}
 Let $X$ and $X'$ be adic spaces which are locally of finite type and partially proper over $\Spa(k,k^+)$.
 \begin{enumerate}
  \item For a support set $\mathcal{C}$ of $X$ and a morphism of adic spaces $f\colon X'\longrightarrow X$,
	let $f^{-1}\mathcal{C}$ be the support set of $X'$ consisting of
	closed subsets $Z'\subset X'$ which are contained in $f^{-1}(Z)$ for some $Z\in \mathcal{C}$.
  \item Let $\mathcal{C}$ and $\mathcal{C}'$ be support sets of $X$ and $X'$, respectively.
	A morphism of pairs $f\colon (X',\mathcal{C}')\longrightarrow (X,\mathcal{C})$
	is a morphism $f\colon X'\longrightarrow X$ satisfying $f^{-1}\mathcal{C}\subset\mathcal{C}'$.
	Such a morphism induces a morphism
	$f^*\colon R\Gamma_{\mathcal{C}}(X,\mathcal{F})\longrightarrow R\Gamma_{\mathcal{C}'}(X',f^*\mathcal{F})$.
	If moreover $f$ is an isomorphism and $f^{-1}\mathcal{C}=\mathcal{C}'$, $f$ is said to be
	an isomorphism of pairs.
 \end{enumerate}
\end{defn}

A formal model naturally gives a support set of the rigid generic fiber.

\begin{defn}\label{defn:formal-model-support-set}
 Let $\mathcal{X}$ be as in the beginning of this subsection, and assume that
 $\mathcal{X}^\red$ is partially proper over $\kappa$.
 We define a support set $\mathcal{C}_{\mathcal{X}}$ of $X=t(\mathcal{X})_{\overline{\eta}}$ as follows:
 a closed set of $X$ belongs to $\mathcal{C}_{\mathcal{X}}$ if it is contained in $\spp^{-1}(Z)$ for
 some quasi-compact closed subset of $\mathcal{X}^{\red}$.
\end{defn}

\begin{prop}\label{prop:formal-scheme-support}
 Let $\mathcal{X}$ be as in the previous definition.
 Then we have an isomorphism $R\Gamma_{c,\mathcal{X}}(X,-)\cong R\Gamma_{\mathcal{C_{\mathcal{X}}}}(X,-)$.
\end{prop}

\begin{prf}
 For a sheaf $\mathcal{F}$ on $X$, we have
 \[
  \Gamma_{c,\mathcal{X}}(X,\mathcal{F})=\varinjlim_{Z\subset\mathcal{X}^\red}\Gamma_Z(\mathcal{X}^\red,\spp_*\mathcal{F})=\varinjlim_{Z\subset\mathcal{X}^\red}\Gamma_{\spp^{-1}(Z)}(X,\mathcal{F})
 =\Gamma_{\mathcal{C}_{\mathcal{X}}}(X,\mathcal{F}),
 \]
 where $Z$ runs through quasi-compact closed subsets of $\mathcal{X}^\red$.
 This concludes the proof.
\end{prf}

\begin{lem}\label{lem:formal-finite-type}
 Let $\mathcal{X}$ and $\mathcal{X}'$ be special formal schemes over $\Spf R$ such that
 $\mathcal{X}^\red$ and $\mathcal{X}'^\red$ are partially proper over $\kappa$.
 Put $X=t(\mathcal{X})_{\overline{\eta}}$ and $X'=t(\mathcal{X}')_{\overline{\eta}}$, respectively.
 Let $f\colon \mathcal{X}'\longrightarrow \mathcal{X}$ be a morphism of finite type over $\Spf R$
 and write $f_{\overline{\eta}}$ for the induced morphism $X'\longrightarrow X$.
 Then, we have $f_{\overline{\eta}}^{-1}\mathcal{C}_{\mathcal{X}}=\mathcal{C}_{\mathcal{X}'}$.
 In particular, if moreover $f_{\overline{\eta}}$ is an isomorphism, then $f_{\overline{\eta}}$ induces
 an isomorphism of pairs $(X',\mathcal{C}_{\mathcal{X}'})\longrightarrow (X,\mathcal{C}_{\mathcal{X}})$.
\end{lem}

\begin{prf}
 For a quasi-compact closed subset $Z$ of $\mathcal{X}^\red$, $Z'=f^{-1}(Z)$ is a quasi-compact closed subset of 
 $\mathcal{X}'^\red$. Therefore $f_{\overline{\eta}}^{-1}(\spp^{-1}(Z))=\spp^{-1}(Z')$ lies in $\mathcal{C}_{\mathcal{X}'}$. This implies that $f^{-1}\mathcal{C}_{\mathcal{X}}$ is contained in $\mathcal{C}_{\mathcal{X}'}$.
 Conversely, let $Z'$ be a quasi-compact closed subset of $\mathcal{X}'^\red$. 
 As $\mathcal{X}^\red$ and $\mathcal{X}'^\red$ are partially proper over $\kappa$ and
 $f$ is of finite type, the induced morphism 
 $f\colon \mathcal{X}'^\red\longrightarrow \mathcal{X}^\red$ is proper.
 Therefore $Z=f(Z')$ is a quasi-compact closed subset of $\mathcal{X}^\red$.
 Hence, $\spp^{-1}(Z')$, being contained in $f_{\overline{\eta}}^{-1}(\spp^{-1}(Z))$,
 lies in $f^{-1}_{\overline{\eta}}\mathcal{C}_{\mathcal{X}}$. Thus we have
 $\mathcal{C}_{\mathcal{X}'}\subset f_{\overline{\eta}}^{-1}\mathcal{C}_{\mathcal{X}}$.
\end{prf}

\subsection{Smooth equivariant sheaves on adic spaces}
In \cite[\S IV.9]{MR2441311}, the theory of smooth equivariant sheaves on Berkovich spaces is developed.
In this subsection, we will adapt it to the framework of adic spaces.

\subsubsection{Basic definitions}
Let $k$ be a non-archimedean field and $X$ an adic space locally of finite type over $\Spa(k,k^+)$.
We denote by $X_\et$ the \'etale site of $X$, and $X_\qcet$ the full subcategory consisting of
\'etale morphisms $f\colon Y\longrightarrow X$ where $Y$ is quasi-compact.
The category $X_\qcet$ has a natural induced structure of a site and the associated topos
$\widetilde{X}_\qcet$ can be identified with the \'etale topos $\widetilde{X}_\et$ for $X$.

We write $p$ for the residue characteristic of $k^+$.
Fix a prime $\ell\neq p$ and a truncated discrete valuation ring $\Lambda$
with residue characteristic $\ell$.
Namely, $\Lambda$ can be written as $\mathcal{O}/(\lambda^n)$, where $\mathcal{O}$ is
a discrete valuation ring with residue characteristic $\ell$, $\lambda$ is a uniformizer of $\mathcal{O}$
and $n\ge 1$ is an integer. Such a ring has the following property:

\begin{lem}\label{lem:TDVR-injective}
 For a $\Lambda$-module $M$, the following are equivalent:
 \begin{enumerate}
  \item $M$ is an injective $\Lambda$-module.
  \item $M$ is a flat $\Lambda$-module.
  \item $M$ is a free $\Lambda$-module.
 \end{enumerate}
\end{lem}

\begin{prf}
 We write $\Lambda=\mathcal{O}/(\lambda^n)$ as above.
 For every $0\le k\le n$, we have an exact sequence
 $\Lambda\yrightarrow{\times\lambda^{n-k}}\Lambda\yrightarrow{\times\lambda^{k}}\Lambda\longrightarrow \Lambda/(\lambda^k)\longrightarrow 0$. 
 Thus we obtain
 \[
  \Tor_1^{\Lambda}\bigl(\Lambda/(\lambda^k),M\bigr)=\Ker(M\yrightarrow{\times\lambda^k}M)/\Imm(M\yrightarrow{\times\lambda^{n-k}}M)
 =\Ext^1_\Lambda\bigl(\Lambda/(\lambda^{n-k}),M\bigr).
 \]
 Therefore, if $M$ is injective, then $\Tor_1^\Lambda(\Lambda/(\lambda),M)=0$ and $M$ is flat.
 Conversely, if $M$ is flat, then $\Ext^1_\Lambda(\Lambda/(\lambda^k),M)=0$ for every $0\le k\le n$.
 Hence Baer's criterion on injectivity tells us that $M$ is an injective $\Lambda$-module.

 On the other hand, as $\Lambda$ is an Artinian local ring,
 ii) and iii) are equivalent.
\end{prf}
We write $\Lambda\text{-}\widetilde{X}_\et$ for the category of $\Lambda$-sheaves on $X_\et$.

Let $G$ be a locally pro-$p$ group, and assume that $X$ is equipped with
a continuous action of $G$ in the following sense. 

\begin{defn}\label{defn:action-continuous}
 An action of $G$ on $X$ is said to be continuous if the following conditions are satisfied:
 \begin{quote}
  For every affinoid open subspace $U=\Spa(A,A^+)$ of $X$, $f\in A$ and
  $a\in k^\times$, there exists an open subgroup $G'\subset G$ such that
  each element $g\in G'$ satisfies $gU=U$ and $\lvert (g^*f-f)(x)\rvert\le \lvert a(x)\rvert$ for $x\in U$.
 \end{quote}
\end{defn}
Under the continuity condition, we know the following result due to Berkovich:

\begin{thm}[{{\cite[Key Lemma 7.2]{MR1262943}}}]
 For every object $f\colon Y\longrightarrow X$ of $X_\qcet$, there exists a compact open subgroup $K_Y$
 of $G$ such that the action of $K_Y$ on $X$ lifts canonically and continuously to $Y$.
\end{thm}

\begin{defn}
 A $G$-equivariant sheaf $\mathcal{F}$ on $X_\et$ is said to be smooth if
 for every object $Y\longrightarrow X$ of $X_\qcet$ the action of $K_Y$ on $\Gamma(Y,\mathcal{F})$
 is smooth.
 We write $\widetilde{X}_\et/G$ for the category of smooth $G$-equivariant sheaves on $X_\et$,
 and $\Lambda\text{-}\widetilde{X}_\et/G$ for the category of smooth $G$-equivariant $\Lambda$-sheaves
 on $X_\et$.
\end{defn}

As $\widetilde{X}_\et=\widetilde{X}_\qcet$, the definition above is a special case of
\cite[D\'efinition IV.8.1]{MR2441311}. In particular, $\widetilde{X}_\et/G$ is a topos
(\cf \cite[Proposition IV.8.12]{MR2441311}).
The forgetful functor $\Lambda\text{-}\widetilde{X}_\et/G\longrightarrow \Lambda\text{-}\widetilde{X}_\et$
is exact (\cf \cite[Corollaire IV.8.6]{MR2441311}).

\begin{defn}
 Let $H$ be a closed subgroup of $G$. 
 \begin{enumerate}
  \item A smooth $G$-equivariant $\Lambda$-sheaf on $X_\et$ can be obviously regarded as a smooth $H$-equivariant 
	$\Lambda$-sheaf on $X_\et$.
	Therefore we get a functor $\Lambda\text{-}\widetilde{X}_\et/G\longrightarrow \Lambda\text{-}\widetilde{X}_\et/H$,
	which is denoted by $\Res_{X/H}^{X/G}$.
  \item Let $\Omega_H$ be a system of representatives of $H\backslash G$.
	For a smooth $H$-equivariant $\Lambda$-sheaf $\mathcal{F}$ on $X_\et$, the $\Lambda$-sheaf
	$\prod_{g\in \Omega_H}g^*\mathcal{F}$ has a natural $G$-equivariant structure.
	We put $\Ind_{X/H}^{X/G}\mathcal{F}=(\prod_{g\in \Omega_H}g^*\mathcal{F})^\infty$,
	where $(-)^\infty$ denotes the smoothification functor
	(\cf \cite[\S IV.8.3.3]{MR2441311}).
	This gives a functor $\Ind_{X/H}^{X/G}\colon \Lambda\text{-}\widetilde{X}_\et/H\longrightarrow \Lambda\text{-}\widetilde{X}_\et/G$.
  \item If $H$ is an open subgroup of $G$, for $\mathcal{F}\in\Lambda\text{-}\widetilde{X}_\et/H$ we put
	\[
	 \cInd_{X/H}^{X/G}\mathcal{F}=\bigoplus_{g\in \Omega_H}g^*\mathcal{F}.
	\]
	It can be naturally regarded as an object of $\Lambda\text{-}\widetilde{X}_\et/G$.
	This gives a functor $\cInd_{X/H}^{X/G}\colon \Lambda\text{-}\widetilde{X}_\et/H\longrightarrow \Lambda\text{-}\widetilde{X}_\et/G$.
 \end{enumerate}
\end{defn}

The following proposition is obvious.

\begin{prop}\label{prop:Ind-adjoint}
 The functor $\Ind_{X/H}^{X/G}$ is the right adjoint of $\Res_{X/H}^{X/G}$, and
 $\cInd_{X/H}^{X/G}$ is the left adjoint of $\Res_{X/H}^{X/G}$.
 The functors $\Res_{X/H}^{X/G}$ and $\cInd_{X/H}^{X/G}$ are exact, and $\Ind_{X/H}^{X/G}$ is left exact.
\end{prop}

Let $f\colon U\longrightarrow X$ be an \'etale morphism. Assume that there exists a compact open subgroup
$K$ of $G$ whose action on $X$ lifts to $U$ continuously (if $U$ is quasi-compact, this is always the case).
We fix such a compact open subgroup $K$.

\begin{lem}
 For $\mathcal{F}\in\Lambda\text{-}\widetilde{U}_\et/K$, 
 the induced $K$-equivariant structure on $f_!\mathcal{F}$ is smooth.
\end{lem}

\begin{prf}
 Let $\mathcal{G}$ be the presheaf on $X_\qcet$ defined as follows:
 for an object $Y\longrightarrow X$ of $X_\qcet$, we put
 \[
  \Gamma(Y,\mathcal{G})=\bigoplus_{\phi\in \Hom_X(Y,U)}\Gamma(Y\yrightarrow{\phi}U,\mathcal{F}).
 \]
 We will prove that the action of $K\cap K_Y$ on $\Gamma(Y,\mathcal{G})$ is smooth.
 Applying \cite[Key Lemma 7.2]{MR1262943} to
 $\phi\in\Hom_X(Y,U)$, we obtain a compact open subgroup $K_\phi$ of $K\cap K_Y$ such that
 $g\in K_\phi$ satisfies $g\circ \phi=\phi\circ g$. 
 Such $K_\phi$ acts on $\Gamma(Y\yrightarrow{\phi}U,\mathcal{F})$, and 
 by the smoothness of $\mathcal{F}$, this action is smooth.
 This concludes the smoothness of the action of $K\cap K_Y$ on $\Gamma(Y,\mathcal{G})$.

 In other words, $\mathcal{G}$ is a smooth $K$-equivariant $\Lambda$-presheaf on $X_{\qcet}$ 
 (\cf \cite[D\'efinition IV.8.1]{MR2441311}).
 Hence \cite[Lemme IV.8.4]{MR2441311} tells us that the sheafification $f_!\mathcal{F}$ of $\mathcal{G}$
 is smooth.
\end{prf}

\begin{defn}
 \begin{enumerate}
  \item Let $\Res_{U/K}^{X/K}\colon \Lambda\text{-}\widetilde{X}_\et/K\longrightarrow \Lambda\text{-}\widetilde{U}_\et/K$ be the functor $\mathcal{F}\longmapsto f^*\mathcal{F}$. In fact, it is easy to see that the induced $K$-equivariant
	structure on $f^*\mathcal{F}$ is smooth. We put $\Res_{U/K}^{X/G}=\Res_{U/K}^{X/K}\circ \Res_{X/K}^{X/G}$,
	which is a functor from $\Lambda\text{-}\widetilde{X}_\et/G$ to $\Lambda\text{-}\widetilde{U}_\et/K$.
  \item Let $\Ind_{U/K}^{X/K}\colon \Lambda\text{-}\widetilde{U}_\et/K\longrightarrow \Lambda\text{-}\widetilde{X}_\et/K$ be the functor $\mathcal{F}\longmapsto (f_*\mathcal{F})^\infty$;
	note that $f_*\mathcal{F}$ carries a $K$-equivariant structure,
	but it is not necessarily smooth.
	We put $\Ind_{U/K}^{X/G}=\Ind_{X/K}^{X/G}\circ \Ind_{U/K}^{X/K}$,
	which is a functor from $\Lambda\text{-}\widetilde{U}_\et/K$ to $\Lambda\text{-}\widetilde{X}_\et/G$.
  \item Let $\cInd_{U/K}^{X/K}\colon \Lambda\text{-}\widetilde{U}_\et/K\longrightarrow \Lambda\text{-}\widetilde{X}_\et/K$ be the functor $\mathcal{F}\longmapsto f_!\mathcal{F}$;
	(\cf the lemma above).
	We put $\cInd_{U/K}^{X/G}=\cInd_{X/K}^{X/G}\circ \cInd_{U/K}^{X/K}$,
	which is a functor from $\Lambda\text{-}\widetilde{U}_\et/K$ to $\Lambda\text{-}\widetilde{X}_\et/G$.
 \end{enumerate}
\end{defn}

The following proposition is also immediate:

\begin{prop}\label{prop:Ind-adjoint-open}
 The functor $\Ind_{U/K}^{X/G}$ is the right adjoint of $\Res_{U/K}^{X/G}$, and
 $\cInd_{U/K}^{X/G}$ is the left adjoint of $\Res_{U/K}^{X/G}$.
 The functors $\Res_{U/K}^{X/G}$ and $\cInd_{U/K}^{X/G}$ are exact, and $\Ind_{U/K}^{X/G}$ is left exact.
 In particular, $\Res_{U/K}^{X/G}$ and $\Ind_{U/K}^{X/G}$ preserve injective objects.
\end{prop}

Let $U$ be an open subset of $X$ which is stable under a compact open subgroup $K$ of $G$.
If $\{gU\}_{g\in G}$ covers $X$, we can construct canonical resolutions of a smooth $G$-equivariant $\Lambda$-sheaf
by using the functors $\Res_{U/K}^{X/G}$, $\Ind_{U/K}^{X/G}$ and $\cInd_{U/K}^{X/G}$
(\cf \cite[Th\'eor\`eme IV.9.31]{MR2441311}).

\begin{prop}\label{prop:functorial-resolution}
 Let $U$ be an open subset of $X$ which is stable under a compact open subgroup $K$ of $G$.
 Assume that $X=\bigcup_{g\in G}gU$. For $\alpha=(\overline{g}_1,\ldots,\overline{g}_s)\in (G/K)^s$,
 we put $U_\alpha=g_1U\cap\cdots\cap g_sU$ and $K_\alpha=g_1Kg_1^{-1}\cap\cdots\cap g_sKg_s^{-1}$.
 Let $\mathcal{F}$ be an object of $\Lambda\text{-}\widetilde{X}_\et/G$.
 \begin{enumerate}
  \item We have a functorial resolution $C^\bullet(\mathcal{F})\longrightarrow \mathcal{F}$ in 
	$\Lambda\text{-}\widetilde{X}_\et/G$ where $C^m(\mathcal{F})$ is given by
	\[
	 C^m(\mathcal{F})=\bigoplus_{\alpha\in G\backslash (G/K)^{-m+1}}\cInd_{U_\alpha/K_\alpha}^{X/G}\Res_{U_\alpha/K_\alpha}^{X/G}\mathcal{F}.
	\]
	Here we abuse notation: $U_\alpha$ and $K_\alpha$ depends on the choice of a lift of
	$\alpha\in G\backslash (G/K)^{-m+1}$ to $(G/K)^{-m+1}$, but 
	$\cInd_{U_\alpha/K_\alpha}^{X/G}\Res_{U_\alpha/K_\alpha}^{X/G}\mathcal{F}$ does not.
  \item Assume moreover that the covering $\{gU\}_{g\in G/K}$ is locally finite; namely,
	each point $x\in X$ has an open neighborhood which intersects only finitely many of $gU$ with $g\in G/K$.
	Then we have a functorial resolution $\mathcal{F}\longrightarrow D^\bullet(\mathcal{F})$ in 
	$\Lambda\text{-}\widetilde{X}_\et/G$ where $D^m(\mathcal{F})$ is given by
	\[
	D^m(\mathcal{F})=\bigoplus_{\alpha\in G\backslash (G/K)^{m+1}}\Ind_{U_\alpha/K_\alpha}^{X/G}\Res_{U_\alpha/K_\alpha}^{X/G}\mathcal{F}.
	\]
 \end{enumerate}
\end{prop}

\begin{prf}
 For $\lambda\in (G/K)^s$, we write $j_\lambda$ for the natural open immersion $U_\lambda\hooklongrightarrow X$.
 
 i)  We have a well-known exact sequence
 \[
 \cdots\longrightarrow \bigoplus_{\lambda\in (G/K)^2}j_{\lambda!}j_\lambda^*\mathcal{F}\longrightarrow \bigoplus_{\lambda\in G/K}j_{\lambda!}j_\lambda^*\mathcal{F}\longrightarrow \mathcal{F}\longrightarrow 0
 \]
 of $\Lambda$-sheaves over $X_\et$.
 It is easy to see that $\bigoplus_{\lambda\in (G/K)^{-m+1}}j_{\lambda!}j_\lambda^*\mathcal{F}$ coincides with
 the underlying $\Lambda$-sheaf of $\bigoplus_{\alpha\in G\backslash (G/K)^{-m+1}}\cInd_{U_\alpha/K_\alpha}^{X/G}\Res_{U_\alpha/K_\alpha}^{X/G}\mathcal{F}$,
 and each homomorphism in the complex above is $G$-equivariant. Hence we have a resolution
 \[
  \bigoplus_{\alpha\in G\backslash (G/K)^{-\bullet+1}}\cInd_{U_\alpha/K_\alpha}^{X/G}\Res_{U_\alpha/K_\alpha}^{X/G}\mathcal{F}\longrightarrow\mathcal{F}\longrightarrow 0.
 \]

 ii) Consider the complex
 \begin{equation*}
 0\longrightarrow \mathcal{F}\longrightarrow \prod_{\lambda\in G/K}j_{\lambda*}j_\lambda^*\mathcal{F}\longrightarrow \prod_{\lambda\in (G/K)^2}j_{\lambda*}j_\lambda^*\mathcal{F}\longrightarrow \cdots\tag{$*$}
 \end{equation*}
 of $\Lambda$-sheaves over $X_\et$. Each term has a $G$-equivariant structure and each homomorphism 
 is $G$-equivariant. 
 For each $g\in G$ and $m\ge 0$, the map
 \[
 \Bigl(\prod_{\lambda\in (G/K)^{m+1}}j_{\lambda*}j_\lambda^*\mathcal{F}\Bigr)\Bigr\rvert_{gU}
 \longrightarrow \Bigl(\prod_{\lambda\in (G/K)^m}j_{\lambda*}j_\lambda^*\mathcal{F}\Bigr)\Bigr\rvert_{gU};\
 (s_\lambda)_{\lambda\in (G/K)^{m+1}}\longmapsto (s_{(g,\lambda)})_{\lambda\in (G/K)^m},
 \]
 where $s_\lambda$ is a local section of $(j_{\lambda*}j_\lambda^*\mathcal{F})\vert_{gU}$,
 is $gKg^{-1}$-equivariant and gives a homotopy between $\id$ and $0$ on the complex $(*)$
 restricted on $gU$. In particular, the smoothification
 \[
   0\longrightarrow \mathcal{F}\longrightarrow \Bigl(\prod_{\lambda\in G/K}j_{\lambda*}j_\lambda^*\mathcal{F}\Bigr)^\infty\longrightarrow \Bigl(\prod_{\lambda\in (G/K)^2}j_{\lambda*}j_\lambda^*\mathcal{F}\Bigr)^\infty\longrightarrow \cdots
 \]
 of the complex $(*)$ is exact.

 Now, by the assumption on the covering $\{gU\}_{g\in G/K}$, we have
 \[
 \prod_{\lambda\in (G/K)^{m+1}}j_{\lambda*}j_\lambda^*\mathcal{F}
 =\bigoplus_{\alpha\in G\backslash(G/K)^{m+1}}\prod_{g\in G}g^*j_{\alpha*}j_\alpha^*\mathcal{F}.
 \]
 Therefore the smoothification of the $G$-equivariant sheaf
 $\prod_{\lambda\in (G/K)^{m+1}}j_{\lambda*}j_\lambda^*\mathcal{F}$ coincides with
 $\bigoplus_{\alpha\in G\backslash(G/K)^{m+1}}\Ind_{U_\alpha/K_\alpha}^{X/G}\Res_{U_\alpha/K_\alpha}^{X/G}\mathcal{F}$. Hence we have a resolution
 \[
 0\longrightarrow \mathcal{F}\longrightarrow \bigoplus_{\alpha\in G\backslash(G/K)^{\bullet+1}}\Ind_{U_\alpha/K_\alpha}^{X/G}\Res_{U_\alpha/K_\alpha}^{X/G}\mathcal{F}.
 \]
\end{prf}

\begin{rem}
 In this paper, we only use the part i) of the proposition above.
 Later we will consider a variant of ii) (see Lemma \ref{lem:non-smooth-functorial-resolution}).
\end{rem}

\subsubsection{Acyclicity}
In the following we will give several acyclicity results for injective objects in
$\Lambda\text{-}\widetilde{X}_\et/G$.

\begin{prop}\label{prop:injective-acyclic}
 For an injective object $\mathcal{F}$ in $\Lambda\text{-}\widetilde{X}_\et/G$
 and an object $Y\longrightarrow X$ of $X_{\et}$ such that $Y$ is quasi-compact and quasi-separated,
 we have $H^i(Y,\mathcal{F})=0$ for $i\ge 1$.
\end{prop}

\begin{prf}
 Let $X_\cohet$ be the full subcategory of $X_\et$ consisting of
 \'etale morphisms $Y\longrightarrow X$ where $Y$ is quasi-compact and quasi-separated.
 It is naturally equipped with a structure of a site, 
 and the natural morphism of sites $X_\et\longrightarrow X_\cohet$ induces an isomorphism of toposes
 $\widetilde{X}_\et\cong \widetilde{X}_\cohet$.

 Let $\mathcal{F}$ be a $G$-equivariant sheaf on $X_\et$ such that for every object $Y\longrightarrow X$
 in $X_\cohet$ the action of $K_Y\subset G$ on $\Gamma(Y,\mathcal{F})$ is smooth.
 Then $\mathcal{F}$ is a smooth $G$-equivariant sheaf; note that every object $Y\longrightarrow X$ in $X_\qcet$
 can be covered by finitely many objects $(Y_\alpha\longrightarrow X)_\alpha$ in $X_\cohet$, and then
 $\Gamma(Y,\mathcal{F})\longrightarrow \prod_{\alpha}\Gamma(Y_\alpha,\mathcal{F})$ is injective.
 Therefore we have an isomorphism of toposes
 \[
  \widetilde{X}_\et/G=\widetilde{X}_\qcet/G\cong \widetilde{X}_\cohet/G
 \]
 (for the definition of $\widetilde{X}_\cohet/G$, see \cite[D\'efinition IV.8.2]{MR2441311}).
 
 We prove that fiber products exist in the category $X_\cohet$. Let $Y\longrightarrow X$ and
 $Z_i\longrightarrow X$ $(i=1,2)$ be objects in $X_\cohet$, and assume that we are given
 morphisms $Z_1\longrightarrow Y$ and $Z_2\longrightarrow Y$ over $X$.
 Then, these morphisms are quasi-compact quasi-separated,
 and so is $Z_1\times_YZ_2\longrightarrow Z_2$. (In the case of schemes, see \cite[IV, \S 1.1, \S 1.2]{EGA}.
 The arguments there can be applied to adic spaces.)
 Hence $Z_1\times_YZ_2$ is quasi-compact and quasi-separated, as desired.

 Now we can apply \cite[Th\'eor\`eme IV.8.15]{MR2441311} (or \cite[Th\'eor\`eme IV.8.17]{MR2441311})
 to conclude the proposition.
\end{prf}

\begin{prop}\label{prop:restr-surjective}
 Let $Y\longrightarrow X$ be an object of $X_\qcet$, and $U$ a quasi-compact open subset of $Y$.
 For an injective object $\mathcal{F}$ in $\Lambda\text{-}\widetilde{X}_\et/G$,
 the homomorphism $\Gamma(Y,\mathcal{F})\longrightarrow \Gamma(U,\mathcal{F})$ is surjective.
\end{prop}

\begin{prf}
 Let $K$ be a compact open subgroup of $K_Y$ which stabilizes $U$.
 We write $\Lambda_Y$ (resp.\ $\Lambda_U$) for the constant sheaf on $Y$ (resp.\ $U$) with values in $\Lambda$.
 They can be regarded as a smooth $K$-equivariant sheaves by the trivial $K$-actions.
 Therefore we can form $\cInd_{Y/K}^{X/G}\Lambda_Y$ and $\cInd_{U/K}^{X/G}\Lambda_U$.
 By Proposition \ref{prop:Ind-adjoint-open}, we have
 \[
  \Hom(\cInd_{Y/K}^{X/G}\Lambda_Y,\mathcal{F})=\Gamma(Y,\mathcal{F})^K,\quad 
  \Hom(\cInd_{U/K}^{X/G}\Lambda_U,\mathcal{F})=\Gamma(U,\mathcal{F})^K.
 \]
 If we denote by $j$ the natural open immersion $U\hooklongrightarrow Y$, we have an injection
 $j_!\Lambda_U\longrightarrow \Lambda_Y$. This gives an injection
 $\cInd_{U/K}^{Y/K}\Lambda_U\longrightarrow \Lambda_Y$ in $\Lambda\text{-}\widetilde{Y}_\et/K$,
 and thus an injection
 $\cInd_{U/K}^{X/G}\Lambda_U\longrightarrow \cInd_{Y/K}^{X/G}\Lambda_Y$ in $\Lambda\text{-}\widetilde{X}_\et/G$.
 Since $\mathcal{F}$ is an injective object of $\Lambda\text{-}\widetilde{X}_\et/G$, the induced homomorphism
 \[
  \Hom(\cInd_{Y/K}^{X/G}\Lambda_Y,\mathcal{F})\longrightarrow \Hom(\cInd_{U/K}^{X/G}\Lambda_U,\mathcal{F})
 \]
 is surjective. Therefore the map $\Gamma(Y,\mathcal{F})^K\longrightarrow \Gamma(U,\mathcal{F})^K$ is surjective.

 As $\mathcal{F}$ is smooth, we have $\Gamma(U,\mathcal{F})=\varinjlim_{K}\Gamma(U,\mathcal{F})^K$.
 Hence the homomorphism $\Gamma(Y,\mathcal{F})\longrightarrow \Gamma(U,\mathcal{F})$ is also surjective.
\end{prf}

\begin{cor}\label{cor:local-coh-vanish}
 Let $U$ and $V$ be quasi-compact quasi-separated open subsets of $X$
 such that $V\subset U$.
 For an injective object $\mathcal{F}$ in $\Lambda\text{-}\widetilde{X}_\et/G$, we have
 $H^i_{U\setminus V}(U,\mathcal{F})=0$ for $i\ge 1$.
\end{cor}

\begin{prf}
 By Proposition \ref{prop:injective-acyclic} and Proposition \ref{prop:restr-surjective},
 we have $H^i_{U\setminus V}(U,\mathcal{F})=0$ for $i\ge 2$,
 and $H^1_{U\setminus V}(U,\mathcal{F})=\Coker(\Gamma(U,\mathcal{F})\to\Gamma(V,\mathcal{F}))=0$.
\end{prf}

\begin{prop}\label{prop:countable-union}
 Let $Y\longrightarrow X$ be an \'etale morphism. Assume that $Y$ is quasi-separated, and
 is the union of countably many quasi-compact open subsets of $Y$. 
 Then, for an injective object $\mathcal{F}$ in $\Lambda\text{-}\widetilde{X}_\et/G$, we have
 $H^i(Y,\mathcal{F})=0$ for $i\ge 1$.
\end{prop}

\begin{prf}
 By the assumption on $Y$, there exists an increasing series $U_1\subset U_2\subset\cdots$ of quasi-compact
 open subsets of $Y$ such that $Y=\bigcup_{n=1}^\infty U_n$.
 For a $\Lambda$-sheaf $\mathcal{G}$ on $X_\et$, we will construct the following spectral sequence:
 \[
  E_2^{i,j}=\varprojlim_n\!{}^i H^j(U_n,\mathcal{G})\longrightarrow H^{i+j}(Y,\mathcal{G}).
 \]
 To show the existence of this spectral sequence, it suffices to show that for an injective object $\mathcal{I}$
 in $\Lambda\text{-}\widetilde{X}_\et$, the projective system $(\Gamma(U_n,\mathcal{I}))_n$ is
 an injective object in the category of projective systems of $\Lambda$-modules.
 By \cite[Proposition 1.1]{MR929536}, we should prove the following two properties:
 \begin{enumerate}
  \item[(a)] $\Gamma(U_n,\mathcal{I})$ is an injective $\Lambda$-module for every $n\ge 1$.
  \item[(b)] The transition map $\Gamma(U_{n+1},\mathcal{I})\longrightarrow \Gamma(U_n,\mathcal{I})$
	     is a split surjection for every $n$.
 \end{enumerate}
 (a) is easy. For (b), note the following exact sequence:
 \[
  0\longrightarrow \Gamma_{U_{n+1}\setminus U_n}(U_{n+1},\mathcal{I})\longrightarrow 
 \Gamma(U_{n+1},\mathcal{I})\longrightarrow \Gamma(U_n,\mathcal{I})\longrightarrow 0.
 \]
 As $\Gamma_{U_{n+1}\setminus U_n}(U_{n+1},\mathcal{I})$ is an injective $\Lambda$-module,
 the map $\Gamma(U_{n+1},\mathcal{I})\longrightarrow \Gamma(U_n,\mathcal{I})$ is
 a split surjection, as desired. 

 As $\varinjlim_n^i=0$ for $i\ge 2$, we obtain the following exact sequence:
 \[
  0\longrightarrow \varprojlim_n\!{}^1 H^{i-1}(U_n,\mathcal{G})\longrightarrow H^i(Y,\mathcal{G})
 \longrightarrow \varprojlim_n H^i(U_n,\mathcal{G})\longrightarrow 0.
 \]
 If $i\ge 1$, Proposition \ref{prop:injective-acyclic} tells us that $H^i(U_n,\mathcal{F})=0$ for every $n$.
 Therefore we have $H^i(Y,\mathcal{F})=0$ for $i\ge 2$, 
 and $H^1(Y,\mathcal{F})=\varprojlim_n^1 \Gamma(U_n,\mathcal{F})$.
 By Proposition \ref{prop:restr-surjective}, the map
 $\Gamma(U_{n+1},\mathcal{F})\longrightarrow \Gamma(U_n,\mathcal{F})$ is a surjection.
 Hence we have $H^1(Y,\mathcal{F})=\varprojlim_n^1 \Gamma(U_n,\mathcal{F})=0$.
 This completes the proof.
\end{prf}

\begin{cor}\label{cor:j_*-acyclic}
 Assume that $X$ is quasi-separated. Let $U$ be an open subset of $X$ which is the union of countably many
 quasi-compact open subsets of $X$. 
 We write $j$ for the natural open immersion $U\hooklongrightarrow X$.
 Then, for an injective object $\mathcal{F}$ in $\Lambda\text{-}\widetilde{X}_\et/G$, we have
 $R^ij_*j^*\mathcal{F}=0$ for $i\ge 1$.
\end{cor}

\begin{prf}
 Note that $R^ij_*j^*\mathcal{F}$ is the sheafification of the presheaf
 \[
  (Y\longrightarrow X)\longmapsto H^i(Y\times_XU,\mathcal{F})
 \]
 on $X_\et$. Thus it suffices to show that $H^i(Y\times_XU,\mathcal{F})=0$
 for $i\ge 1$ and an object $Y\longrightarrow X$ of $X_\cohet$
  (\cf the proof of Proposition \ref{prop:injective-acyclic}).
 Take an increasing series $U_1\subset U_2\subset\cdots$ of quasi-compact
 open subsets of $U$ such that $U=\bigcup_{n=1}^\infty U_n$.
 Then, as $X$ is quasi-separated, $Y\times_XU_n$ is a quasi-compact open subset of $Y$.
 Therefore,  $Y\times_XU$ is the union of countably many quasi-compact open subsets $(Y\times_XU_n)_{n\ge 1}$.
 Hence Proposition \ref{prop:countable-union} tells us that $H^i(Y\times_XU,\mathcal{F})=0$ for $i\ge 1$.
 This concludes the proof.
\end{prf}

Recall the result in \cite[Lemme 2.3]{Fargues-Zel}:

\begin{lem}\label{lem:section-D_c-action}
 Let $U$ be an open subset of $X$ which is stable under a compact open subgroup $K$ of $G$.
 For a smooth $G$-equivariant $\Lambda$-sheaf $\mathcal{F}$, $\Gamma(U,\mathcal{F})$ has a natural structure
 of a $\mathcal{D}_c(K)$-module.
\end{lem}

\begin{prf}
 We have $\Gamma(U,\mathcal{F})=\varprojlim_{V\subset U}\Gamma(V,\mathcal{F})$, where $V$ runs through
 quasi-compact open subsets of $U$ which are stable under $K$. As the action of $K$ on $\Gamma(V,\mathcal{F})$
 is smooth, $\Gamma(V,\mathcal{F})$ has a structure of a $\mathcal{D}_c(K)$-module, and the transition maps
 of the projective system $(\Gamma(V,\mathcal{F}))_V$ are compatible with the actions of $\mathcal{D}_c(K)$.
 Hence $\Gamma(U,\mathcal{F})$ has a structure of a $\mathcal{D}_c(K)$-module.
\end{prf}

\begin{defn}
 Let $U$ be an open subset of $X$ which is stable under a compact open subgroup $K$ of $G$.
 By the lemma above, a left exact functor
 \[
 \Gamma(U/K,-)\colon \Lambda\text{-}\widetilde{X}_\et/G\longrightarrow \Mod\bigl(\mathcal{D}_c(K)\bigr);
 \quad\mathcal{F}\longmapsto\Gamma(U,\mathcal{F})
 \]
 is induced. We denote by $R\Gamma(U/K,-)$ the right derived functor of $\Gamma(U/K,-)$,
 and by $H^i(U/K,-)$ the $i$th cohomology of $R\Gamma(U/K,-)$.
\end{defn}

\begin{rem}\label{rem:quasi-compact-smooth}
 As the functors
 \[
  i_\mathcal{D}\colon \Rep_\Lambda(K)\longrightarrow \Mod\bigl(\mathcal{D}_c(K)\bigr),\quad
 \infty_\mathcal{D}\colon \Mod\bigl(\mathcal{D}_c(K)\bigr)\longrightarrow \Rep_\Lambda(K)
 \]
 are exact, they induce functors between derived categories
 \begin{align*}
  i_\mathcal{D}&\colon D^+\bigl(\Rep_\Lambda(K)\bigr)\longrightarrow D^+\bigl(\Mod(\mathcal{D}_c(K))\bigr),\\
 \infty_\mathcal{D}&\colon D^+\bigl(\Mod(\mathcal{D}_c(K))\bigr)\longrightarrow D^+\bigl(\Rep_\Lambda(K)\bigr).
 \end{align*}
 These two functors are adjoint to each other, and satisfy $\infty_\mathcal{D}\circ i_{\mathcal{D}}=\id$.
 Therefore $D^+(\Rep_\Lambda(K))$ can be regarded as a full subcategory of 
 $D^+(\Mod(\mathcal{D}_c(K)))$ by $i_\mathcal{D}$.

 Under this setting, if $U$ in the previous definition is quasi-compact, the image of $R\Gamma(U/K,-)$
 lies in $D^+(\Rep_\Lambda(K))$.
\end{rem}

\begin{cor}\label{cor:derived-same-Gamma}
 Let $U$ and $K$ be as in the definition above. Assume that $U$ is quasi-separated,
 and is the union of countably many quasi-compact open subsets.
 Then the following diagram is 2-commutative:
 \[
 \xymatrix{%
 D^+(\Lambda\text{-}\widetilde{X}_\et/G)\ar[rr]^-{R\Gamma(U/K,-)}\ar[d]&& D^+\bigl(\Mod(\mathcal{D}_c(K))\bigr)\ar[d]\\
 D^+(\Lambda\text{-}\widetilde{X}_\et)\ar[rr]^-{R\Gamma(U,-)}&& D^+\bigl(\Mod(\Lambda)\bigr)\lefteqn{.}
 }
 \]
 Here $\Mod(\Lambda)$ denotes the category of $\Lambda$-modules, and the vertical arrows denote
 the forgetful functors.
\end{cor}

\begin{prf}
 Clear from Proposition \ref{prop:countable-union}.
\end{prf}

\begin{prop}\label{prop:constructible-smooth}
 Let $U$ be an open subset of $X$ which is stable under a compact open subgroup $K$ of $G$.
 Let $\mathcal{F}$ be an object of $\Lambda\text{-}\widetilde{X}_\et/G$ and $i\ge 0$ an integer.
 Assume the following conditions:
 \begin{itemize}
  \item the base field $k$ is separably closed.
  \item the characteristic of $k$ is zero, or $X$ is smooth over $\Spa(k,k^+)$.
  \item $U$ is quasi-separated, and is the union of countably many quasi-compact open subsets of $X$.
  \item $\mathcal{F}$ is constructible as a $\Lambda$-sheaf.
  \item $H^i(U,\mathcal{F})$ is a finitely generated $\Lambda$-module.
 \end{itemize}
 Then the action of $K$ on $H^i(U,\mathcal{F})$
 is smooth and the $\mathcal{D}_c(K)$-module structure on $H^i(U/K,\mathcal{F})$ can be identified with
 that on $i_\mathcal{D}H^i(U,\mathcal{F})$. In particular,
 we have $i_\mathcal{D}\infty_\mathcal{D}H^i(U/K,\mathcal{F})=H^i(U/K,\mathcal{F})$.
\end{prop}

\begin{prf}
 We can take an increasing sequence $U_1\subset U_2\subset \cdots$ of quasi-compact open subsets of $X$
 which are stable under $K$ such that $U=\bigcup_{n=1}^\infty U_n$. As in the proof of
 Proposition \ref{prop:countable-union}, we have the following exact sequence:
 \[
 0\longrightarrow \varprojlim_n\!{}^1 H^{i-1}(U_n,\mathcal{F})\longrightarrow H^i(U,\mathcal{F})
 \longrightarrow \varprojlim_n H^i(U_n,\mathcal{F})\longrightarrow 0.
 \]
 As $\mathcal{F}$ is constructible, $H^{i-1}(U_n,\mathcal{F})$ is a finitely generated $\Lambda$-module for
 every $n$ by \cite[Proposition 3.1]{MR1620118} (in the case where the characteristic of $k$ is $0$)
 or \cite[Proposition 6.1.1, (1.7.7)]{MR1734903} (in the case where $X$ is smooth over $\Spa(k,k^+)$).
 Therefore, as $\Lambda$ is Artinian, the projective system $(H^{i-1}(U_n,\mathcal{F}))_n$ satisfies
 the Mittag-Leffler condition, and thus $\varprojlim_n^1 H^{i-1}(U_n,\mathcal{F})=0$.
 Hence we have an isomorphism $H^i(U,\mathcal{F})\yrightarrow{\cong} \varprojlim_n H^i(U_n,\mathcal{F})$.
 
 On the other hand, by the assumption, $H^i(U,\mathcal{F})$ is an Artinian $\Lambda$-module.
 Therefore the decreasing series of $\Lambda$-submodules
 $(\Ker(H^i(U,\mathcal{F})\to H^i(U_n,\mathcal{F})))_n$ is stationary.
 Hence, for a large enough $n$, the map $H^i(U,\mathcal{F})\longrightarrow H^i(U_n,\mathcal{F})$ is
 injective.
 By Corollary \ref{cor:derived-same-Gamma} and Remark \ref{rem:quasi-compact-smooth},
 the action of $K$ on $H^i(U_n,\mathcal{F})=H^i(U_n/K,\mathcal{F})$ is smooth.
 Thus, the action of $K$ on $H^i(U,\mathcal{F})$ is also smooth.

 Consider a $\mathcal{D}_c(K)$-homomorphism $H^i(U/K,\mathcal{F})\longrightarrow H^i(U_n/K,\mathcal{F})$,
 which fits into the following diagram by Corollary \ref{cor:derived-same-Gamma}:
 \[
  \xymatrix{%
 H^i(U/K,\mathcal{F})\ar[r]\ar@{=}[d]& H^i(U_n/K,\mathcal{F})\ar@{=}[d]\\
 H^i(U,\mathcal{F})\ar[r]& H^i(U_n,\mathcal{F})\lefteqn{.}
 }
 \]
 This map is injective if $n$ is large enough.
 Hence $H^i(U/K,\mathcal{F})$ satisfies 
 \[
  \infty_\mathcal{D}H^i(U/K,\mathcal{F})=H^i(U,\mathcal{F}),\quad
 i_\mathcal{D}\infty_\mathcal{D}H^i(U/K,\mathcal{F})=H^i(U/K,\mathcal{F}).
 \]
 This concludes the proof.
\end{prf}

\begin{prop}\label{prop:injective-Gamma_c-acyclic}
 Let $U$ be an open subspace of $X$ which is partially proper over $\Spa(k,k^+)$.
 For an injective object $\mathcal{F}$ in $\Lambda\text{-}\widetilde{X}_\et/G$,
 we have $H^i_c(U,\mathcal{F})=0$ for $i\ge 1$.
\end{prop}

\begin{prf}
 By the same way as in the proof of Proposition \ref{prop:H_c-limit} i), we can prove that
 $H^i_c(U,\mathcal{F})\cong \varinjlim_{Z\subset U}H^i_Z(U,\mathcal{F})$, where $Z$ runs through
 quasi-compact closed subsets of $U$.
 Let $\mathcal{C}$ be the set consisting of closed subsets of $U$ of the form $V\setminus W$ where
 $V$ and $W$ are quasi-compact open subsets of $U$ with $W\subset V$.
 We prove that this set is cofinal in the set of
 all quasi-compact closed subsets of $U$. Let $Z$ be an arbitrary quasi-compact closed subset of $U$.
 We can take a quasi-compact open subset $V$ of $U$ containing $Z$.
 As $U$ is partially proper over $\Spa(k,k^+)$, the closure $\overline{V}$ of $V$ in $U$ is quasi-compact
 (\cf \cite[Lemma 1.3.13]{MR1734903}).
 Since $Z$ is closed, for each point $x\in \overline{V}\setminus V$, there exists a quasi-compact
 open neighborhood $W_x$ of $x$ in $U$ such that $W_x\cap Z=\varnothing$.
 As $\overline{V}\setminus V$ is quasi-compact,
 we can find finitely many $x_1,\ldots,x_n\in \overline{V}\setminus V$ so that $W_{x_1},\ldots,W_{x_n}$
 cover $\overline{V}\setminus V$. Put $W=V\cap (W_{x_1}\cup\cdots\cup W_{x_n})$, which is a quasi-compact
 open subset of $V$.
 Since $V\setminus W=\overline{V}\setminus (W_{x_1}\cup\cdots\cup W_{x_n})$ is closed in $U$,
 it gives an element of $\mathcal{C}$ containing $Z$.

 Therefore we have $H^i_c(U,\mathcal{F})\cong \varinjlim_{Z\in\mathcal{C}}H^i_Z(U,\mathcal{F})$. 
 On the other hand, Proposition \ref{cor:local-coh-vanish} tells us that
 $H^i_Z(U,\mathcal{F})=0$ for $Z\in\mathcal{C}$ and $i\ge 1$. This concludes the proof.
\end{prf}

\begin{lem}
 Assume that $X$ is partially proper over $\Spa(k,k^+)$.
 For an object $\mathcal{F}$ in $\Lambda\text{-}\widetilde{X}_\et/G$,
 the induced $G$-action on $\Gamma_c(X,\mathcal{F})$ is smooth.
\end{lem}

\begin{prf}
 We have $\Gamma_c(X,\mathcal{F})=\varinjlim_{V\subset X}\Gamma_{\overline{V}}(X,\mathcal{F})$, where
 $V$ runs through quasi-compact open subsets of $X$ and $\overline{V}$ denotes the closure of $V$
 (note that $\overline{V}$ is quasi-compact).
 Therefore it suffices to show that the action of $K_V$ on $\Gamma_{\overline{V}}(X,\mathcal{F})$
 is smooth. 
 As $\overline{V}$ is quasi-compact, we can take a quasi-compact open subset $U$ of $X$ containing $\overline{V}$.
 Put $K=K_U\cap K_V$. Then the homomorphisms
 $\Gamma_{\overline{V}}(X,\mathcal{F})\yrightarrow{\cong}\Gamma_{\overline{V}}(U,\mathcal{F})\hooklongrightarrow \Gamma(U,\mathcal{F})$ are $K$-equivariant. Since the action of $K$ on $\Gamma(U,\mathcal{F})$ is smooth,
 so is the action on $\Gamma_{\overline{V}}(X,\mathcal{F})$.
\end{prf}

\begin{defn}
 Assume that $X$ is partially proper over $\Spa(k,k^+)$.
 By the lemma above,
 a left exact functor
 \[
 \Gamma_c(X/G,-)\colon \Lambda\text{-}\widetilde{X}_\et/G\longrightarrow \Rep_\Lambda(G);\quad \mathcal{F}\longmapsto\Gamma_c(X,\mathcal{F})
 \]
 is induced. We denote by $R\Gamma_c(X/G,-)$ the right derived functor of $\Gamma_c(X/G,-)$,
 and by $H^i_c(X/G,-)$ the $i$th cohomology of $R\Gamma_c(X/G,-)$.
\end{defn}

\begin{cor}\label{cor:derived-same}
 Assume that $X$ is partially proper over $\Spa(k,k^+)$.
 The following diagram is 2-commutative:
 \[
 \xymatrix{%
 D^+(\Lambda\text{-}\widetilde{X}_\et/G)\ar[rr]^-{R\Gamma_c(X/G,-)}\ar[d]&& D^+\bigl(\Rep_\Lambda(G)\bigr)\ar[d]\\
 D^+(\Lambda\text{-}\widetilde{X}_\et)\ar[rr]^-{R\Gamma_c(X,-)}&& D^+\bigl(\Mod(\Lambda)\bigr)\lefteqn{.}
 }
 \]
 Here the vertical arrows denote the forgetful functors.
\end{cor}

\begin{prf}
 It follows immediately from Proposition \ref{prop:injective-Gamma_c-acyclic}.
\end{prf}

\subsubsection{The Godement resolution}
Here we introduce the Godement resolution for a smooth $G$-equivariant sheaf.

\begin{defn}\label{defn:Godement}
 For each $x\in X$, fix a geometric point $i_{\overline{x}}\colon \overline{x}\longrightarrow X$
 lying over $x$.
 For an arbitrary $\Lambda$-sheaf $\mathcal{F}$ on $X_\et$,
 consider the smooth $G$-equivariant $\Lambda$-sheaf 
 \[
  \mathcal{C}(\mathcal{F})=\Ind_{X/\{1\}}^{X/G}\Bigl(\prod_{x\in X}i_{\overline{x}*}\mathcal{F}_{\overline{x}}\Bigr).
 \]
 If $\mathcal{F}$ is an object of $\Lambda\text{-}\widetilde{X}_\et/G$, the canonical morphism
 $\mathcal{F}\longrightarrow \prod_{x\in X}i_{\overline{x}*}\mathcal{F}_{\overline{x}}$
 in $\Lambda\text{-}\widetilde{X}_\et$ induces a morphism
 $\mathcal{F}\longrightarrow \mathcal{C}(\mathcal{F})$ in $\Lambda\text{-}\widetilde{X}_\et/G$.
 It is an injection; indeed, $\mathcal{C}(\mathcal{F})$ is a subsheaf of
 $\prod_{g\in G}\prod_{x\in X}g^*i_{\overline{x}*}\mathcal{F}_{\overline{x}}$, and the natural morphism
 $\mathcal{F}\longrightarrow \prod_{g\in G}\prod_{x\in X}g^*i_{\overline{x}*}\mathcal{F}_{\overline{x}}$
 is obviously injective.

 By repeating this construction, we have the following functorial resolution
 \[
 0\longrightarrow \mathcal{F}\longrightarrow \mathcal{C}^0(\mathcal{F})\longrightarrow \mathcal{C}^1(\mathcal{F})\longrightarrow \cdots,
 \]
 which is called the Godement resolution of $\mathcal{F}$.
\end{defn}

\begin{prop}\label{prop:Godement-injective}
 Let $\mathcal{F}$ be an object of $\Lambda\text{-}\widetilde{X}_\et/G$ which is flat as a $\Lambda$-sheaf.
 For every $i\ge 0$, $\mathcal{C}^i(\mathcal{F})$ is flat and injective
 in the category $\Lambda\text{-}\widetilde{X}_\et/G$.
 Moreover, for the maximal ideal $\mathfrak{m}$ of $\Lambda$, we have
 $\mathcal{C}^i(\mathcal{F})\otimes_{\Lambda}\Lambda/\mathfrak{m}\cong \mathcal{C}^i(\mathcal{F}\otimes_{\Lambda}\Lambda/\mathfrak{m})$.
\end{prop}

\begin{prf}
 By Lemma \ref{lem:TDVR-injective}, $\mathcal{F}_{\overline{x}}$ is an injective $\Lambda$-module.
 Therefore $\prod_{x\in X}i_{\overline{x}*}\mathcal{F}_{\overline{x}}$ is an injective $\Lambda$-sheaf.
 As $\Ind_{X/\{1\}}^{X/G}$ preserves injective objects, $\mathcal{C}(\mathcal{F})$ is an injective object
 in the category $\Lambda\text{-}\widetilde{X}_\et/G$.

 Next we prove that $\mathcal{C}(\mathcal{F})$ is flat. For an object $Y\longrightarrow X$ in $\widetilde{X}_{\qcet}$, we have
 \[
  \Gamma\bigl(Y,\mathcal{C}(\mathcal{F})\bigr)=\varinjlim_{K\subset K_Y} \Bigl(\prod_{g\in G}\prod_{x\in X}\Gamma(Y,g^*i_{\overline{x}*}\mathcal{F}_{\overline{x}})\Bigr)^K
 =\varinjlim_{K\subset K_Y} \Bigl(\prod_{g\in \Omega_K}\prod_{x\in X}\Gamma(Y,g^*i_{\overline{x}*}\mathcal{F}_{\overline{x}})\Bigr),
 \]
 where $K$ runs through compact open subgroups of $K_Y$ and $\Omega_K$ is a system of representatives of $G/K$.
 As $\Gamma(Y,g^*i_{\overline{x}*}\mathcal{F}_{\overline{x}})$ is a finite direct sum of $\mathcal{F}_{\overline{x}}$, it is a flat $\Lambda$-module. Since flatness of $\Lambda$-modules is preserved by arbitrary direct products and filtered inductive limits (\cf Lemma \ref{lem:TDVR-injective}), we can conclude that
 $\Gamma(Y,\mathcal{C}(\mathcal{F}))$ is flat. Therefore each stalk of $\mathcal{C}(\mathcal{F})$ is also flat.

 By the description of $\Gamma(Y,\mathcal{C}(\mathcal{F}))$ above,
 the functor $\mathcal{F}\longmapsto \mathcal{C}(\mathcal{F})$ is exact.
 If we take a generator $\lambda$ of $\mathfrak{m}$, we have an exact sequence
 $\mathcal{F}\yrightarrow{\times\lambda}\mathcal{F}\longrightarrow \mathcal{F}\otimes_{\Lambda}\Lambda/\mathfrak{m}\longrightarrow 0$.
 Therefore the sequence
 $\mathcal{C}(\mathcal{F})\yrightarrow{\times\lambda}\mathcal{C}(\mathcal{F})\longrightarrow \mathcal{C}(\mathcal{F}\otimes_{\Lambda}\Lambda/\mathfrak{m})\longrightarrow 0$ is exact, 
 and we obtain $\mathcal{C}(\mathcal{F})\otimes_{\Lambda}\Lambda/\mathfrak{m}\cong \mathcal{C}(\mathcal{F}\otimes_{\Lambda}\Lambda/\mathfrak{m})$.

 Put $\mathcal{G}=\Coker(\mathcal{F}\to\mathcal{C}(\mathcal{F}))$.
 For each $x\in X$, we have an exact sequence 
 $0\longrightarrow \mathcal{F}_{\overline{x}}\longrightarrow \mathcal{C}(\mathcal{F})_{\overline{x}}\longrightarrow \mathcal{G}_{\overline{x}}\longrightarrow 0$. As $\mathcal{F}_{\overline{x}}$ and $\mathcal{C}(\mathcal{F})_{\overline{x}}$ are flat ($=$ injective), so is $\mathcal{G}_{\overline{x}}$. 
 Thus $\mathcal{G}$ is flat.
 Moreover, we have $\mathcal{G}\otimes_\Lambda\Lambda/\mathfrak{m}=\Coker(\mathcal{F}\otimes_\Lambda\Lambda/{\mathfrak{m}}\rightarrow \mathcal{C}(\mathcal{F}\otimes_\Lambda\Lambda/{\mathfrak{m}}))$.
 Now we can repeat the same argument.
\end{prf}

The following proposition is also needed in the next section (\cf \cite[Lemme 2.7]{Fargues-Zel}).

\begin{prop}\label{prop:truncated-Godement}
 Assume that $k$ is separably closed and $X$ is $d$-dimensional.
 Let $\mathcal{F}$ be an object in $\Lambda\text{-}\widetilde{X}_\et/G$ which is flat as a $\Lambda$-sheaf,
 and $0\longrightarrow \mathcal{F}\longrightarrow \mathcal{C}^\bullet(\mathcal{F})$ the Godement resolution
 of $\mathcal{F}$. Put $\mathcal{G}=\Ker(\mathcal{C}^{2d}(\mathcal{F})\to\mathcal{C}^{2d+1}(\mathcal{F}))$.
 
 Then, for every open subspace $U$ of $X$ which is partially proper over $\Spa(k,k^+)$,
 $\Gamma_c(U,\mathcal{C}^m(\mathcal{F}))$ and $\Gamma_c(U,\mathcal{G})$ are free $\Lambda$-modules
 and $H^i_c(U,\mathcal{G})=0$ for $i\ge 1$.
\end{prop}

To prove this proposition, we use the following lemma.

\begin{lem}\label{lem:Gamma_c-flat}
 Assume that $k$ is separably closed and $X$ is finite-dimensional.
 Let $U$ be an open subspace of $X$ which is partially proper over $\Spa(k,k^+)$.
 Let $\mathcal{F}$ be a flat $\Lambda$-sheaf on $X_\et$ satisfying $H^i_c(U,\mathcal{F})=0$ for every $i\ge 1$.
 Then the following are equivalent:
 \begin{enumerate}
  \item $\Gamma_c(U,\mathcal{F})$ is a free $\Lambda$-module.
  \item $H^1_c(U,\mathcal{F}\otimes_{\Lambda}\Lambda/\mathfrak{m})=0$, where $\mathfrak{m}$ is the maximal
	ideal of $\Lambda$.
  \item $H^i_c(U,\mathcal{F}\otimes_{\Lambda}\Lambda/\mathfrak{m})=0$ for every $i\ge 1$.
 \end{enumerate}
\end{lem}

\begin{prf}
 As $R\Gamma_c(U,-)$ is bounded, we have
 \[
  R\Gamma_c(U,\mathcal{F})\Lotimes_\Lambda\Lambda/\mathfrak{m}\cong R\Gamma_c(U,\mathcal{F}\Lotimes_{\Lambda}\Lambda/\mathfrak{m}).
 \]
 By the conditions on $\mathcal{F}$, the left hand side is equal to 
 $\Gamma_c(U,\mathcal{F})\Lotimes_\Lambda\Lambda/\mathfrak{m}$, and the right hand side is equal to
 $R\Gamma_c(U,\mathcal{F}\otimes_{\Lambda}\Lambda/\mathfrak{m})$. In particular we have
 \[
  \Tor_i^\Lambda\bigl(\Lambda/\mathfrak{m},\Gamma_c(U,\mathcal{F})\bigr)\cong H^i_c(U,\mathcal{F}\otimes_{\Lambda}\Lambda/\mathfrak{m})
 \]
 for every $i$.
 Therefore we have i)$\implies$iii)$\implies$ii)$\implies$i), as desired.
\end{prf}

\begin{prf}[of Proposition \ref{prop:truncated-Godement}]
 First we prove that $\Gamma_c(U,\mathcal{C}^m(\mathcal{F}))$ is free.
 By Proposition \ref{prop:Godement-injective}, $\mathcal{C}^m(\mathcal{F}\otimes_\Lambda\Lambda/\mathfrak{m})$
 is an injective object of $(\Lambda/\mathfrak{m})\text{-}\widetilde{X}_\et/G$.
 Therefore, by Proposition \ref{prop:injective-Gamma_c-acyclic} and Proposition \ref{prop:Godement-injective}
 we have
 \[
 H^i_c\bigl(U,\mathcal{C}^m(\mathcal{F})\otimes_\Lambda\Lambda/\mathfrak{m}\bigr)\cong 
 H^i_c\bigl(U,\mathcal{C}^m(\mathcal{F}\otimes_\Lambda\Lambda/\mathfrak{m})\bigr)=0
 \]
 for $i\ge 1$. Hence Lemma \ref{lem:Gamma_c-flat} tells us that $\Gamma_c(U,\mathcal{C}^m(\mathcal{F}))$ is free.
 
 Next we show that $\Gamma_c(U,\mathcal{G})$ is free and $H^i_c(U,\mathcal{G})=0$ for $i\ge 1$.
 For simplicity we put $\mathcal{J}^m=\mathcal{C}^m(\mathcal{F})$ for $0\le m\le 2d-1$.
 We have an exact sequence
\[
 0\longrightarrow \mathcal{F}\longrightarrow \mathcal{J}^0\longrightarrow \cdots\longrightarrow \mathcal{J}^{2d-1}\longrightarrow \mathcal{G}\longrightarrow 0.
\]
 Thus, by Proposition \ref{prop:Godement-injective} and Proposition \ref{prop:injective-Gamma_c-acyclic} we have
 $H^i_c(U,\mathcal{G})\cong H^{i+2d}_c(U,\mathcal{F})=0$ for every $i\ge 1$
(\cf \cite[Proposition 5.3.11]{MR1734903}).
 On the other hand, by the proof of Proposition \ref{prop:Godement-injective},
 we know that $\mathcal{G}$ is a flat $\Lambda$-sheaf.
 Therefore we have an exact sequence
 \[
  0\to \mathcal{F}\otimes_\Lambda\Lambda/\mathfrak{m}\to \mathcal{J}^0\otimes_\Lambda\Lambda/\mathfrak{m}\to \cdots\to \mathcal{J}^{2d-1}\otimes_\Lambda\Lambda/\mathfrak{m}\to \mathcal{G}\otimes_\Lambda\Lambda/\mathfrak{m}\to 0.
 \]
 In the beginning of the proof, we have proved that
 $H^i_c(U,\mathcal{J}^m\otimes_\Lambda\Lambda/\mathfrak{m})=0$ for $i\ge 1$ and $0\le m\le 2d-1$.
 This implies that $H^1_c(U,\mathcal{G}\otimes_\Lambda\Lambda/\mathfrak{m})\cong H^{1+2d}_c(U,\mathcal{F}\otimes_\Lambda\Lambda/\mathfrak{m})=0$. Thus Lemma \ref{lem:Gamma_c-flat} tells us that $\Gamma_c(U,\mathcal{G})$ is a free
 $\Lambda$-module.
\end{prf}

\section{Duality theorem}\label{sec:duality}
Let $R$, $\kappa$, $F$ and $k$ be as in Section \ref{subsec:formal-adic}.
We denote the residue characteristic of $R$ by $p$.
Fix a truncated discrete valuation ring $\Lambda$ with residue characteristic $\neq p$
and a locally pro-$p$ group $G$.

\begin{thm}\label{thm:duality-thm}
 Let $\mathcal{X}$ be a special formal scheme over $\Spf R$ equipped with a continuous action of $G$
 in the sense of \cite[D\'efinition 2.3.10]{MR2074714}.
 Assume the following:
 \begin{enumerate}
  \item[(a)] The rigid geometric generic fiber $X=t(\mathcal{X})_{\overline{\eta}}$ of $\mathcal{X}$
	is purely $d$-dimensional and smooth over $\Spa(k,k^+)$.
  \item[(b)] $\mathcal{X}$ is locally algebraizable (\cf \cite[Definition 3.19]{formalnearby}).
  \item[(c)] $\mathcal{X}^\red$ is partially proper over $\kappa$.
  \item[(d)] There exists a quasi-compact open subset $V$ of $\mathcal{X}^\red$ such that
	$\mathcal{X}^\red=\bigcup_{g\in G}gV$ and $\{g\in G\mid g\overline{V}\cap \overline{V}\neq \varnothing\}$ is compact.
 \end{enumerate}
 Then, for each integer $i$, we have a $G$-equivariant isomorphism
 \[
  H^{2d+i}_{c,\mathcal{X}}(X,\Lambda)(d)\yrightarrow{\cong}R^i\!\D\bigl(R\Gamma_c(X/G,\Lambda)\bigr).
 \]
 Note that the condition (c) ensures that $X$ is partially proper over $\Spa(k,k^+)$ (\cf Proposition \ref{prop:generic-fiber-partially-proper}), and thus Corollary \ref{cor:derived-same} and \cite[Proposition 5.3.11]{MR1734903} imply that $R\Gamma_c(X/G,\Lambda)$ lies in the category $D^b(\Rep_\Lambda(G))$.
\end{thm}

Fix $V$ satisfying the condition (d) in the theorem above and put $Z=\overline{V}$, $U=\spp^{-1}(Z)^\circ$.
Take a pro-$p$ open subgroup $K$ which stabilizes $V$. Then $Z$ and $U$ are also stable under $K$.
The condition (d) tells us that $\{gZ\}_{g\in G/K}$ is a locally finite covering of $\mathcal{X}^\red$.
Note also that $\{gU\}_{g\in G/K}$ is a locally finite covering of $X$.
By Proposition \ref{prop:tube-property} iii), $U$ is partially proper over $\Spa(k,k^+)$. 

Let $0\longrightarrow \Lambda\longrightarrow \mathcal{I}^\bullet$ be the Godement resolution introduced
in Definition \ref{defn:Godement} for the smooth $G$-equivariant constant sheaf $\Lambda$ on $X_\et$.
By Proposition \ref{prop:Godement-injective}, this is an injective resolution in the category
$\Lambda\text{-}\widetilde{X}_\et/G$.
Put $\mathcal{J}^\bullet=\tau_{\le 2d}\mathcal{I}^\bullet$.
For $\alpha=(\overline{g}_1,\ldots,\overline{g}_m)\in (G/K)^m$, let $U_\alpha$ and $K_\alpha$ be
as in Proposition \ref{prop:functorial-resolution}. 
We write $j_\alpha$ for the natural open immersion $U_\alpha\hooklongrightarrow X$. 

The following lemma gives a complex which computes $R\!\D(R\Gamma_c(X/G,\Lambda))$.

\begin{lem}\label{lem:complex-representative-C}
 Let $C^{\bullet\bullet}$ denote the the double complex 
 \[
 C^{\bullet\bullet}=\Gamma_c\bigl(X/G,C^\bullet(\mathcal{J}^\bullet)\bigr)=\bigoplus_{\alpha\in G\backslash(G/K)^{-\bullet+1}}\cInd_{K_\alpha}^G\Gamma_c(U_\alpha,\mathcal{J}^\bullet)
 \]
  in $\Rep_\Lambda(G)$
 (\cf Proposition \ref{prop:functorial-resolution} i)).
 Then we have an isomorphism 
 \[
 R\Gamma_c(X/G,\Lambda)\cong \Tot C^{\bullet\bullet}
 \]
 in $D^-(\Rep_\Lambda(G))$.
 Moreover, we have an isomorphism
 \[
  R\!\D\bigl(R\Gamma_c(X/G,\Lambda)\bigr)\cong \D(\Tot C^{\bullet\bullet})
 \]
 in $D^+(\Rep_\Lambda(G))$.
\end{lem}

\begin{prf}
 By Proposition \ref{prop:functorial-resolution} i), $\Tot C^\bullet(\mathcal{J}^\bullet)$ gives
 a resolution of $\Lambda$ in $\Lambda\text{-}\widetilde{X}_\et/G$.
 By Corollary \ref{cor:derived-same} and Proposition \ref{prop:truncated-Godement}, for integers $i\ge 1$ and $m,n\ge 0$ we have
 \begin{align*}
  H^i_c\bigl(X/G,C^m(\mathcal{J}^n)\bigr)&=\bigoplus_{\alpha\in G\backslash(G/K)^{-m+1}}\cInd_{K_\alpha}^G H^i_c(X,j_{\alpha!}j_\alpha^*\mathcal{J}^n)\\
  &=\bigoplus_{\alpha\in G\backslash(G/K)^{-m+1}}\cInd_{K_\alpha}^G H^i_c(U_\alpha,\mathcal{J}^n)\\
  &=0.
 \end{align*}
 Therefore, each component of the complex $\Tot C^\bullet(\mathcal{J}^\bullet)$ is acyclic with respect to
 $\Gamma_c(X/G,-)$. Hence we have
 \[
  R\Gamma_c(X/G,\Lambda)\cong\Gamma_c\bigl(X/G,\Tot C^\bullet(\mathcal{J}^\bullet)\bigr)
 =\Tot \Gamma_c\bigl(X/G,C^\bullet(\mathcal{J}^\bullet)\bigr)=\Tot C^{\bullet\bullet}
 \]
 in $D^-(\Rep_\Lambda(G))$.

 Furthermore, Proposition \ref{prop:truncated-Godement} tells us that $\Gamma_c(U_\alpha,\mathcal{J}^n)$
 is a projective object in $\Rep_\Lambda(K_\alpha)$. Therefore 
 $C^{mn}=\bigoplus_{\alpha\in G\backslash(G/K)^{-m+1}}\cInd_{K_\alpha}^G\Gamma_c(U_\alpha,\mathcal{J}^n)$
 is a projective object in $\Rep_\Lambda(G)$.
 Hence we have $R\!\D(R\Gamma_c(X/G,\Lambda))\cong \D(\Tot C^{\bullet\bullet})$. 
\end{prf}

\begin{lem}\label{lem:Gamma_c-sp-acyclic}
 For an injective object $\mathcal{F}$ of $\Lambda\text{-}\widetilde{X}_\et/G$,
 we have 
 \[
  R\Gamma_{c,\mathcal{X}}(X,j_{\alpha*}j_\alpha^*\mathcal{F})=\Gamma(U_\alpha,\mathcal{F}).
 \]
 In particular, $j_{\alpha*}j_\alpha^*\mathcal{F}$ is acyclic with respect to 
 $\Gamma_{c,\mathcal{X}}(X,-)$.
\end{lem}

\begin{prf}
 For $\alpha=(\overline{g}_1,\ldots,\overline{g}_m)\in (G/K)^m$, we set $Z_\alpha=g_1Z\cap \cdots\cap g_mZ$. 
 Then we have $U_\alpha=\spp^{-1}(Z_\alpha)^\circ$. Therefore $U_\alpha$ is the union of countably many
 quasi-compact open subsets by Proposition \ref{prop:tube-property} ii).

 Note that $\Gamma_{c,\mathcal{X}}(X,j_{\alpha*}(-))=\Gamma_c(\mathcal{X}^\red,-)\circ \spp_*\circ j_{\alpha*}=\Gamma(\mathcal{X}^\red,-)\circ \spp_*\circ j_{\alpha*}$;
 indeed, for a sheaf $\mathcal{G}$ on $(U_\alpha)_\et$, all elements of
 $\Gamma(\mathcal{X}^\red,\spp_*j_{\alpha*}\mathcal{G})$
 are supported on the quasi-compact closed subset $Z_\alpha\subset\mathcal{X}^\red$.
 On the other hand, by Proposition \ref{prop:countable-union} and Corollary \ref{cor:j_*-acyclic},
 we have $R\Gamma(U_\alpha,\mathcal{F})=\Gamma(U_\alpha,\mathcal{F})$ and 
 $Rj_{\alpha*}j^*_\alpha\mathcal{F}=j_{\alpha*}j^*_\alpha\mathcal{F}$, respectively.
 Hence we have
 \begin{align*}
  R\Gamma_{c,\mathcal{X}}(X,j_{\alpha*}j_\alpha^*\mathcal{F})&=R\Gamma_c(\mathcal{X}^\red,R\spp_*j_{\alpha*}j^*_\alpha\mathcal{F})=R\Gamma_c(\mathcal{X}^\red,R\spp_*Rj_{\alpha*}j^*_\alpha\mathcal{F})\\
  &=R\Gamma(\mathcal{X}^\red,R\spp_*Rj_{\alpha*}j^*_\alpha\mathcal{F})=R\Gamma(U_\alpha,\mathcal{F})
  =\Gamma(U_\alpha,\mathcal{F}),
 \end{align*}
 as desired.
\end{prf}

\begin{lem}\label{lem:non-smooth-functorial-resolution}
 For an object $\mathcal{F}$ of $\Lambda\text{-}\widetilde{X}_\et/G$,
 we have a functorial resolution $\mathcal{F}\longrightarrow D^\bullet(\mathcal{F})$ in
 $\Lambda\text{-}\widetilde{X}_\et/{G^{\mathrm{disc}}}$ where $D^m(\mathcal{F})$ is given by
 \[
  D^m(\mathcal{F})=\prod_{\alpha\in G\backslash(G/K)^{m+1}}\prod_{g\in K_\alpha\backslash G}g^*j_{\alpha*}j_\alpha^*\mathcal{F}.
 \]
 Recall that $G^{\mathrm{disc}}$ denotes the group $G$ with discrete topology.

 Moreover, $\Gamma_{c,\mathcal{X}}(X,D^m({\mathcal{F}}))$ is naturally equipped with a structure of
 a $\mathcal{D}_c(G)$-module under which
 \[
  \Gamma_{c,\mathcal{X}}\bigl(X,D^m({\mathcal{F}})\bigr)=\bigoplus_{\alpha\in G\backslash(G/K)^{m+1}}\cInd_{\mathcal{D}_c(K_\alpha)}^{\mathcal{D}_c(G)}\Gamma(U_\alpha,\mathcal{F})
 \]
 and $\Gamma_{c,\mathcal{X}}(X,D^\bullet({\mathcal{F}}))$ is a complex in $\Mod(\mathcal{D}_c(G))$.
\end{lem}

\begin{prf}
  Consider the complex
 \begin{equation*}
 0\longrightarrow \mathcal{F}\longrightarrow \prod_{\lambda\in G/K}j_{\lambda*}j_\lambda^*\mathcal{F}\longrightarrow \prod_{\lambda\in (G/K)^2}j_{\lambda*}j_\lambda^*\mathcal{F}\longrightarrow \cdots\tag{$*$}
 \end{equation*}
 of $\Lambda$-sheaves over $X_\et$. It is well-known that this complex is exact
 (\cf the proof of Proposition \ref{prop:functorial-resolution} ii)).
 Each term has a $G$-equivariant structure and each homomorphism 
 is $G$-equivariant. 
 Clearly we have
 \[
 \prod_{\lambda\in (G/K)^{m+1}}j_{\lambda*}j_\lambda^*\mathcal{F}
 =\prod_{\alpha\in G\backslash(G/K)^{m+1}}\prod_{g\in K_\alpha\backslash G}g^*j_{\alpha*}j_\alpha^*\mathcal{F}.
 \]
 Hence we have a desired resolution.

 As $\{gZ\}_{g\in G/K}$ is a locally finite covering of $\mathcal{X}_\red$,
 we have
 \begin{align*}
  \spp_*\Bigl(\prod_{\alpha\in G\backslash(G/K)^{m+1}}\prod_{g\in K_\alpha\backslash G}g^*j_{\alpha*}j_\alpha^*\mathcal{F}\Bigr)
  &=\prod_{\alpha\in G\backslash(G/K)^{m+1}}\prod_{g\in K_\alpha\backslash G}\spp_*g^*j_{\alpha*}j_\alpha^*\mathcal{F}\\
  &=\bigoplus_{\alpha\in G\backslash(G/K)^{m+1}}\bigoplus_{g\in K_\alpha\backslash G}\spp_*g^*j_{\alpha*}j_\alpha^*\mathcal{F}.
 \end{align*}
 Therefore, by Corollary \ref{cor:H_c-direct-sum} and Lemma \ref{lem:Gamma_c-sp-acyclic}, we obtain
 \begin{align*}
  \Gamma_{c,\mathcal{X}}\bigl(X,D^m(\mathcal{F})\bigr)&=\Gamma_c\bigl(\mathcal{X}_\red,\spp_*D^m(\mathcal{F})\bigr)\\
  &=\bigoplus_{\alpha\in G\backslash(G/K)^{m+1}}\bigoplus_{g\in K_\alpha\backslash G}\Gamma_c(\mathcal{X}_\red,\spp_*g^*j_{\alpha*}j_\alpha^*\mathcal{F})\\
  &=\bigoplus_{\alpha\in G\backslash(G/K)^{m+1}}\bigoplus_{g\in K_\alpha\backslash G}\Gamma(g^{-1}U_\alpha,\mathcal{F}).
 \end{align*}
 By Lemma \ref{lem:section-D_c-action}, $\Gamma(U_\alpha,\mathcal{F})$ has a natural structure of
 a $\mathcal{D}_c(K_\alpha)$-module and we have
 \[
 \Gamma_{c,\mathcal{X}}\bigl(X,D^m(\mathcal{F})\bigr)=\bigoplus_{\alpha\in G\backslash(G/K)^{m+1}}\cInd_{\mathcal{D}_c(K_\alpha)}^{\mathcal{D}_c(G)}\Gamma(U_\alpha,\mathcal{F})
 \]
 as $G$-modules. Hence $\Gamma_{c,\mathcal{X}}(X,D^m(\mathcal{F}))$ can be regarded as
 a $\mathcal{D}_c(G)$-module.

 It is easy to see that the homomorphism
 \[
  \Gamma_{c,\mathcal{X}}\bigl(X,D^m(\mathcal{F})\bigr)\longrightarrow \Gamma_{c,\mathcal{X}}(X,D^{m+1}(\mathcal{F}))
 \]
 is $\mathcal{D}_c(G)$-linear. Indeed, it follows from the following simple fact:
 \begin{quote}
  for $\alpha=(\overline{g}_1,\ldots,\overline{g}_{m+2})\in (G/K)^{m+2}$ and $\beta\in (G/K)^{m+1}$
  which can be obtained by removing one of the entries $\overline{g}_1,\ldots,\overline{g}_{m+2}$
  from $\alpha$, the restriction map
  $\Gamma(U_\beta,\mathcal{F})\longrightarrow \Gamma(U_\alpha,\mathcal{F})$
  is $\mathcal{D}_c(K_\alpha)$-linear.
 \end{quote}
 This completes the proof.
\end{prf}

By these lemmas, we can give a complex which represents $R\Gamma_{c,\mathcal{X}}(X,\Lambda)$.

\begin{cor}\label{cor:complex-representative-D_0}
 We write $D_0^{\bullet\bullet}$ for the double complex 
 \[
  \Gamma_{c,\mathcal{X}}\bigl(X,D^\bullet(\mathcal{I}^\bullet)\bigr)
 =\bigoplus_{\alpha\in G\backslash(G/K)^{\bullet+1}}\cInd_{\mathcal{D}_c(K_\alpha)}^{\mathcal{D}_c(G)}\Gamma(U_\alpha,\mathcal{I}^\bullet)
 \]
 in $\Mod(\mathcal{D}_c(G))$.
 Then we have a $G$-equivariant isomorphism 
 \[
  R\Gamma_{c,\mathcal{X}}(X,\Lambda)\cong \Tot D_0^{\bullet\bullet}
 \]
 in $D^+(\Mod(\Lambda))$.
\end{cor}

\begin{prf}
 By Lemma \ref{lem:non-smooth-functorial-resolution}, $\Tot D^\bullet(\mathcal{I}^\bullet)$ gives
 a resolution of $\Lambda$ in $\Lambda\text{-}\widetilde{X}_\et/G^{\mathrm{disc}}$.
 In the same way as in the proof of Lemma \ref{lem:non-smooth-functorial-resolution},
 we can deduce from Lemma \ref{lem:Gamma_c-sp-acyclic} that
 \[
  R\Gamma_{c,\mathcal{X}}\bigl(X,D^m(\mathcal{I}^n)\bigr)
 =\bigoplus_{\alpha\in G\backslash(G/K)^{m+1}}\cInd_{\mathcal{D}_c(K_\alpha)}^{\mathcal{D}_c(G)}\Gamma(U_\alpha,\mathcal{I}^n)
 \]
 (see also Corollary \ref{cor:H_c-direct-sum}).
 Therefore, each component of $\Tot D^\bullet(\mathcal{I}^\bullet)$ is
 acyclic with respect to $\Gamma_{c,\mathcal{X}}(X,-)$.
 Hence we have
 \[
 R\Gamma_{c,\mathcal{X}}(X,\Lambda)\cong \Gamma_{c,\mathcal{X}}\bigl(X,\Tot D^\bullet(\mathcal{I}^\bullet)\bigr)=\Tot D_0^{\bullet\bullet}
 \]
 in $D^+(\Mod(\Lambda))$.
\end{prf}

\begin{lem}\label{lem:complex-representative-D}
 We put $D^{\bullet\bullet}=\infty_{\mathcal{D}}D_0^{\bullet\bullet}$, which is a double complex in
 $\Rep_\Lambda(G)$.
 Then we have a $G$-equivariant isomorphism 
 \[
  R\Gamma_{c,\mathcal{X}}(X,\Lambda)\cong \Tot D^{\bullet\bullet}
 \]
 in $D^+(\Mod(\Lambda))$.
\end{lem}

\begin{prf}
 We need to prove that the natural homomorphism of complexes 
 $\Tot D^{\bullet\bullet}\longrightarrow \Tot D_0^{\bullet\bullet}$ is a quasi-isomorphism.
 Consider the following morphism of spectral sequences:
 \[
  \xymatrix{%
 E_1^{m,n}=H^n(D^{m\bullet})\ar@{=>}[r]\ar[d]&H^{m+n}(\Tot D^{\bullet\bullet})\ar[d]\\
 E_1^{m,n}=H^n(D_0^{m\bullet})\ar@{=>}[r]&H^{m+n}(\Tot D_0^{\bullet\bullet})\lefteqn{.}
 }
 \]
 It suffices to show that $H^n(D^{m\bullet})\longrightarrow H^n(D_0^{m\bullet})$ is an isomorphism.
 By definition, the $i$th cohomology of $\Gamma(U_\alpha,\mathcal{I}^\bullet)$
 is $H^i(U_\alpha/K_\alpha,\Lambda)$. Therefore we have
 \begin{align*}
  H^n(D_0^{m\bullet})&=\bigoplus_{\alpha\in G\backslash(G/K)^{m+1}}\cInd_{\mathcal{D}_c(K_\alpha)}^{\mathcal{D}_c(G)}H^n(U_\alpha/K_\alpha,\Lambda),\\
  H^n(D^{m\bullet})&=\bigoplus_{\alpha\in G\backslash(G/K)^{m+1}}\infty_{\mathcal{D}}\bigl(\cInd_{\mathcal{D}_c(K_\alpha)}^{\mathcal{D}_c(G)}H^n(U_\alpha/K_\alpha,\Lambda)\bigr)\\
  &=\bigoplus_{\alpha\in G\backslash(G/K)^{m+1}}\cInd_{K_\alpha}^G\bigl(\infty_{\mathcal{D}}H^n(U_\alpha/K_\alpha,\Lambda)\bigr)\\
  &=\bigoplus_{\alpha\in G\backslash(G/K)^{m+1}}\cInd_{\mathcal{D}_c(K_\alpha)}^{\mathcal{D}_c(G)}\bigl(i_{\mathcal{D}}\infty_{\mathcal{D}}H^n(U_\alpha/K_\alpha,\Lambda)\bigr).
 \end{align*}
 By Proposition \ref{prop:constructible-smooth} and Proposition \ref{prop:tube-property} iv), we have
 $i_{\mathcal{D}}\infty_{\mathcal{D}}H^n(U_\alpha/K_\alpha,\Lambda)=H^n(U_\alpha/K_\alpha,\Lambda)$.
 This concludes the proof.
\end{prf}

\begin{lem}\label{lem:C-tilde}
 Let $\widetilde{C}^{\bullet\bullet}$ be the double complex in $\Mod(\mathcal{D}_c(G))$ given by
 \[
  \bigoplus_{\alpha\in G\backslash(G/K)^{\bullet+1}}\cInd_{\mathcal{D}_c(K_\alpha)}^{\mathcal{D}_c(G)}\Gamma_c(U_\alpha,\mathcal{J}^\bullet)^*.
 \]
 Then there exists a natural morphism of double complexes $\widetilde{C}^{\bullet\bullet}\longrightarrow \D^m(C^{\bullet\bullet})$ that induces a quasi-isomorphism
 $\Tot \widetilde{C}^{\bullet\bullet}\longrightarrow \D^m(\Tot C^{\bullet\bullet})$.
\end{lem}

\begin{prf}
 By Proposition \ref{prop:dual-cInd}, we have a natural morphism
 $\widetilde{C}^{\bullet\bullet}\longrightarrow \D^m(C^{\bullet\bullet})$.
 To observe that $\Tot \widetilde{C}^{\bullet\bullet}\longrightarrow \D^m(\Tot C^{\bullet\bullet})$ is
 a quasi-isomorphism, it suffices to show that
 the morphism $\widetilde{C}^{n\bullet}\longrightarrow \D^m(C^{-n,\bullet})$
 is a quasi-isomorphism for each $n$. We have
 \[
  H^{-i}(\widetilde{C}^{n\bullet})=\bigoplus_{\alpha\in G\backslash(G/K)^{n+1}}\cInd_{\mathcal{D}_c(K_\alpha)}^{\mathcal{D}_c(G)}H^i_c(U_\alpha,\Lambda)^*.
 \]
 To compute the cohomology of $\D^m(C^{-n,\bullet})$, note the following points:
 \begin{itemize}
  \item Since the set $\{\overline{g}\in G/K\mid U\cap gU\neq \varnothing\}$ is finite, 
	for a fixed $n$ there exist only finitely many $\alpha\in G\backslash(G/K)^{n+1}$ such that
	$U_\alpha\neq\varnothing$. Therefore the direct sum $\bigoplus_{\alpha\in G\backslash(G/K)^{n+1}}$ 
	in $C^{-n,\bullet}$ is finite and commutes with $\D^m$.
  \item By Proposition \ref{prop:truncated-Godement}, we have
	\[
	\D^m\bigl(\cInd_{K_\alpha}^G\Gamma_c(U_\alpha,\mathcal{J}^\bullet)\bigr)
	=R\!\D^m\bigl(\cInd_{K_\alpha}^G\Gamma_c(U_\alpha,\mathcal{J}^\bullet)\bigr).
	\]
  \item By Proposition \ref{prop:tube-property} iv), 
	$H^i_c(U_\alpha,\Lambda)$ is a finitely generated $\Lambda$-module. 
	Therefore, by Proposition \ref{prop:dual-cInd},
	$\cInd_{K_\alpha}^G H^i_c(U_\alpha,\Lambda)$ is acyclic
	with respect to $\D^m$.
 \end{itemize}
 Therefore we have
 \begin{align*}
  H^{-i}\bigl(\D^m(C^{-n,\bullet})\bigr)&=\bigoplus_{\alpha\in G\backslash(G/K)^{n+1}}R^{-i}\!\D^m\bigl(\cInd_{K_\alpha}^G\Gamma_c(U_\alpha,\mathcal{J}^\bullet)\bigr)\\
  &=\bigoplus_{\alpha\in G\backslash(G/K)^{n+1}}\D^m\bigl(\cInd_{K_\alpha}^G H^i_c(U_\alpha,\Lambda)\bigr).
 \end{align*}
 Since $H^i_c(U_\alpha,\Lambda)$ is a finitely generated $\Lambda$-module, 
 Proposition \ref{prop:dual-cInd} tells us that the homomorphism
 \[
  \cInd_{\mathcal{D}_c(K_\alpha)}^{\mathcal{D}_c(G)}H^i_c(U_\alpha,\Lambda)^*\longrightarrow 
 \D^m\bigl(\cInd_{K_\alpha}^GH^i_c(U_\alpha,\Lambda)\bigr)
 \]
 is an isomorphism. Hence the homomorphism 
 $H^{-i}(\widetilde{C}^{n\bullet})\longrightarrow H^{-i}(\D^m(C^{-n,\bullet}))$
 is an isomorphism, as desired.
\end{prf}

\begin{lem}\label{lem:Poincare-duality}
 Let $D_0^{\bullet\bullet}$ be as in Corollary \ref{cor:complex-representative-D_0}, and
 $\widetilde{C}^{\bullet\bullet}$ as in Lemma \ref{lem:C-tilde}.
 We have a natural isomorphism
 \[
  \Tot(D_0^{\bullet\bullet})(d)[2d]\yrightarrow{\cong} \Tot(\widetilde{C}^{\bullet\bullet})
 \]
 in $D^+(\Mod(\mathcal{D}_c(G)))$.
\end{lem}

\begin{prf}
 First we fix $\alpha\in (G/K)^{m+1}$, and construct a morphism
 \[
  \Gamma(U_\alpha,\mathcal{I}^\bullet)(d)[2d]\longrightarrow \Gamma_c(U_\alpha,\mathcal{J}^\bullet)^*
 \]
 in $D^+(\Mod(\mathcal{D}_c(K_\alpha)))$.
 Since $\mathcal{I}^n$ is flat for each $n$, the complex $\Tot(\mathcal{J}^\bullet\otimes\mathcal{I}^\bullet)$
 gives a resolution of $\Lambda$. Therefore we have a morphism of complexes
 $\Tot(\mathcal{J}^\bullet\otimes\mathcal{I}^\bullet)\longrightarrow \mathcal{I}^\bullet$,
 which is determined up to homotopy. The cup product
 \[
  \Gamma_c(U_\alpha,\mathcal{J}^\bullet)\otimes_\Lambda\Gamma(U_\alpha,\mathcal{I}^\bullet)(d)[2d]
 \longrightarrow \Gamma_c(U_\alpha,\mathcal{J}^\bullet\otimes\mathcal{I}^\bullet)(d)[2d]
 \]
 gives a morphism of complexes
 \[
  \Gamma(U_\alpha,\mathcal{I}^\bullet)(d)[2d]\longrightarrow 
 \Hom\Bigl(\Gamma_c(U_\alpha,\mathcal{J}^\bullet),\Gamma_c\bigl(U_\alpha,\Tot(\mathcal{J}^\bullet\otimes\mathcal{I}^\bullet)\bigr)(d)[2d]\Bigr).
 \]
 Consider the following morphisms of complexes:
 \begin{align*}
  &\Gamma_c\bigl(U_\alpha,\Tot(\mathcal{J}^\bullet\otimes\mathcal{I}^\bullet)\bigr)(d)[2d]
  \longrightarrow \Gamma_c(U_\alpha,\mathcal{I}^\bullet)(d)[2d]\longrightarrow \Gamma_c(X,\mathcal{I}^\bullet)(d)[2d]\\
  &\qquad\longrightarrow \tau_{\ge 0}\Bigl(\Gamma_c(X,\mathcal{I}^\bullet)(d)[2d]\Bigr)\yleftarrow{(*)}H^{2d}_c(X,\Lambda)(d)\yrightarrow{\Tr_X}\Lambda.
 \end{align*}
 The morphism $(*)$ is a quasi-isomorphism, as $H_c^i(X,\Lambda)(d)=0$ for $i>2d$.
 By composing these morphisms, we get morphisms of complexes
 \begin{align*}
  \Gamma(U_\alpha,\mathcal{I}^\bullet)(d)[2d]&\longrightarrow \Hom\bigl(\Gamma_c(U_\alpha,\mathcal{J}^\bullet),\tau_{\ge 0}(\Gamma_c(X,\mathcal{I}^\bullet)(d)[2d])\bigr)\\
  &\yleftarrow{(*)}\Hom\bigl(\Gamma_c(U_\alpha,\mathcal{J}^\bullet),H^{2d}_c(X,\Lambda)(d)\bigr)\\
  &\longrightarrow \Hom\bigl(\Gamma_c(U_\alpha,\mathcal{J}^\bullet),\Lambda\bigr)=\Gamma_c(U_\alpha,\mathcal{J}^\bullet)^*.
 \end{align*}
 Since $\Gamma_c(U_\alpha,\mathcal{J}^\bullet)$ consists of free $\Lambda$-modules
 (\cf Proposition \ref{prop:truncated-Godement}), $(*)$ is a quasi-isomorphism.
 As in \cite[Lemme 2.6]{Fargues-Zel}, it is easy to show that these morphisms are
 $\mathcal{D}_c(K_\alpha)$-linear.
 Moreover, by Corollary \ref{cor:derived-same-Gamma} and Proposition \ref{prop:truncated-Godement},
 the $(-i)$th cohomology of the composite morphism
 $\Gamma(U_\alpha,\mathcal{I}^\bullet)(d)[2d]\longrightarrow \Gamma_c(U_\alpha,\mathcal{J}^\bullet)^*$
 in the derived category $D^+(\Mod(\mathcal{D}_c(K_\alpha)))$ 
 is by definition the isomorphism of the Poincar\'e duality
 \[
 H^{2d-i}(U_\alpha,\Lambda)(d)\yrightarrow{\cong} H^i_c(U_\alpha,\Lambda)^*.
 \]
 Hence we obtain an isomorphism 
 \[
  \Gamma(U_\alpha,\mathcal{I}^\bullet)(d)[2d]\yrightarrow{\cong}\Gamma_c(U_\alpha,\mathcal{J}^\bullet)^*
 \]
 in $D^+(\Mod(\mathcal{D}_c(K_\alpha)))$.

 Put
 \begin{align*}
 \widetilde{C}^{\bullet\bullet}_1&=\bigoplus_{\alpha\in G\backslash(G/K)^{\bullet+1}}\cInd_{\mathcal{D}_c(K_\alpha)}^{\mathcal{D}_c(G)}\Hom\bigl(\Gamma_c(U_\alpha,\mathcal{J}^\bullet),\tau_{\ge 0}(\Gamma_c(X,\mathcal{I}^\bullet)(d)[2d])\bigr),\\
  \widetilde{C}^{\bullet\bullet}_2&=\bigoplus_{\alpha\in G\backslash(G/K)^{\bullet+1}}\cInd_{\mathcal{D}_c(K_\alpha)}^{\mathcal{D}_c(G)}\Hom\bigl(\Gamma_c(U_\alpha,\mathcal{J}^\bullet),H^{2d}_c(X,\Lambda)(d)\bigr),
 \end{align*}
 which are double complexes in $\Mod(\mathcal{D}_c(G))$.
 Then, by the construction above, we obtain morphisms of double complexes
 \[
 D_0^{\bullet\bullet}(d)[2d]\longrightarrow \widetilde{C}^{\bullet\bullet}_1\yleftarrow{(*)}\widetilde{C}^{\bullet\bullet}_2\longrightarrow \widetilde{C}^{\bullet\bullet}.
 \]
 As $(*)$ induces a quasi-isomorphism $\widetilde{C}^{n\bullet}_1\longleftarrow\widetilde{C}^{n\bullet}_2$
 for each $n$, the morphism $\Tot(*)$ is also a quasi-isomorphism. Similarly,
 we can conclude that the composite morphism
 \[
  \Tot D_0^{\bullet\bullet}(d)[2d]\longrightarrow \Tot\widetilde{C}^{\bullet\bullet}
 \]
 in $D^+(\Mod(\mathcal{D}_c(G)))$ is an isomorphism.
\end{prf} 

\begin{prf}[of Theorem \ref{thm:duality-thm}]
 By Lemma \ref{lem:C-tilde} and Lemma \ref{lem:Poincare-duality}, we have isomorphisms
 \[
  \Tot D^{\bullet\bullet}(d)[2d]\cong \infty_{\mathcal{D}}(\Tot\widetilde{C}^{\bullet\bullet})
 \cong \infty_{\mathcal{D}}\bigl(\D^m(\Tot C^{\bullet\bullet})\bigr)=\D(\Tot C^{\bullet\bullet})
 \]
 in $D^+(\Rep_\Lambda(G))$. Therefore, by Lemma \ref{lem:complex-representative-C} and
 Lemma \ref{lem:complex-representative-D}, we have a $G$-equivariant isomorphism
 \[
 R\Gamma_{c,\mathcal{X}}(X,\Lambda)(d)[2d]\cong R\!\D\bigl(R\Gamma_c(X/G,\Lambda)\bigr)
 \]
 in $D^+(\Mod(\Lambda))$.
 By taking cohomology, we get the desired isomorphism
 \[
  H^{2d+i}_{c,\mathcal{X}}(X,\Lambda)(d)\yrightarrow{\cong}R^i\!\D\bigl(R\Gamma_c(X/G,\Lambda)\bigr).
 \]
\end{prf}

To construct the isomorphism above, we chose $V$ and $K$.
Next we prove that the isomorphism is independent of these choices.

\begin{prop}\label{prop:independence}
 The isomorphism
 \[
  H^{2d+i}_{c,\mathcal{X}}(X,\Lambda)(d)\yrightarrow{\cong}R^i\!\D\bigl(R\Gamma_c(X/G,\Lambda)\bigr)
 \]
 in Theorem \ref{thm:duality-thm} is independent of the choice of $V$ and $K$.
\end{prop}

\begin{prf}
 We denote the isomorphism attached to $V$ and $K$ by $f_{V,K}$.
 Let $V'$ and $K'$ be another choice. We should prove that $f_{V,K}=f_{V',K'}$.
 We may assume either $V=V'$ or $K=K'$. Indeed, if we obtain the equality
 in this case, then in general we have
 \[
  f_{V,K}=f_{V,K\cap K'}=f_{V\cup V',K\cap K'}=f_{V',K\cap K'}=f_{V',K'}
 \]
 as desired. In particular, we may assume that $V\subset V'$ and $K\subset K'$.
 
 We put $Z'=\overline{V'}$ and $U'=\spp^{-1}(Z')^\circ$. For $\beta\in (G/K')^m$, 
 we define $U'_\beta$ and $K'_\beta$ in the same way as in
 Proposition \ref{prop:functorial-resolution}.
 Since the open covering $\{gU\}_{g\in G/K}$ is a refinement of $\{gU'\}_{g\in G/K'}$,
 there exists a natural morphism of double complexes
 \[
  \prod_{\beta\in (G/K')^{\bullet+1}}j'_{\beta*}j'^*_\beta\mathcal{I}^\bullet
 \longrightarrow \prod_{\alpha\in (G/K)^{\bullet+1}}j_{\alpha*}j_\alpha^*\mathcal{I}^\bullet.
 \]
 Here $j'_\beta$ denotes the open immersion $U'_\beta\hooklongrightarrow X$.
 This morphism turns out to be a $G$-equivariant morphism of double complexes
 $D'^\bullet(\mathcal{I}^\bullet)\longrightarrow D^\bullet(\mathcal{I}^\bullet)$,
 where $D'^\bullet(\mathcal{I}^n)$ denotes the complex $D^\bullet(\mathcal{I}^n)$ in
 Lemma \ref{lem:non-smooth-functorial-resolution} attached to $V'$ and $K'$.
 Put $D'^{\bullet\bullet}_0=\Gamma_{c,\mathcal{X}}(X,D'^\bullet(\mathcal{I}^\bullet))$,
 which is a double complex in $\Mod(\mathcal{D}_c(G))$
 by Lemma \ref{lem:non-smooth-functorial-resolution}.
 It is easy to see that the induced morphism 
 $D'^{\bullet\bullet}_0\longrightarrow D^{\bullet\bullet}_0$ is $\mathcal{D}_c(G)$-equivariant.
 Note that the following diagram is commutative, where $(*)$ (resp.\ $(**)$) is induced from 
 the augmentation morphism $\mathcal{I}^\bullet\longrightarrow D'^\bullet(\mathcal{I}^\bullet)$
 (resp.\ $\mathcal{I}^\bullet\longrightarrow D^\bullet(\mathcal{I}^\bullet)$)
 in Lemma \ref{lem:non-smooth-functorial-resolution}:
 \[
  \xymatrix{%
 &\Gamma_{c,\mathcal{X}}(X,\mathcal{I}^\bullet)\ar[ld]_-{(*)}\ar[rd]^-{(**)}\\
 \Tot D'^{\bullet\bullet}_0\ar[rr]&&\Tot D^{\bullet\bullet}_0\lefteqn{.}
 }
 \]
 Hence we obtain the commutative diagram below, where we put
 $D'^{\bullet\bullet}=\infty_{\mathcal{D}}D'^{\bullet\bullet}_0$
 (\cf Lemma \ref{lem:complex-representative-D}):
 \[
 \xymatrix{%
 &H^i_{c,\mathcal{X}}(X,\Lambda)\ar[ld]_-{\cong}\ar[rd]^-{\cong}\\
 H^i(\Tot D'^{\bullet\bullet})\ar[rr]&&H^i(\Tot D^{\bullet\bullet})\lefteqn{.}
 }
 \]

 Let $C'^\bullet(\mathcal{J}^\bullet)$ denote the double complex 
 \[
  \bigoplus_{\beta\in G\backslash(G/K')^{-\bullet+1}}\cInd_{U'_\beta/K'_\beta}^{X/G}\Res_{U'_\beta/K'_\beta}^{X/G}\mathcal{J}^\bullet
 \]
 attached to $U'$ and $K'$ defined in Proposition \ref{prop:functorial-resolution} i).
 Then we have a morphism $C^\bullet(\mathcal{J}^\bullet)\longrightarrow C'^\bullet(\mathcal{J}^\bullet)$
 of double complexes in $\Lambda\text{-}\widetilde{X}_\et/G$.
 This morphism is clearly compatible with the augmentation morphisms
 $C^\bullet(\mathcal{J}^\bullet)\longrightarrow \mathcal{J}^\bullet$ and 
 $C'^\bullet(\mathcal{J}^\bullet)\longrightarrow \mathcal{J}^\bullet$ (\cf Proposition \ref{prop:functorial-resolution} i)).
 Put
 \[
  C'^{\bullet\bullet}=\Gamma_c\bigl(X/G,C'^\bullet(\mathcal{J}^\bullet)\bigr)
 =\bigoplus_{\beta\in G\backslash(G/K')^{-\bullet+1}}\cInd_{K'_\beta}^G\Gamma_c(U'_\beta,\mathcal{J}^\bullet).
 \]
 Then, a $G$-equivariant morphism $C'^{\bullet\bullet}\longrightarrow C^{\bullet\bullet}$ is induced, and the following diagrams are commutative (\cf Lemma \ref{lem:complex-representative-C}):
 \[
  \xymatrix{%
 \Tot C^{\bullet\bullet}\ar[rr]\ar[rd]_-{\cong}&&\Tot C'^{\bullet\bullet}\ar[ld]^-{\cong}\\
 & R\Gamma_c(X/G,\Lambda)\lefteqn{,}
 }
 \]
 \[
  \xymatrix{%
 &R\!\D\bigl(R\Gamma_c(X/G,\Lambda)\bigr)\ar[ld]_-{\cong}\ar[rd]^-{\cong}\\
 \D(\Tot C'^{\bullet\bullet})\ar[rr]&&\D(\Tot C^{\bullet\bullet})\lefteqn{.}
 }
 \]
 As in Lemma \ref{lem:C-tilde} and the proof of Lemma \ref{lem:Poincare-duality}, we put
 \begin{align*}
  \widetilde{C}'^{\bullet\bullet}&=\bigoplus_{\beta\in G\backslash(G/K')^{\bullet+1}}\cInd_{\mathcal{D}_c(K'_\beta)}^{\mathcal{D}_c(G)}\Gamma_c(U'_\beta,\mathcal{J}^\bullet)^*,\\
  \widetilde{C}_1'^{\bullet\bullet}&=\bigoplus_{\beta\in G\backslash(G/K')^{\bullet+1}}\cInd_{\mathcal{D}_c(K'_\beta)}^{\mathcal{D}_c(G)}\Hom\bigl(\Gamma_c(U'_\beta,\mathcal{J}^\bullet),\tau_{\ge 0}(\Gamma_c(X,\mathcal{I}^\bullet)(d)[2d])\bigr),\\
  \widetilde{C}_2'^{\bullet\bullet}&=\bigoplus_{\beta\in G\backslash(G/K')^{\bullet+1}}\cInd_{\mathcal{D}_c(K'_\beta)}^{\mathcal{D}_c(G)}\Hom\bigl(\Gamma_c(U'_\beta,\mathcal{J}^\bullet),H^{2d}_c(X,\Lambda)(d)\bigr).
 \end{align*}
 Then, $\mathcal{D}_c(G)$-equivariant morphisms
$\widetilde{C}'^{\bullet\bullet}\longrightarrow \widetilde{C}^{\bullet\bullet}$, $\widetilde{C}_1'^{\bullet\bullet}\longrightarrow \widetilde{C}_1^{\bullet\bullet}$ and $\widetilde{C}_2'^{\bullet\bullet}\longrightarrow \widetilde{C}_2^{\bullet\bullet}$ are naturally induced. Furthermore, we can easily check
 that the following diagram is commutative (the horizontal arrows are the morphisms appeared in Lemma \ref{lem:C-tilde} and the proof of Lemma \ref{lem:Poincare-duality}):
 \[
  \xymatrix{%
 D_0'^{\bullet\bullet}(d)[2d]\ar[r]\ar[d]&\widetilde{C}'^{\bullet\bullet}_1\ar[d]&\widetilde{C}'^{\bullet\bullet}_2\ar[r]\ar[l]\ar[d]&\widetilde{C}'^{\bullet\bullet}\ar[r]\ar[d]&\D^m(C'^{\bullet\bullet})\ar[d]\\
 D_0^{\bullet\bullet}(d)[2d]\ar[r]&
 \widetilde{C}^{\bullet\bullet}_1&\widetilde{C}^{\bullet\bullet}_2\ar[r]\ar[l]
 &\widetilde{C}^{\bullet\bullet}\ar[r]&\D^m(C^{\bullet\bullet})\lefteqn{.}
 }
 \]
 Putting all together, we obtain the commutative diagram
 \[
  \xymatrix{%
 H^{2d+i}_{c,\mathcal{X}}(X,\Lambda)(d)\ar[d]^-{\cong}\ar@/_5pc/[dd]_-{\cong}&R^i\!\D\bigl(R\Gamma_c(X/G,\Lambda)\bigr)\ar[d]_-{\cong}\ar@/^5pc/[dd]^-{\cong}\\
 H^i(\Tot D'^{\bullet\bullet}(d)[2d])\ar[r]\ar[d]& H^i\bigl(\D(\Tot C'^{\bullet\bullet})\bigr)\ar[d]\\
 H^i(\Tot D^{\bullet\bullet}(d)[2d])\ar[r]& H^i\bigl(\D(\Tot C^{\bullet\bullet})\bigr)\lefteqn{,}
 }
 \]
 which gives the desired equality $f_{V,K}=f_{V',K'}$.
\end{prf}

\begin{cor}\label{cor:functoriality}
 For an isomorphism $\varphi\colon \mathcal{X}\yrightarrow{\cong}\mathcal{X}$ of formal schemes
 which is compatible with the action of $G$ on $\mathcal{X}$, we have the following commutative diagram:
 \[
  \xymatrix{%
 H^{2d+i}_{c,\mathcal{X}}(X,\Lambda)(d)\ar[r]^-{\cong}\ar[d]^-{\varphi^*}_-{\cong}&
 R^i\!\D\bigl(R\Gamma_c(X/G,\Lambda)\bigr)\ar[d]^-{R^i\!\D(\varphi_*)}_-{\cong}\\
 H^{2d+i}_{c,\mathcal{X}}(X,\Lambda)(d)\ar[r]^-{\cong}&
 R^i\!\D\bigl(R\Gamma_c(X/G,\Lambda)\bigr)\lefteqn{.}
 }
 \]
\end{cor}

\begin{prf}
 We use the notation $f_{V,K}$ in the proof of Proposition \ref{prop:independence}.
 It is immediate to see that the following diagram is commutative:
 \[
  \xymatrix{%
 H^{2d+i}_{c,\mathcal{X}}(X,\Lambda)(d)\ar[rr]^-{f_{V,K}}_-{\cong}\ar[d]^-{\varphi^*}_-{\cong}&&
 R^i\!\D\bigl(R\Gamma_c(X/G,\Lambda)\bigr)\ar[d]^-{R^i\!\D(\varphi_*)}_-{\cong}\\
 H^{2d+i}_{c,\mathcal{X}}(X,\Lambda)(d)\ar[rr]^-{f_{\varphi^{-1}(V),K}}_-{\cong}&&
 R^i\!\D\bigl(R\Gamma_c(X/G,\Lambda)\bigr)\lefteqn{.}
 }
 \]
 Hence the corollary follows from Proposition \ref{prop:independence}.
\end{prf}

By standard argument as in \cite[\S 3]{Fargues-Zel}, we can show the similar results
for $\ell$-adic coefficients.

\begin{thm}\label{thm:duality-thm-l-adic}
 Let $\mathcal{X}$ and $G$ be as in Theorem \ref{thm:duality-thm}, 
 and $L_\lambda$ a finite extension of $\Q_\ell$. Then, for each integer $i$,
 we have a $G$-equivariant isomorphism
 \[
  H^{2d+i}_{c,\mathcal{X}}(X,L_\lambda)(d)\yrightarrow{\cong}R^i\!\D\bigl(R\Gamma_c(X/G,L_\lambda)\bigr).
 \]
 This isomorphism is functorial with respect to an automorphism of $\mathcal{X}$ which is compatible with
 the action of $G$ on $\mathcal{X}$.
\end{thm}

\begin{rem}\label{rem:coh-finite-type}
 Assume that $G$ is a quotient of the group $\mathbf{H}(\Q_p)$ for some connected reductive group $\mathbf{H}$
 over $\Q_p$. Then, for every field $L$ of characteristic $0$, the category $\Rep_L(G)$ is noetherian
 and has finite projective dimension.

 In this case, for a finite extension $L_\lambda$ of $\Q_\ell$, the \v{C}ech spectral sequence
 \[
  E_1^{i,j}=\bigoplus_{\alpha\in G\backslash (G/K)^{-i+1}}\cInd_{K_\alpha}^G H^j_c(U_\alpha,L_\lambda)\Longrightarrow H^{i+j}_c(X,L_\lambda)
 \]
 (\cf Lemma \ref{lem:complex-representative-C}) and Proposition \ref{prop:tube-property} iv) tell us that 
 $R\Gamma_c(X/G,L_\lambda)$ is an object of $D^b_{\mathrm{fg}}(\Rep_{L_\lambda}(G))$
 (recall that $U_\alpha=\varnothing$ for all but finitely many $\alpha\in G\backslash(G/K)^{-i+1}$).
 Hence, by Corollary \ref{cor:derived-dual} iii), we have a $G$-equivariant functorial isomorphism
 \[
  H^{2d+i}_{c,\mathcal{X}}(X,\overline{\Q}_\ell)(d)\yrightarrow{\cong}R^i\!\D\bigl(R\Gamma_c(X/G,\overline{\Q}_\ell)\bigr).
 \]
\end{rem}

\section{Application to the Rapoport-Zink tower}\label{sec:application}
\subsection{Rapoport-Zink tower for $\GSp(2n)$}\label{subsec:RZ-tower}
Let $n\ge 1$ be an integer. 
For a ring $A$, let $\langle\ ,\ \rangle\colon A^{2n}\times A^{2n}\longrightarrow A$
be the symplectic pairing defined by
$\langle (x_i),(y_i)\rangle=x_1y_{2n}+x_2y_{2n-1}+\cdots +x_ny_{n+1}-x_{n+1}y_n-\cdots-x_{2n}y_1$,
and $\GSp_{2n}(A)$ the symplectic similitude group with respect to $\langle\ ,\ \rangle$.

Here we briefly recall the definition of the Rapoport-Zink tower for $\GSp(2n)$.
See \cite[\S 3.1]{RZ-LTF} for details. In this section, we assume that $p\neq 2$.

We fix a $n$-dimensional $p$-divisible group $\X$ over $\overline{\F}_p$
which is isoclinic of slope $1/2$ and a polarization $\lambda_0\colon \X\yrightarrow{\cong}\X^\vee$.
We write $\Z_{p^\infty}$ for the completion of the maximal unramified extension of $\Z_p$
and $\Nilp$ for the category of $\Z_{p^\infty}$-schemes on which $p$ is locally nilpotent.
Let $\M\colon \Nilp\longrightarrow \mathbf{Set}$ be the moduli functor of
deformations by quasi-isogenies $(X,\rho)$ of $(\X,\lambda_0)$
(for the precise definition, see \cite[\S 3.1]{RZ-LTF}).
It is known that $\M$ is represented by a special formal scheme over $\Spf\Z_{p^\infty}$, which is
also denoted by $\M$. By \cite[Proposition 2.32]{MR1393439}, every irreducible component of
$\M^\red$ is projective over $\overline{\F}_p$. In particular, $\M^\red$ is partially proper
over $\overline{\F}_p$.
We write $M$ for the rigid generic fiber $t(\M)_\eta$ of $\M$.
The adic space $M$ is purely $n(n+1)/2$-dimensional and smooth over $\Spa(\Q_{p^\infty},\Z_{p^\infty})$,
where $\Q_{p^\infty}=\Frac \Z_{p^\infty}$. By Proposition \ref{prop:generic-fiber-partially-proper},
$M$ is partially proper over $\Spa(\Q_{p^\infty},\Z_{p^\infty})$.

By adding level structures on the universal polarized $p$-divisible group on $M$, we can construct
a projective system of \'etale coverings $\{M_K\}_{K\subset K_0}$ of $M$, where $K$ runs through
open subgroups of $K_0=\GSp_{2n}(\Z_p)$. This projective system is called the Rapoport-Zink tower
for $\GSp(2n)$.
For each $K$ and $g\in G=\GSp_{2n}(\Q_p)$ with $g^{-1}Kg\subset K_0$,
we have a natural isomorphism $[g]\colon M_K\yrightarrow{\cong} M_{g^{-1}Kg}$ called the Hecke operator.
In particular, the group $G$ acts on the pro-object $\{M_K\}_{K\subset K_0}$ on the right.

Let $J$ be the group of self-quasi-isogenies of $\X$ preserving $\lambda_0$ up to multiplication by $\Q_p^\times$.
We can construct a connected reductive algebraic group $\mathbf{J}$ over $\Q_p$
in a natural way such that $\mathbf{J}(\Q_p)=J$.
In particular, $J$ is naturally equipped with a topology. Concretely, $\mathbf{J}$ is isomorphic to 
$\mathrm{GU}(n,D)$, where $D$ is the quaternion division algebra over $\Q_p$ (\cf \cite[Remark 3.11]{RZ-LTF}).
By definition, $J$ acts on $\M$ and $M$. This action naturally extends to $M_K$ for each $K$ and
transition maps in the projective system $\{M_K\}_{K\subset K_0}$ are
compatible with the actions of $J$.
Further, the Hecke operators also commute with the actions of $J$.
By \cite[Corollaire 4.4.1]{MR2074714}, the action of $J$ on $M_K$ is continuous in the sense of
Definition \ref{defn:action-continuous}.
Sometimes it is convenient to consider the quotient $M_K/p^\Z$ of $M_K$ by the discrete subgroup $p^\Z$
of the center of $J$.

Fix a prime number $\ell$ which is different from $p$. Put
\begin{align*}
 H^i_c(M_K)&=H^i_c(M_K\otimes_{\Q_{p^\infty}}\overline{\Q}_{p^\infty},\overline{\Q}_\ell),& H^i_c(M_\infty)&=\varinjlim_{K\subset K_0}H^i_c(M_K),\\
 H^i_c(M_K/p^\Z)&=H^i_c\bigl((M_K/p^\Z)\otimes_{\Q_{p^\infty}}\overline{\Q}_{p^\infty},\overline{\Q}_\ell\bigr),& H^i_c(M_\infty/p^\Z)&=\varinjlim_{K\subset K_0}H^i_c(M_K/p^\Z).
\end{align*}
As $G\times J$ acts on the tower $\{M_K\}_{K\subset K_0}$, the $\overline{\Q}_\ell$-vector spaces
$H^i_c(M_\infty)$ and $H^i_c(M_\infty/p^\Z)$ are equipped with actions of $G\times J$.
The actions of $G$ are obviously smooth. The actions of $J$ are also smooth, as the action of $J$ on $M_K$
is smooth (\cf \cite[Corollaire 4.4.7]{MR2074714}).
By \cite[Lemma 5.36]{MR1393439}, the action of $G$ on $H^i_c(M_\infty/p^\Z)$ factors through $G/p^\Z$, 
where $p^\Z\subset \Q_p^\times\subset G$ is a discrete subgroup of the center of $G$.
For an open subgroup $K$ of $K_0$, we have 
\[
 H^i_c(M_\infty)^K=H^i_c(M_K),\quad H^i_c(M_\infty/p^\Z)^K=H^i_c(M_K/p^\Z).
\]
By using the Weil descent datum on $\M$ (\cf \cite[3.48]{MR1393439}), we can define actions
of the Weil group $W_{\Q_p}$ of $\Q_p$ on $H^i_c(M_K)$, $H^i_c(M_\infty)$, $H^i_c(M_K/p^\Z)$, 
and $H^i_c(M_\infty/p^\Z)$. Hence, $H^i_c(M_\infty)$ is a representation of $G\times J\times W_{\Q_p}$,
and $H^i_c(M_\infty/p^\Z)$ is a representation of $G/p^\Z\times J/p^\Z\times W_{\Q_p}$.

In the following, we also write $M$ and $M_K$ for $M\otimes_{\Q_{p^\infty}}\overline{\Q}_{p^\infty}$ and 
$M_K\otimes_{\Q_{p^\infty}}\overline{\Q}_{p^\infty}$ by abuse of notation.
Recall that the formal model $\M$ of $M$ gives a support set $\mathcal{C}_\M$ of $M$
(\cf Definition \ref{defn:formal-model-support-set}). We denote by the same symbol $\mathcal{C}_\M$
the support set of $M_K$ induced by the morphism $M_K\longrightarrow M$
(\cf Definition \ref{defn:support-set-mor} i)).
Similarly, the formal model $\M/p^\Z$ of $M/p^\Z$ determines a support set of $M_K/p^\Z$,
which is also denoted by $\mathcal{C}_\M$ for simplicity.

\begin{prop}\label{prop:Hecke-support-set}
 Let $K$ be an open subgroup of $K_0$, and $g$ an element of $G$ such that $g^{-1}Kg\subset K_0$.
 Then, the Hecke operator $[g]\colon M_K\yrightarrow{\cong} M_{g^{-1}Kg}$ 
 induces an isomorphism of pairs
 $(M_K,\mathcal{C}_\M)\yrightarrow{\cong} (M_{g^{-1}Kg},\mathcal{C}_\M)$
 (\cf Definition \ref{defn:support-set-mor} ii)).
 Similarly, we have an isomorphism of pairs
 $[g]\colon (M_K/p^\Z,\mathcal{C}_\M)\yrightarrow{\cong} (M_{g^{-1}Kg}/p^\Z,\mathcal{C}_\M)$.
\end{prop}

\begin{prf}
 It suffices to show that if $Z\subset M_K$ belongs to $\mathcal{C}_\M$,
 then $[g](Z)$ belongs to $\mathcal{C}_\M$.
 If $g$ lies in the center $Z_G$ of $G$, then the claim holds, because \cite[Lemma 5.36]{MR1393439}
 tells us that the Hecke action $[g]$ on $M=M_{K_0}$ extends to an automorphism of $\M$.
 Therefore, replacing $g$ by $zg$ with a suitable element $z\in p^\Z$, we may assume that
 $\Z_p^{2n}\subset g\Z_p^{2n}$. Take an integer $N\ge 0$ such that
 $g\Z_p^{2n}\subset p^{-N}\Z_p^{2n}$.

 For an integer $k\ge 0$, let $\M^{\red,k}$ be the subfunctor of $\M^{\red}$
 consisting of $(X,\rho)$ such that $p^k\rho$ and $p^k\rho^{-1}$ are isogenies.
 By \cite[Proposition 2.9]{MR1393439}, it is represented by a closed subscheme of $\M^{\red}$,
 which is quasi-compact (\cf \cite[Corollary 2.31]{MR1393439}).
 Clearly we have $\M^{\red}=\bigcup_{k\ge 0}\M^{\red,k}$. Since $Z\in \mathcal{C}_\M$,
 we can find $k\ge 0$ such that $Z\subset \spp_K^{-1}(\M^{\red,k})$,
 where $\spp_K$ denotes the composite $M_K\yrightarrow{p_K} M\yrightarrow{\spp}\M^\red$.
 It suffices to prove that $[g](Z)\subset \spp_{g^{-1}Kg}^{-1}(\M^{\red,k+N})$.

 Take a geometric point $x$ in $Z$. It corresponds to a triple $(X,\rho,\alpha)$,
 where $X$ is a $p$-divisible group over the valuation ring $\kappa(x)^+$,
 $\rho\colon \X\otimes_{\overline{\F}_p}(\kappa(x)^+/p)\longrightarrow X\otimes_{\kappa(x)^+}(\kappa(x)^+/p)$ is a quasi-isogeny, and $\alpha$ is a level structure 
 $\Z_p^{2d}\yrightarrow{\cong} T_p(X\otimes_{\kappa(x)^+}\kappa(x))\pmod{K}$.
 Put $y=[g](x)$. The pair $(X',\rho')$ corresponding to the point $p_{g^{-1}Kg}(y)$ in $M$ can be described as
 follows.
 Note that $\alpha$ gives a homomorphisms
 $(\Q_p/\Z_p)^{2n}\longrightarrow X\otimes_{\kappa(x)^+}\kappa(x)$, which is well-defined up to 
 $K$-action.
 Let $H$ be the scheme-theoretic closure in $X$ of the image of $g\Z_p^{2n}/\Z_p^{2n}$ 
 under this homomorphism.
 Then, it is a finite locally free subgroup scheme of $X[p^N]$, since $\kappa(x)^+$ is a valuation ring.
 By the definition of Hecke operators (\cf \cite[5.34]{MR1393439}), we have $X'=X/H$ and
 $\rho'=(\phi\bmod p)\circ\rho$, where $\phi\colon X\longrightarrow X'$ is the canonical isogeny.
 As $\spp_K(x)$ lies in $\M^{\red,k}$, $p^k\rho$ and $p^k\rho^{-1}$ are isogenies.
 Therefore, so is $p^k\rho'$. On the other hand, since $H\subset X[p^N]$,
 $\Ker \phi$ is killed by $p^N$, and thus $p^N\phi^{-1}$ is an isogeny. Hence
 $p^{k+N}\rho'^{-1}$ is also an isogeny. 
 Namely, $y=[g](x)$ lies in $\spp^{-1}_{g^{-1}Kg}(\M^{\red,k+N})$. 
 Now we conclude that $[g](Z)\subset \spp_{g^{-1}Kg}^{-1}(\M^{\red,k+N})$.
\end{prf}

Put $H^i_{\mathcal{C}_\M}(M_\infty)=\varinjlim_{K\subset K_0}H^i_{\mathcal{C}_\M}(M_K)$
and $H^i_{\mathcal{C}_\M}(M_\infty/p^\Z)=\varinjlim_{K\subset K_0}H^i_{\mathcal{C}_\M}(M_K/p^\Z)$.
By the previous proposition, $G$ acts naturally on $H^i_{\mathcal{C}_\M}(M_\infty)$
and $H^i_{\mathcal{C}_\M}(M_\infty/p^\Z)$. For an open subgroup $K$ of $K_0$, we have
\[
 H^i_{\mathcal{C}_\M}(M_\infty)^K=H^i_{\mathcal{C}_\M}(M_K),\quad
 H^i_{\mathcal{C}_\M}(M_\infty/p^\Z)^K=H^i_{\mathcal{C}_\M}(M_K/p^\Z).
\]
Obviously the groups $J$ and $W_{\Q_p}$ act on $H^i_{\mathcal{C}_\M}(M_\infty)$ and
$H^i_{\mathcal{C}_\M}(M_\infty/p^\Z)$.
We have natural homomorphisms 
$H^i_c(M_\infty)\longrightarrow H^i_{\mathcal{C}_\M}(M_\infty)$ and
$H^i_c(M_\infty/p^\Z)\longrightarrow H^i_{\mathcal{C}_\M}(M_\infty/p^\Z)$, which are
$G\times J\times W_{\Q_p}$-equivariant.

In the sequel, we will describe the action of $G$ on $H^i_{\mathcal{C}_\M}(M_\infty)$
by using some formal models (\cf \cite[\S 5.2]{RZ-GSp4}).
For an integer $m\ge 0$, let $\M_m$ be the formal scheme
classifying Drinfeld $m$-level structures on the universal $p$-divisible group on $\M$
(for a precise definition, see \cite[\S 3.2]{RZ-GSp4}). The formal scheme $\M_m$ is finite over $\M$,
and satisfies $t(\M_m)_\eta\cong M_{K_m}$, where $K_m$ is the kernel of
$\GSp_{2n}(\Z_p)\longrightarrow \GSp_{2n}(\Z/p^m\Z)$.

\begin{lem}\label{lem:formal-support}
 There exists a natural isomorphism
 \[
  H^i_{\mathcal{C}_\M}(M_{K_m})\cong H^i_c(\M_m^\red,R\Psi\overline{\Q}_\ell),
 \]
 where $R\Psi\overline{\Q}_\ell=R\spp_*\overline{\Q}_\ell$ denotes the formal nearby cycle complex
 (see \cite{MR1395723}).
\end{lem}

\begin{prf}
 By Proposition \ref{prop:formal-scheme-support}, we have
 \[
  H^i_c(\M_m^\red,R\Psi\overline{\Q}_\ell)=H^i_{c,\M_m}(M_{K_m})=H^i_{\mathcal{C}_{\M_m}}(M_{K_m}).
 \]
 Hence it suffices to show the equality $\mathcal{C}_{\M_m}=\mathcal{C}_\M$ of support sets of
 $M_{K_m}$. This immediately follows from Lemma \ref{lem:formal-finite-type} and the fact that
 $\M_m\longrightarrow \M$ is finite.
\end{prf}

Let $G^+$ denote the submonoid $\{g\in G\mid \Z_p^{2n}\subset g\Z_p^{2n}\}$ of $G$.
For $g\in G^+$, let $e(g)$ be the minimal non-negative integer such that $g\Z_p^{2n}\subset p^{-e(g)}\Z_p^{2n}$.
Following \cite[\S 6]{MR2169874}, for $g\in G^+$ with $m\ge e=e(g)$, we can define a formal scheme
$\M_{m,g}$ over $\M_m$ satisfying the following properties (\cf \cite[\S 5.2]{RZ-GSp4}):
\begin{itemize}
 \item The structure morphism $\pr\colon \M_{m,g}\longrightarrow \M_m$ is proper,
       and the induced morphism
       on rigid generic fibers $t(\M_{m,g})_\eta\longrightarrow t(\M_m)_\eta$ is an isomorphism.
 \item There exists a proper morphism $[g]\colon \M_{m,g}\longrightarrow \M_{m-e}$ such that
       the composite
       \[
       M_{K_m}=t(\M_m)_\eta\yrightarrow{\pr^{-1}}t(\M_{m,g})_\eta\yrightarrow{[g]}t(\M_{m-e})_\eta
       =M_{K_{m-e}}
       \]
       coincides with the Hecke operator $M_{K_m}\yrightarrow{[g]} M_{g^{-1}K_mg}\longrightarrow M_{K_{m-e}}$
       attached to $g\in G$
       (note that $\Z_p^{2n}\subset g\Z_p^{2n}\subset p^{-e}\Z_p^{2n}$ implies that
       $g^{-1}K_mg\subset K_{m-e}$).
\end{itemize}

\begin{lem}\label{lem:action-compare}
 For $g\in G^+$ and $m\ge e=e(g)$, the composite of
 \[
  H^i_c(\M^\red_{m-e},R\Psi\overline{\Q}_\ell)\yrightarrow{[g]^*}H^i_c(\M^\red_{m,g},R\Psi\overline{\Q}_\ell)\yrightarrow[\cong]{\pr_*}H^i_c(\M^\red_m,R\Psi\overline{\Q}_\ell)
 \]
 corresponds to the composite of
 \[
  H^i_{\mathcal{C}_\M}(M_{K_{m-e}})\longrightarrow H^i_{\mathcal{C}_\M}(M_{g^{-1}K_mg})\yrightarrow{[g]^*}H^i_{\mathcal{C}_\M}(M_{K_m})
 \]
 under the isomorphism in Lemma \ref{lem:formal-support}.
\end{lem}

\begin{prf}
 As in Lemma \ref{lem:formal-support}, we can see that
 \[
  H^i_c(\M^\red_{m,g},R\Psi\overline{\Q}_\ell)\cong H^i_{\mathcal{C}_\M}\bigl(t(\M_{m,g})_{\overline{\eta}}\bigr)\yrightarrow[\cong]{(\pr^*)^{-1}} H^i_{\mathcal{C}_\M}(M_{K_m}).
 \]
 Hence the lemma immediately follows from the property of $[g]\colon \M_{m,g}\longrightarrow \M_{m-e}$.
\end{prf}

By this description and the main theorem in \cite{RZ-GSp4}, we can compare the cuspidal parts
of $H^i_c(M_\infty)$ and $H^i_{\mathcal{C}_\M}(M_\infty)$ in the case $n=2$.

\begin{cor}\label{cor:compare-cuspidal-parts}
 Assume that $n=2$. Then, no supercuspidal representation of $G$ appears
 as a subquotient of the kernel and the cokernel of the map 
 $H^i_c(M_\infty)\longrightarrow H^i_{\mathcal{C}_\M}(M_\infty)$.
 The same holds for the map $H^i_c(M_\infty/p^\Z)\longrightarrow H^i_{\mathcal{C}_\M}(M_\infty/p^\Z)$.
\end{cor}

\begin{prf}
 We will use the notation in \cite[\S 5]{RZ-GSp4} freely.
 By \cite[Proposition 5.11, Proposition 5.18]{RZ-GSp4}, we have a $G$-equivariant isomorphism
 $H^i_c(M_\infty)\cong H^i_c(\M^\red_\infty,\mathcal{F}^{[0]})_{\overline{\Q}_\ell}$.
 On the other hand, \cite[Theorem 3.1]{MR1395723} tells us that 
 $H^i_{\mathcal{C}_{\M}}(M_\infty)\cong H^i_c(\M^\red_\infty,\mathcal{F}^{[2]})_{\overline{\Q}_\ell}$.
 The right hand side $H^i_c(\M^\red_\infty,\mathcal{F}^{[2]})_{\overline{\Q}_\ell}$ is endowed with an action of $G$
 (\cf \cite[\S 5.2]{RZ-GSp4}), and Lemma \ref{lem:action-compare} ensures that the isomorphism above
 is $G$-equivariant.
 By \cite[Proposition 5.6 i)]{RZ-GSp4}, for each $h\in \{1,2\}$ we have an exact sequence of
 smooth $G$-representations
 \[
 H^{i-1}_c(\M^\red_\infty,\mathcal{F}^{(h)})_{\overline{\Q}_\ell}\to
 H^i_c(\M^\red_\infty,\mathcal{F}^{[h-1]})_{\overline{\Q}_\ell}\to
 H^i_c(\M^\red_\infty,\mathcal{F}^{[h]})_{\overline{\Q}_\ell}\to
 H^i_c(\M^\red_\infty,\mathcal{F}^{(h)})_{\overline{\Q}_\ell}.
 \] 
 Hence the kernel and the cokernel of 
 $H^i_c(\M^\red_\infty,\mathcal{F}^{[0]})_{\overline{\Q}_\ell}\longrightarrow H^i_c(\M^\red_\infty,\mathcal{F}^{[2]})_{\overline{\Q}_\ell}$ is a successive extension of subquotients of
 $H^j_c(\M^\red_\infty,\mathcal{F}^{(h)})_{\overline{\Q}_\ell}$ for $i-1\le j\le i$ and $1\le h\le 2$.
 On the other hand, by \cite[Theorem 5.21]{RZ-GSp4},
 $H^j_c(\M^\red_\infty,\mathcal{F}^{(h)})_{\overline{\Q}_\ell}$ has no supercuspidal subquotient.
 Therefore the kernel and the cokernel of $H^i_c(M_\infty)\longrightarrow H^i_{\mathcal{C}_\M}(M_\infty)$ have no supercuspidal subquotient.
 Similar argument can be applied to the map 
 $H^i_c(M_\infty/p^\Z)\longrightarrow H^i_{\mathcal{C}_\M}(M_\infty/p^\Z)$.
\end{prf}

The following result on the vanishing of cohomology is also essentially obtained
in \cite{RZ-GSp4}.

\begin{prop}\label{prop:coh-vanishing}
 \begin{enumerate}
  \item For an integer $i<\dim M-\dim\M^\red$, we have $H^i_c(M_\infty)=0$.
  \item For an integer $i>\dim M+\dim\M^\red$, we have $H^i_{\mathcal{C}_\M}(M_\infty)=0$.
 \end{enumerate}
\end{prop}

\begin{prf}
 Fix an integer $m\ge 0$. Let $\mathcal{U}$ be a quasi-compact open formal subscheme of $\M_m$.
 By the $p$-adic uniformization theorem, there exist a scheme $U$ which is separated of finite type over $\Z_{p^\infty}$ and a closed subscheme $Z$ of the special fiber $U_s$ of $U$ such that
 $\mathcal{U}$ is isomorphic to the formal completion of $U$ along $Z$
 (\cf \cite[Corollary 4.4]{RZ-GSp4}; we can take $U$ as an open subscheme of a suitable integral model of the Siegel modular variety).
 We denote the closed immersion $Z\hooklongrightarrow U_s$ by $\iota$.
 By \cite[Theorem 4.35]{formalnearby} and \cite[Proposition 3.13]{formalnearby}, we have
 \[
  H^i_c\bigl(t(\mathcal{U})_{\overline{\eta}},\overline{\Q}_\ell\bigr)\cong H^i_c(Z,R\Psi_{\mathcal{U},c}\overline{\Q}_\ell)\cong H^i_c(Z,R\iota^!R\psi_U\overline{\Q}_\ell).
 \]
 Therefore, \cite[Lemma 5.26]{RZ-GSp4} tells us that 
 $H^i_c(t(\mathcal{U})_{\overline{\eta}},\overline{\Q}_\ell)=0$ if $i<\dim U_\eta-\dim Z=\dim M-\dim\M^\red$. Since
 $ H^i_c(M_\infty)=\varinjlim_m\varinjlim_{\mathcal{U}\subset\M_m}H^i_c(t(\mathcal{U})_{\overline{\eta}},\overline{\Q}_\ell)$, we conclude i).

 We prove ii). By \cite[Theorem 3.1]{MR1395723} and \cite[Lemma 5.26]{RZ-GSp4}, we have
 \[
  H^i_c(\mathcal{U}^\red,R\Psi\overline{\Q}_\ell)\cong H^i_c(Z,\iota^*R\psi_U\overline{\Q}_\ell)=0
 \]
 for $i>\dim M-\dim\M^\red$. Therefore, by Lemma \ref{lem:formal-support} and Proposition \ref{prop:H_c-limit} ii), we have
 \[
  H^i_{\mathcal{C}_\M}(M_{K_m})\cong H^i_c(\M_m^\red,R\Psi\overline{\Q}_\ell)\cong\varinjlim_{\mathcal{U}\subset \M_m}H^i_c(\mathcal{U}^\red,R\Psi\overline{\Q}_\ell)=0
 \]
 for $i>\dim M-\dim\M^\red$. 
 Hence $H^i_{\mathcal{C}_\M}(M_\infty)=\varinjlim_mH^i_{\mathcal{C}_\M}(M_{K_m})$ also vanishes
 for $i>\dim M-\dim\M^\red$. 
\end{prf}

\subsection{Application of the duality theorem}

Fix an isomorphism $\overline{\Q}_\ell\cong \C$ and identify them.
Every representation in this subsection is considered over $\C$.
Let $\widetilde{K}$ be an open compact-mod-center subgroup of $G$
and $\tau$ an irreducible smooth representation of $\widetilde{K}$.
Denote by $\chi\colon Z_G\longrightarrow \C^\times$ the central character of $\tau^\vee$.
For a smooth $G$-representation $V$, put $V_\tau=\Hom_{\widetilde{K}}(\tau,V\otimes_{\mathcal{H}(Z_G)}\chi^{-1})$.
Then, $H^i_c(M_\infty)_\tau$ and $H^i_{\mathcal{C}_\M}(M_\infty)_\tau$ are
representations of $J\times W_{\Q_p}$.
By \cite[Lemma 5.36]{MR1393439}, the actions of the center $Z_J=Z_G$ of $J$
on $H^i_c(M_\infty)_\tau$ and $H^i_{\mathcal{C}_\M}(M_\infty)_\tau$ are given by $\chi$.
Hence we can consider the Bernstein decomposition with respect to the action of $J$
(\cf Section \ref{sec:Zelevinsky}):
\[
H^i_c(M_\infty)_\tau=\bigoplus_{\mathfrak{s}\in I_\chi}H^i_c(M_\infty)_{\tau,\mathfrak{s}},\qquad
 H^i_{\mathcal{C}_\M}(M_\infty)_\tau=\bigoplus_{\mathfrak{s}\in I_\chi}H^i_{\mathcal{C}_\M}(M_\infty)_{\tau,\mathfrak{s}}.
\]

\begin{thm}\label{thm:RZ-Zel}
 Fix $\mathfrak{s}\in I_\chi$.
 Assume that $H^q_c(M_\infty)_{\tau,\mathfrak{s}}$ is a finite length $J$-representation
 for every integer $q$.
 Then, for each integer $i$, we have an isomorphism of $J\times W_{\Q_p}$-representations
 \[
  H^{2d+\iota(\mathfrak{s})-i}_{\mathcal{C}_\M}(M_\infty)_{\tau^\vee,\mathfrak{s}^\vee}(d)\cong
 \Zel\bigl(H^i_c(M_\infty)_{\tau,\mathfrak{s}}^\vee\bigr),
 \]
 where $d=n(n+1)/2$ is the dimension of $M$.
\end{thm}

\begin{rem}
 For the Drinfeld tower, a similar result is proved by Fargues \cite[Th\'eor\`eme 4.6]{Fargues-Zel}.
 In that case $\M$ is a $p$-adic formal scheme, thus $H^i_{\mathcal{C}_\M}(M_\infty)$ coincides
 with $H^i_c(M_\infty)$ (\cf Remark \ref{rem:p-adic-case}).
\end{rem}

\begin{prf}[of Theorem \ref{thm:RZ-Zel}]
 First we will reduce the theorem to the case where $\tau$ is trivial on $p^\Z$.
 Take $c\in \C$ such that $c^2=\chi(p)$ and define the character $\omega\colon \Q_p^\times\longrightarrow \C^\times$ by $\omega(a)=c^{v_p(a)}$, where $v_p$ is the $p$-adic valuation.
 We denote the composite of the similitude character
 $G\longrightarrow \Q_p^\times$ (resp.\ $J\longrightarrow \Q_p^\times$) and $\omega$
 by $\omega_G$ (resp.\ $\omega_J$). Then, $\tau'=\tau\otimes(\omega_G\vert_{\widetilde{K}})$
 is trivial on $p^\Z$.
 Moreover, as in the proof of \cite[Lemma 3.5]{RZ-LTF}, we have
 $H^q_c(M_\infty)\otimes\omega_G\otimes\omega_J^{-1}\cong H^q_c(M_\infty)$ as
 $G\times J\times W_{\Q_p}$-representations for every integer $q$. Therefore we obtain
 \begin{align*}
  H^q_c(M_\infty)_\tau&=\Hom_{\widetilde{K}}\bigl(\tau,H^q_c(M_\infty)\otimes_{\mathcal{H}(Z_G)}\chi^{-1}\bigr)\\
  &=\Hom_{\widetilde{K}}\bigl(\tau',(H^q_c(M_\infty)\otimes_{\mathcal{H}(Z_G)}\chi^{-1})\otimes\omega_G\bigr)\\
  &=\Hom_{\widetilde{K}}\bigl(\tau',(H^q_c(M_\infty)\otimes\omega_G)\otimes_{\mathcal{H}(Z_G)}\omega^2\chi^{-1}\bigr)\\
  &\cong \Hom_{\widetilde{K}}\bigl(\tau',(H^q_c(M_\infty)\otimes\omega_J)\otimes_{\mathcal{H}(Z_G)}\omega^2\chi^{-1}\bigr)
  =H^q_c(M_\infty)_{\tau'}\otimes\omega_J.
 \end{align*}
 For $\mathfrak{s}=[(\mathbf{M},\sigma)]\in I_\chi$, put $\mathfrak{s}'=[(\mathbf{M},\sigma\otimes (\omega_J^{-1}\vert_{\mathbf{M}(\Q_p)}))]$. Then, we conclude that
 $H^q_c(M_\infty)_{\tau,\mathfrak{s}}\cong H^q_c(M_\infty)_{\tau',\mathfrak{s}'}\otimes\omega_J$
 as $J\times W_{\Q_p}$-representations. In particular, the $J$-representation
 $H^q_c(M_\infty)_{\tau',\mathfrak{s}'}$ has finite length.
 Similarly we have $H^q_{\mathcal{C}_\M}(M_\infty)_{\tau^\vee,\mathfrak{s}^\vee}\cong H^q_{\mathcal{C}_\M}(M_\infty)_{\tau'^\vee,\mathfrak{s}'^\vee}\otimes\omega_J^{-1}$.
 Suppose that the theorem holds for $\tau'$ and $\mathfrak{s}'$. Then, by the isomorphisms above,
 we have
 \begin{align*}
 H^{2d+\iota(\mathfrak{s})-i}_{\mathcal{C}_\M}(M_\infty)_{\tau^\vee,\mathfrak{s}^\vee}(d)&\cong
 H^{2d+\iota(\mathfrak{s}')-i}_{\mathcal{C}_\M}(M_\infty)_{\tau'^\vee,\mathfrak{s}'^\vee}(d)\otimes\omega_J^{-1}\cong \Zel\bigl(H^i_c(M_\infty)_{\tau',\mathfrak{s}'}^\vee\bigr)\otimes\omega_J^{-1}\\
  &\cong \Zel\bigl(H^i_c(M_\infty)_{\tau,\mathfrak{s}}^\vee\otimes\omega_J\bigr)\otimes\omega_J^{-1}\stackrel{(*)}{\cong} \Zel\bigl(H^i_c(M_\infty)_{\tau,\mathfrak{s}}^\vee\bigr).
 \end{align*}
 For the isomorphism $(*)$, see Lemma \ref{lem:Zel-twist} i).

 Thus, in the following we may assume that $\tau$ is trivial on $p^\Z$.
 Then, we have
 \begin{gather*}
  H^q_c(M_\infty)_{\tau,\mathfrak{s}}=H^q_c(M_\infty/p^\Z)_{\tau,\mathfrak{s}}=\Hom_{\widetilde{K}}\bigl(\tau,H^q_c(M_\infty/p^\Z)_\mathfrak{s}\bigr),\\
  H^q_{\mathcal{C}_\M}(M_\infty)_{\tau^\vee,\mathfrak{s}^\vee}=H^q_{\mathcal{C}_\M}(M_\infty/p^\Z)_{\tau^\vee,\mathfrak{s}^\vee}=\Hom_{\widetilde{K}}\bigl(\tau^\vee,H^q_{\mathcal{C}_\M}(M_\infty/p^\Z)_{\mathfrak{s}^\vee}\bigr).
 \end{gather*}
 Since $\widetilde{K}$ is compact-mod-center, there exists a self-dual chain of lattices $\mathscr{L}$
 of $\Q_p^{2d}$ (\cf \cite[Definition 3.1, Definition 3.13]{MR1393439}) such that
 every $g\in \widetilde{K}$ and $L\in\mathscr{L}$ satisfy $gL\in \mathscr{L}$.
 As in \cite[\S 4.1]{RZ-LTF}, we write $K_{\!\mathscr{L}}$ for the stabilizer of $\mathscr{L}$ in $G$.
 For an integer $m\ge 0$, we put
 \[
  K_{\!\mathscr{L},m}=\{g\in K_{\!\mathscr{L}}\mid \text{for every $L\in\mathscr{L}$, $g$ acts trivially on $L/p^mL$}\}.
 \]
 It is an open normal subgroup of $K_{\!\mathscr{L}}$.
 We denote by $N_{\!\mathscr{L}}$ the subgroup of $G$ consisting of $g\in G$ satisfying
 $g\mathscr{L}=\mathscr{L}$.
 We have $K_{\!\mathscr{L},m}\subset K_{\!\mathscr{L}}\subset N_{\!\mathscr{L}}$ and
 $K_{\!\mathscr{L},m}$ is normal in $N_{\!\mathscr{L}}$.
 By definition, $\widetilde{K}$ is contained in $N_{\!\mathscr{L}}$.

 Take an integer $m\ge 0$ large enough so that  $K_{\!\mathscr{L},m}\subset K_0$ and
 $\tau\vert_{K_{\!\mathscr{L},m}}$ is trivial.
 In \cite[Definition 4.3]{RZ-LTF}, the author constructed a formal scheme $\M^\flat_{\!\mathscr{L},m}$
 over $\Spf \Z_{p^\infty}$ satisfying the following:
\begin{itemize}
 \item $t(\M^\flat_{\!\mathscr{L},m})_\eta$ is naturally isomorphic to $M_{K_{\!\mathscr{L}},m}$.
 \item $\M^\flat_{\!\mathscr{L},m}$ is naturally endowed with an action of $\widetilde{K}\times J$
       and a Weil descent datum, and they are compatible with those structures 
       on $M_{K_{\!\mathscr{L}},m}$ under the isomorphism above.
\end{itemize}
 We shall apply Theorem \ref{thm:duality-thm-l-adic} to $\M^\flat_{\!\mathscr{L},m}/p^\Z$ and $J/p^\Z$.
 We should verify the conditions in Theorem \ref{thm:duality-thm}. The conditions (a) and (c)
 are satisfied, as explained in the previous subsection. The condition (b) is satisfied by
 \cite[Remark 4.12]{RZ-LTF}. For the condition (d), let $\mathcal{I}$ denote the set of
 irreducible components of $\M^{\flat,\red}_{\!\mathscr{L},m}/p^\Z$. 
 For $\alpha\in \mathcal{I}$, put $V_\alpha=(\M^{\flat,\red}_{\!\mathscr{L},m}/p^\Z)\setminus \bigcup_{\beta\in\mathcal{I},\alpha\cap \beta=\varnothing}\beta$. It is a quasi-compact open subset
 of $\M^{\flat,\red}_{\!\mathscr{L},m}/p^\Z$. By \cite[Lemma 5.1 ii)]{RZ-LTF}, the action of $J/p^\Z$
 on $\mathcal{I}$ has finite orbits. Take $\alpha_1,\ldots,\alpha_k\in \mathcal{I}$ so that
 $\mathcal{I}=\bigcup_{j=1}^k (J/p^\Z)\alpha_j$, and put $V=\bigcup_{j=1}^k V_{\alpha_j}$.
 Clearly we have $\M^{\flat,\red}_{\!\mathscr{L},m}/p^\Z=\bigcup_{h\in J/p^\Z}hV$.
 The closure $\overline{V}$ is the union of finitely many irreducible components
 of $\M^{\flat,\red}_{\!\mathscr{L},m}/p^\Z$. Therefore, by the same way as in
 \cite[Corollary 4.3 ii)]{RZ-irr-comp}, we can prove that 
 the set $\{g\in J/p^\Z\mid g\overline{V}\cap \overline{V}\neq\varnothing\}$ is compact.
 Thus the condition (d) is satisfied. Now, by Theorem \ref{thm:duality-thm-l-adic} and
 Remark \ref{rem:coh-finite-type}, we have a $J$-equivariant isomorphism
 \[
  H^{2d+i}_{c,\M^\flat_{\!\mathscr{L},m}}(M_{K_{\!\mathscr{L},m}}/p^\Z)(d)\cong R^i\!\D\bigl(R\Gamma_c\bigl((M_{K_{\!\mathscr{L},m}}/p^\Z)/(J/p^\Z),\overline{\Q}_\ell\bigr)\bigr).
 \]
 By Corollary \ref{cor:functoriality},
 this isomorphism is also $\widetilde{K}\times W_{\Q_p}$-equivariant.

 We prove that the left hand side is equal to 
 $H^{2d+i}_{\mathcal{C}_\M}(M_{K_{\!\mathscr{L},m}}/p^\Z)(d)$.
 For simplicity, we denote the support set $\mathcal{C}_{\M^\flat_{\!\mathscr{L},m}}$ of
 $M_{K_{\!\mathscr{L},m}}$ by $\mathcal{C}_{\mathscr{L},m}$.
 Let $\mathscr{L}_{\mathrm{Iw}}$ be the self-dual chain of lattices
 \[
 \bigl\{(p^m\Z_p)^{\oplus j}\oplus (p^{m+1}\Z_p)^{\oplus (2d-j)}\bigr\}_{0\le j\le 2d,m\in\Z}.
 \]
 The group $K_{\!\mathscr{L}_{\mathrm{Iw}}}$ attached to $\mathscr{L}_{\mathrm{Iw}}$ is
 an Iwahori subgroup of $G$.
 There exists $g_0\in G$ such that $\mathscr{L}\subset g_0\mathscr{L}_{\mathrm{Iw}}$.
 The following morphisms of formal schemes are naturally induced:
 \[
  \M^\flat_{\!\mathscr{L},m}\yleftarrow{(1)}\M^\flat_{\!g_0\mathscr{L}_{\mathrm{Iw}},m}\yrightarrow[\cong]{[g_0]}
  \M^\flat_{\!\mathscr{L}_{\mathrm{Iw}},m}\yrightarrow{(2)}\M_m\yrightarrow{(3)}\M.
 \]
 The morphisms (1), (2) and (3) are proper.
 The rigid generic fiber of the diagram above is identified with the following diagram:
 \[
  M_{K_{\!\mathscr{L},m}}\yleftarrow{\pi} M_{g_0K_{\!\mathscr{L}_{\mathrm{Iw}},m}g_0^{-1}}\yrightarrow[\cong]{[g_0]}
  M_{K_{\!\mathscr{L}_{\mathrm{Iw}},m}}\longrightarrow M_{K_m}\longrightarrow M.
 \]
 Therefore, by Lemma \ref{lem:formal-finite-type} and Proposition \ref{prop:Hecke-support-set},
 we have equalities 
 \[
  \pi^{-1}\mathcal{C}_{\mathscr{L},m}=\mathcal{C}_{g_0\mathscr{L}_{\mathrm{Iw}},m}=[g_0]^{-1}\mathcal{C}_{\mathscr{L}_{\mathrm{Iw}},m}=[g_0]^{-1}\mathcal{C}_\M=\mathcal{C}_\M=\pi^{-1}\mathcal{C}_\M
 \]
 of support sets of $M_{g_0K_{\!\mathscr{L}_{\mathrm{Iw}},m}g_0^{-1}}$
 (recall that we denote by $\mathcal{C}_\M$ the support set of $M_K$ induced from the support set 
 $\mathcal{C}_\M$ of $M$ for various $K\subset K_0$).
 Since $\pi$ is finite and surjective, we conclude that $\mathcal{C}_{\mathscr{L},m}=\mathcal{C}_\M$.
 Hence $H^{2d+i}_{c,\M^\flat_{\!\mathscr{L},m}}(M_{K_{\!\mathscr{L},m}}/p^\Z)(d)=H^{2d+i}_{\mathcal{C}_\M}(M_{K_{\!\mathscr{L},m}}/p^\Z)(d)$, and thus we have a $\widetilde{K}\times J\times W_{\Q_p}$-equivariant isomorphism
 \[
  H^{2d+i}_{\mathcal{C}_\M}(M_{K_{\!\mathscr{L},m}}/p^\Z)(d)\cong R^i\!\D\bigl(R\Gamma_c\bigl((M_{K_{\!\mathscr{L},m}}/p^\Z)/(J/p^\Z),\overline{\Q}_\ell\bigr)\bigr).
 \]
 The $(\tau^\vee,\mathfrak{s}^\vee)$-part of the left hand side is equal to 
 $H^{2d+i}_{\mathcal{C}_\M}(M_\infty)_{\tau^\vee,\mathfrak{s}^\vee}(d)$.
 We will consider the $(\tau^\vee,\mathfrak{s}^\vee)$-part of the right hand side.
 For simplicity, we put 
 \[
   A=R\Gamma_c\bigl((M_{K_{\!\mathscr{L},m}}/p^\Z)/(J/p^\Z),\overline{\Q}_\ell\bigr),
 \]
 which is an object of $D^b(\Rep(J/p^\Z))$ endowed with an action of $\widetilde{K}\times W_{\Q_p}$.
 The action of $\widetilde{K}$ factors through the finite quotient 
 $H=\widetilde{K}/K_{\!\mathscr{L},m}p^\Z$.
 By Corollary \ref{cor:derived-same}, $H^q(A)=H^q_c(M_{K_{\!\mathscr{L},m}}/p^\Z)$.
 Therefore, we have a spectral sequence
 \[
  E_2^{s,t}= R^s\!\D\bigl(H^{-t}_c(M_{K_{\!\mathscr{L},m}}/p^\Z)\bigr)\Longrightarrow R^{s+t}\!\D(A).
 \]
 Take the $(\tau^\vee,\mathfrak{s}^\vee)$-part of this spectral sequence:
 \[
  E_2^{s,t}= R^s\!\D\bigl(H^{-t}_c(M_{K_{\!\mathscr{L},m}}/p^\Z)\bigr)_{\tau^\vee,\mathfrak{s}^\vee}\Longrightarrow R^{s+t}\!\D(A)_{\tau^\vee,\mathfrak{s}^\vee}.
 \]
 We can observe that 
 \[
  R^s\!\D\bigl(H^{-t}_c(M_{K_{\!\mathscr{L},m}}/p^\Z)\bigr)_{\tau^\vee,\mathfrak{s}^\vee}
 =R^s\!\D\bigl(H^{-t}_c(M_{K_{\!\mathscr{L},m}}/p^\Z)_{\tau,\mathfrak{s}}\bigr)
 =R^s\!\D\bigl(H^{-t}_c(M_\infty)_{\tau,\mathfrak{s}}\bigr).
 \]
 Indeed, for the $\tau^\vee$-part, notice the isotypic decomposition 
 \[
 H^{-t}_c(M_{K_{\!\mathscr{L},m}}/p^\Z)=\bigoplus_{\sigma}H^{-t}_c(M_{K_{\!\mathscr{L},m}}/p^\Z)_\sigma\otimes\sigma,
 \]
 where $\sigma$ runs through irreducible representations of $H$. Since each $\sigma$ is
 finite-dimensional and there are only finitely many such $\sigma$, we have
 \[
  R^s\!\D\bigl(H^{-t}_c(M_{K_{\!\mathscr{L},m}}/p^\Z)\bigr)=\bigoplus_{\sigma}R^s\!\D\bigl(H^{-t}_c(M_{K_{\!\mathscr{L},m}}/p^\Z)_\sigma\bigr)\otimes\sigma^\vee,
 \]
 and thus $R^s\!\D(H^{-t}_c(M_{K_{\!\mathscr{L},m}}/p^\Z))_{\tau^\vee}=R^s\!\D(H^{-t}_c(M_{K_{\!\mathscr{L},m}}/p^\Z)_\tau)$. For the $\mathfrak{s}$-part, see \cite[Remarque 1.5]{Fargues-Zel}.

 By the assumption, $H^{-t}_c(M_\infty)_{\tau,\mathfrak{s}}$ has finite length as a $J$-representation. 
 Therefore, Theorem \ref{thm:Schneider-Stuhler} tells us that 
 $R^s\!\D\bigl(H^{-t}_c(M_{K_{\!\mathscr{L},m}}/p^\Z)\bigr)_{\tau^\vee,\mathfrak{s}^\vee}=0$ unless
 $s=\iota(\mathfrak{s})$. Thus, we have
 \[
  R^i\!\D(A)_{\tau^\vee,\mathfrak{s}^\vee}\cong R^{\iota(\mathfrak{s})}\!\D\bigl(H^{-i+\iota(\mathfrak{s})}_c(M_\infty)_{\tau,\mathfrak{s}}\bigr)=\Zel\bigl(H^{-i+\iota(\mathfrak{s})}_c(M_\infty)^\vee_{\tau,\mathfrak{s}}\bigr).
 \]
 Hence we obtain a $J\times W_{\Q_p}$-equivariant isomorphism
 \[
  H^{2d+i}_{\mathcal{C}_\M}(M_\infty)_{\tau^\vee,\mathfrak{s}^\vee}(d)\cong \Zel\bigl(H^{-i+\iota(\mathfrak{s})}_c(M_\infty)^\vee_{\tau,\mathfrak{s}}\bigr).
 \]
 Replacing $i$ by $\iota(\mathfrak{s})-i$, we conclude the theorem.
\end{prf}

To apply Theorem \ref{thm:RZ-Zel}, we need the following technical assumption.

\begin{assump}\label{assump:finiteness}
 For each integer $q$ and each compact open subgroup $H$ of $J$,
 the $G$-representation $H^q_c(M_\infty)^H$ is finitely generated.
\end{assump}

\begin{rem}
 We know that for every compact open subgroup $K$ of $G$, the $J$-representation
 $H^q_c(M_\infty)^K$ is finitely generated (\cf \cite[Proposition 4.4.13]{MR2074714};
 see also \cite[Theorem 4.4]{RZ-irr-comp}).
 Therefore, if we can establish an analogue of Faltings' isomorphism for our tower $\{M_K\}_K$,
 we can prove Assumption \ref{assump:finiteness} by switching the role of $G$ and $J$.
\end{rem}

\begin{lem}\label{lem:rep-theory}
 \begin{enumerate}
  \item For an irreducible supercuspidal representation $\pi$ of $G$, there exists
	a compact-mod-center open subgroup $\widetilde{K}$ of $G$ and an irreducible smooth
	representation $\tau$ of $\widetilde{K}$ such that $\cInd_{\widetilde{K}}^G\tau$ is
	admissible (hence supercuspidal) and $\pi$ is a direct summand of $\cInd_{\widetilde{K}}^G\tau$.
  \item If $n=2$, under the same setting as i), we can take $(\widetilde{K},\tau)$ so that
	$\pi=\cInd_{\widetilde{K}}^G\tau$.
  \item Let $(\widetilde{K},\tau)$ be as in i). Under Assumption \ref{assump:finiteness},
	$H^q_c(M_\infty)_\tau$ is a finite length $J$-representation for every integer $q$.
 \end{enumerate}
\end{lem}

\begin{prf}
 i) Put $G_1=\Sp_{2n}(\Q_p)$ and write $G'$ for the image of $\Q_p^\times\times G_1\longrightarrow G$;
 $(c,g)\longmapsto cg$. Let $\pi'$ be the restriction of $\pi$ to $G'$.
 Since $G'$ is a normal subgroup of index $2$ in $G$,
 $\pi'$ has finite length. 
 Take an irreducible subrepresentation
 $\pi_1$ of $\pi'$, and regard it as a representation of $\Q_p^\times\times G_1$.
 It is irreducible and supercuspidal.
 By \cite[Theorem 7.14]{MR2390287}, there exist a compact open subgroup $K$ of $G_1$,
 an irreducible smooth representation $\sigma$ of $K$ and
 a smooth character $\chi$ of $\Q_p^\times$ such that $\pi_1\cong \chi\otimes\cInd_K^{G_1}\sigma$.
 Let $\widetilde{K}$ be the image of $\Q_p^\times\times K$ in $G_1$. Then, $\chi\otimes\sigma$
 descends to an irreducible smooth representation $\tau$ of $\widetilde{K}$
 and we have $\cInd_{\widetilde{K}}^{G'}\tau\cong \pi_1$.
 Then, $\cInd_{\widetilde{K}}^G\tau\cong \cInd_{G'}^G\pi_1$ is admissible.
 If $\pi_1\neq \pi'$, $\cInd_{\widetilde{K}}^G\tau$ is
 isomorphic to $\pi$; otherwise $\cInd_{\widetilde{K}}^G\tau$ is the direct sum of two supercuspidal
 representations, one of which is isomorphic to $\pi$.

 ii) is proved in \cite{MR948308}; see the final comment in \cite[p.~328]{MR948308}.
 
 iii) Let $\chi$ be the central character of $\tau^\vee$. By the Frobenius reciprocity, we have
 $H^q_c(M_\infty)_\tau\cong \Hom_G(\cInd_{\widetilde{K}}^G\tau,H^q_c(M_\infty)\otimes_{\mathcal{H}(Z_G)}\chi^{-1})$.
 As $\cInd_{\widetilde{K}}^G\tau$ is a supercuspidal representation of finite length 
 with central character $\chi^{-1}$,
 it is a finite direct sum of irreducible supercuspidal representations of $G$.
 Therefore, it suffices to show that for an irreducible supercuspidal representation $\pi''$ of $G$
 with central character $\chi^{-1}$,
 $\Hom_G(\pi'',H^q_c(M_\infty)\otimes_{\mathcal{H}(Z_G)}\chi^{-1})$
 is a $J$-representation of finite length.
 We apply \cite[Lemma 5.2]{LT-LTF} to $H^q_c(M_\infty)\otimes_{\mathcal{H}(Z_G)}\chi^{-1}$.
 By \cite[Proposition 4.4.13]{MR2074714} and Assumption \ref{assump:finiteness},
 all the conditions in \cite[Lemma 5.2]{LT-LTF} are satisfied.
 Hence we conclude that $\Hom_G(\pi'',H^q_c(M_\infty)\otimes_{\mathcal{H}(Z_G)}\chi^{-1})$ has
 finite length (see the final paragraph of the proof of \cite[Lemma 5.2]{LT-LTF}).
\end{prf}

Now we will give some consequences of Theorem \ref{thm:RZ-Zel}.
For simplicity, we focus on the cohomology $H^i_c(M_\infty/p^\Z)$.

\begin{cor}\label{cor:vanishing-noncusp-part}
 Suppose Assumption \ref{assump:finiteness}.
 Let $\pi$ be an irreducible supercuspidal representation of $G$
 and $\rho$ an irreducible non-supercuspidal representation of $J$.
 Then, the representation $\pi\otimes \rho$ of $G\times J$ does not appear as a subquotient of
 $H^{d-\dim\M^\red}_c(M_\infty/p^\Z)$.
\end{cor}

\begin{prf}
 Clearly we may assume that $\pi$ (resp.\ $\rho$) is trivial on $p^\Z\subset G$
 (resp.\ $p^\Z\subset J$).
 Take $(\widetilde{K},\tau)$ as in Lemma \ref{lem:rep-theory} i). Then $\tau$ is trivial on
 $p^\Z\subset G$.
 Let $\mathfrak{s}\in I_{p^\Z}$ be the inertially equivalence class of cuspidal data for $J$
 such that $\rho\in \Rep(J/p^\Z)_{\mathfrak{s}}$. As $\rho$ is non-supercuspidal,
 we have $\iota(\mathfrak{s})\ge 1$.
 
 Put $d_0=\dim\M^\red$. 
 As $\pi$ is projective in the category $\Rep(G/p^\Z)$,
 $\pi\otimes\rho$ appears in $H^{d-d_0}_c(M_\infty/p^\Z)$
 if and only if $\rho$ appears in $\Hom_G(\pi,H^{d-d_0}_c(M_\infty/p^\Z))$. 
 Since $\Hom_G(\pi,H^{d-d_0}_c(M_\infty/p^\Z))$ is embedded into
 \[
 \Hom_G\bigl(\cInd_{\widetilde{K}}^G\tau,H^{d-d_0}_c(M_\infty/p^\Z)\bigr)
 =H^{d-d_0}_c(M_\infty)_\tau,
 \]
 it suffices to show that $H^{d-d_0}_c(M_\infty)_{\tau,\mathfrak{s}}=0$.
 By Theorem \ref{thm:RZ-Zel} and Lemma \ref{lem:rep-theory} iii), we have
 $\Zel(H^{d-d_0}_c(M_\infty)^\vee_{\tau,\mathfrak{s}})\cong H^{d+d_0+\iota(\mathfrak{s})}_{\mathcal{C}_\M}(M_\infty)_{\tau^\vee,\mathfrak{s}^\vee}(d)$. By Proposition \ref{prop:coh-vanishing} ii), 
 the right hand side equals $0$. Hence we conclude that $H^{d-d_0}_c(M_\infty)_{\tau,\mathfrak{s}}=0$.
\end{prf}

\begin{cor}\label{cor:A-packet}
 In addition to Assumption \ref{assump:finiteness}, assume that $n=2$.
 Let $\pi$ and $\rho$ be as in Corollary \ref{cor:vanishing-noncusp-part},
 and $\sigma$ an irreducible $\ell$-adic representation of $W_{\Q_p}$.
 Then, $\pi\otimes \rho\otimes\sigma$ appears as a subquotient of
 $H^3_c(M_\infty/p^\Z)$ if and only if
 $\pi^\vee\otimes \Zel(\rho^\vee)\otimes\sigma^\vee(-3)$ appears as a subquotient of 
 $H^4_c(M_\infty/p^\Z)$.
\end{cor}

\begin{prf}
 Again we may assume that $\pi$ (resp.\ $\rho$) is trivial on $p^\Z\subset G$ (resp.\ $p^\Z\subset J$).
 Let $\mathfrak{s}\in I_{p^\Z}$ be as in the proof of Corollary \ref{cor:vanishing-noncusp-part}.
 In this case we have $\iota(\mathfrak{s})=1$. 
 Take $(\widetilde{K},\tau)$ as in Lemma \ref{lem:rep-theory} ii). 
 Note that, since $\pi=\cInd_{\widetilde{K}}^G\tau$ is irreducible, we have
 $\pi^\vee=(\cInd_{\widetilde{K}}^G\tau)^\vee=\Ind_{\widetilde{K}}^G\tau^\vee=\cInd_{\widetilde{K}}^G\tau^\vee$. 
 In the same way as in the proof of Corollary \ref{cor:vanishing-noncusp-part},
 we can see that $\pi\otimes\rho\otimes\sigma$ appears as a subquotient of $H^3_c(M_\infty/p^\Z)$
 if and only if
 $\rho\otimes\sigma$ appears as a subquotient of $H^3_c(M_\infty)_{\tau,\mathfrak{s}}$.
 Similarly, $\pi^\vee\otimes\Zel(\rho^\vee)\otimes\sigma^\vee(-3)$ appears as a subquotient of
 $H^4_c(M_\infty/p^\Z)$
 if and only if $\Zel(\rho^\vee)\otimes\sigma^\vee(-3)$ appears as a subquotient of $H^4_c(M_\infty)_{\tau^\vee,\mathfrak{s}^\vee}$. 

 On the other hand, by Theorem \ref{thm:RZ-Zel} and Lemma \ref{lem:rep-theory} iii), we have
 \[
  \Zel\bigl(H^3_c(M_\infty)^\vee_{\tau,\mathfrak{s}}\bigr)\cong H^4_{\mathcal{C}_\M}(M_\infty)_{\tau^\vee,\mathfrak{s}^\vee}(3).
 \]
 Corollary \ref{cor:compare-cuspidal-parts} tells us that
 \begin{align*}
  H^4_{\mathcal{C}_\M}(M_\infty)_{\tau^\vee}
  &=\Hom_{\widetilde{K}}\bigl(\tau^\vee,H^4_{\mathcal{C}_\M}(M_\infty/p^\Z)\bigr)
  =\Hom_G\bigl(\pi^\vee,H^4_{\mathcal{C}_\M}(M_\infty/p^\Z)\bigr)\\
  &\cong \Hom_G\bigl(\pi^\vee,H^4_c(M_\infty/p^\Z)\bigr)=H^4_c(M_\infty)_{\tau^\vee}.
 \end{align*}
 Thus we have an isomorphism of $J\times W_{\Q_p}$-representations
 \[
  \Zel\bigl(H^3_c(M_\infty)^\vee_{\tau,\mathfrak{s}}\bigr)\cong H^4_c(M_\infty)_{\tau^\vee,\mathfrak{s}^\vee}(3).
 \]
 In particular, $\rho\otimes\sigma$ appears in $H^3_c(M_\infty)_{\tau,\mathfrak{s}}$ if and only if
 $\Zel(\rho^\vee)\otimes\sigma^\vee(-3)$ appears in $H^4_c(M_\infty)_{\tau^\vee,\mathfrak{s}^\vee}$.
 This concludes the proof.
\end{prf}

\begin{rem}\label{rem:Ito-Mieda}
 \begin{enumerate}
  \item In the case $n=2$, by a global method, Tetsushi Ito and the author proved that
	for a supercuspidal representation $\pi$ of $G$ and
	a supercuspidal representation $\rho$ of $J$,
	$\pi\otimes \rho$ does not appear in $H^i_c(M_\infty/p^\Z)$ unless $i=3$.
	Together with Corollary \ref{cor:vanishing-noncusp-part}, we obtain the vanishing of
	the $G$-cuspidal part $H^2_c(M_\infty/p^\Z)_{\text{$G$-cusp}}=0$
	(note that in this case $\dim\M^\red=1$).
  \item Consider the case $n=2$.
	Let $\pi$ be an irreducible supercuspidal representation of $G$ which is trivial on $p^\Z$,
	$\phi$ the $L$-parameter of $\pi$ (\cf \cite{MR2800725}),
	and $\Pi_\phi^J$ the $L$-packet of $J$ attached to $\phi$ (\cf \cite{Gan-Tantono}).
	Suppose that there exists a non-supercuspidal
	representation $\rho$ in $\Pi^J_\phi$. Then, $\Pi^J_\phi$ consists of two representations
	$\rho$ and $\rho'$, where $\rho'$ is supercuspidal.
	In this case, $\{\Zel(\rho),\rho'\}$ is known to be a non-tempered $A$-packet.

	Motivated by Corollary \ref{cor:A-packet} and the results in \cite{RZ-LTF},
	the author expects the following:
	\begin{itemize}
	 \item $\pi\otimes \rho^\vee$ and $\pi\otimes \rho'^\vee$ appear in $H^3_c(M_\infty/p^\Z)$, and
	 \item $\pi\otimes \Zel(\rho)^\vee$ appears in $H^4_c(M_\infty/p^\Z)$.
	\end{itemize}
	This problem will be considered in a forthcoming joint paper with Tetsushi Ito.
	The speculation above suggests the existence of some relation between local $A$-packets and
	the cohomology of the Rapoport-Zink tower. In the case of $\GL(n)$, a result in this direction
	is obtained by Dat \cite{MR2904195}.
 \end{enumerate}
\end{rem}

\begin{rem}\label{rem:generalizations}
 \begin{enumerate}
  \item The arguments here can be applied to more general Rapoport-Zink towers.
	In \cite{RZ-irr-comp}, the geometric properties used in the proof of Theorem \ref{thm:RZ-Zel}
	are obtained for many Rapoport-Zink spaces.
	See \cite[Theorem 2.6, Corollary 4.3, Theorem 4.4]{RZ-irr-comp} especially.
  \item There is another way to generalize the result in \cite{Fargues-Zel}; the author expects
	some relation between $\Zel_{G\times J}(H^i_c(M_\infty)_{\mathfrak{s}}^\vee)$ and 
	$H^{2d+\iota(\mathfrak{s})-i}_c(M_\infty)_{\mathfrak{s}^\vee}(d)$, where
	$\mathfrak{s}$ is an inertially equivalence class of cuspidal data for $G\times J$.
	In a subsequent paper, the author plans to work on this problem in the case of $\GSp(4)$.
	However, the case where $n\ge 3$ seems much more difficult, and Theorem \ref{thm:RZ-Zel}
	has advantage that it is valid for all the $\GSp(2n)$ cases (and many other cases).
 \end{enumerate}
\end{rem}

\def\cftil#1{\ifmmode\setbox7\hbox{$\accent"5E#1$}\else
  \setbox7\hbox{\accent"5E#1}\penalty 10000\relax\fi\raise 1\ht7
  \hbox{\lower1.15ex\hbox to 1\wd7{\hss\accent"7E\hss}}\penalty 10000
  \hskip-1\wd7\penalty 10000\box7} \def\cprime{$'$} \def\cprime{$'$}
\providecommand{\bysame}{\leavevmode\hbox to3em{\hrulefill}\thinspace}
\providecommand{\MR}{\relax\ifhmode\unskip\space\fi MR }
\providecommand{\MRhref}[2]{%
  \href{http://www.ams.org/mathscinet-getitem?mr=#1}{#2}
}
\providecommand{\href}[2]{#2}


\begin{thebibliography}{Hub98b}

\bibitem[Ber84]{MR771671}
J.~N. Bernstein, \emph{Le ``centre'' de {B}ernstein}, Representations of
  reductive groups over a local field, Travaux en Cours, Hermann, Paris, 1984,
  Edited by P. Deligne, pp.~1--32.

\bibitem[Ber94]{MR1262943}
V.~G. Berkovich, \emph{Vanishing cycles for formal schemes}, Invent. Math.
  \textbf{115} (1994), no.~3, 539--571.

\bibitem[Ber96]{MR1395723}
\bysame, \emph{Vanishing cycles for formal schemes. {II}}, Invent. Math.
  \textbf{125} (1996), no.~2, 367--390.

\bibitem[Boy09]{MR2511742}
P.~Boyer, \emph{Monodromie du faisceau pervers des cycles \'evanescents de
  quelques vari\'et\'es de {S}himura simples}, Invent. Math. \textbf{177}
  (2009), no.~2, 239--280.

\bibitem[Car90]{MR1044827}
H.~Carayol, \emph{Nonabelian {L}ubin-{T}ate theory}, Automorphic forms,
  {S}himura varieties, and {$L$}-functions, {V}ol.\ {II} ({A}nn {A}rbor, {MI},
  1988), Perspect. Math., vol.~11, Academic Press, Boston, MA, 1990,
  pp.~15--39.

\bibitem[Dat12]{MR2904195}
J.-F. Dat, \emph{Op\'erateur de {L}efschetz sur les tours de {D}rinfeld et
  {L}ubin-{T}ate}, Compos. Math. \textbf{148} (2012), no.~2, 507--530.

\bibitem[Fal02]{MR1936369}
G.~Faltings, \emph{A relation between two moduli spaces studied by {V}. {G}.
  {D}rinfeld}, Algebraic number theory and algebraic geometry, Contemp. Math.,
  vol. 300, Amer. Math. Soc., Providence, RI, 2002, pp.~115--129.

\bibitem[Far04]{MR2074714}
L.~Fargues, \emph{Cohomologie des espaces de modules de groupes
  {$p$}-divisibles et correspondances de {L}anglands locales}, Ast\'erisque
  (2004), no.~291, 1--199, Vari{\'e}t{\'e}s de Shimura, espaces de
  Rapoport-Zink et correspondances de Langlands locales.

\bibitem[Far06]{Fargues-Zel}
\bysame, \emph{Dualit\'e de {P}oincar\'e et involution de {Z}elevinsky dans la
  cohomologie \'etale \'equivariante des espaces analytiques rigides},
  preprint, \url{http://www.math.jussieu.fr/~fargues/Prepublications.html},
  2006.

\bibitem[FGL08]{MR2441311}
L.~Fargues, A.~Genestier, and V.~Lafforgue, \emph{L'isomorphisme entre les
  tours de {L}ubin-{T}ate et de {D}rinfeld}, Progress in Mathematics, vol. 262,
  Birkh\"auser Verlag, Basel, 2008.

\bibitem[GT]{Gan-Tantono}
W.~T. Gan and W.~Tantono, \emph{The local {L}anglands conjecture for {${\rm
  GSp}(4)$} {II}: the case of inner forms}, preprint.

\bibitem[GT11]{MR2800725}
W.~T. Gan and S.~Takeda, \emph{The local {L}anglands conjecture for {${\rm
  GSp}(4)$}}, Ann. of Math. (2) \textbf{173} (2011), no.~3, 1841--1882.

\bibitem[Har97]{MR1464867}
M.~Harris, \emph{Supercuspidal representations in the cohomology of
  {D}rinfel\cprime d upper half spaces; elaboration of {C}arayol's program},
  Invent. Math. \textbf{129} (1997), no.~1, 75--119.

\bibitem[Hoc69]{MR0251026}
M.~Hochster, \emph{Prime ideal structure in commutative rings}, Trans. Amer.
  Math. Soc. \textbf{142} (1969), 43--60.

\bibitem[HT01]{MR1876802}
M.~Harris and R.~Taylor, \emph{The geometry and cohomology of some simple
  {S}himura varieties}, Annals of Mathematics Studies, vol. 151, Princeton
  University Press, Princeton, NJ, 2001, With an appendix by Vladimir G.
  Berkovich.

\bibitem[Hub93]{MR1207303}
R.~Huber, \emph{Continuous valuations}, Math. Z. \textbf{212} (1993), no.~3,
  455--477.

\bibitem[Hub94]{MR1306024}
\bysame, \emph{A generalization of formal schemes and rigid analytic
  varieties}, Math. Z. \textbf{217} (1994), no.~4, 513--551.

\bibitem[Hub96]{MR1734903}
\bysame, \emph{\'{E}tale cohomology of rigid analytic varieties and adic
  spaces}, Aspects of Mathematics, E30, Friedr. Vieweg \& Sohn, Braunschweig,
  1996.

\bibitem[Hub98a]{MR1620118}
\bysame, \emph{A finiteness result for direct image sheaves on the
  \'etale site of rigid analytic varieties}, J. Algebraic Geom. \textbf{7}
  (1998), no.~2, 359--403.

\bibitem[Hub98b]{MR1626021}
\bysame, \emph{A comparison theorem for {$l$}-adic cohomology},
  Compositio Math. \textbf{112} (1998), no.~2, 217--235.

\bibitem[IM10]{RZ-GSp4}
T.~Ito and Y.~Mieda, \emph{Cuspidal representations in the $\ell$-adic
  cohomology of the {R}apoport-{Z}ink space for $\mathrm{GSp}(4)$}, preprint,
  arXiv:1005.5619, 2010.

\bibitem[Jan88]{MR929536}
U.~Jannsen, \emph{Continuous \'etale cohomology}, Math. Ann. \textbf{280}
  (1988), no.~2, 207--245.

\bibitem[Man05]{MR2169874}
E.~Mantovan, \emph{On the cohomology of certain {PEL}-type {S}himura
  varieties}, Duke Math. J. \textbf{129} (2005), no.~3, 573--610.

\bibitem[Mie10]{formalnearby}
Y.~Mieda, \emph{Variants of formal nearby cycles}, preprint, arXiv:1005.5616,
  to appear in J. Inst. Math. Jussieu, 2010.

\bibitem[Mie12a]{LT-LTF}
\bysame, \emph{Lefschetz trace formula and $\ell$-adic cohomology of
  {L}ubin-{T}ate tower}, Math. Res. Lett. \textbf{19} (2012), no.~1, 95--107.

\bibitem[Mie12b]{RZ-LTF}
\bysame, \emph{Lefschetz trace formula and {$\ell$}-adic cohomology of
  {R}apoport-{Z}ink tower for {$\mathrm{GSp}(4)$}}, preprint, arXiv:1212.4922,
  2012.

\bibitem[Mie13]{RZ-irr-comp}
\bysame, \emph{On irreducible components of {R}apoport-{Z}ink spaces},
  preprint, arXiv:1312.4784, 2013.

\bibitem[Moy88]{MR948308}
A~Moy, \emph{Representations of {$G\,{\rm Sp}(4)$} over a {$p$}-adic field.
  {I}, {II}}, Compositio Math. \textbf{66} (1988), no.~3, 237--284, 285--328.

\bibitem[Rap95]{MR1403942}
M.~Rapoport, \emph{Non-{A}rchimedean period domains}, Proceedings of the
  {I}nternational {C}ongress of {M}athematicians, {V}ol.\ 1, 2 ({Z}\"urich,
  1994) (Basel), Birkh\"auser, 1995, pp.~423--434.

\bibitem[Ren10]{MR2567785}
D.~Renard, \emph{Repr\'esentations des groupes r\'eductifs {$p$}-adiques},
  Cours Sp\'ecialis\'es, vol.~17, Soci\'et\'e Math\'ematique de France, Paris,
  2010.

\bibitem[RZ96]{MR1393439}
M.~Rapoport and Th. Zink, \emph{Period spaces for {$p$}-divisible groups},
  Annals of Mathematics Studies, vol. 141, Princeton University Press,
  Princeton, NJ, 1996.

\bibitem[SS97]{MR1471867}
P.~Schneider and U.~Stuhler, \emph{Representation theory and sheaves on the
  {B}ruhat-{T}its building}, Inst. Hautes \'Etudes Sci. Publ. Math. (1997),
  no.~85, 97--191.

\bibitem[Ste08]{MR2390287}
S.~Stevens, \emph{The supercuspidal representations of {$p$}-adic classical
  groups}, Invent. Math. \textbf{172} (2008), no.~2, 289--352.

\bibitem[Vig90]{MR1159104}
M.-F. Vign{\'e}ras, \emph{On formal dimensions for reductive {$p$}-adic
  groups}, Festschrift in honor of {I}. {I}. {P}iatetski-{S}hapiro on the
  occasion of his sixtieth birthday, {P}art {I} ({R}amat {A}viv, 1989), Israel
  Math. Conf. Proc., vol.~2, Weizmann, Jerusalem, 1990, pp.~225--266.

\bibitem[EGA]{EGA}
A.~Grothendieck and J.~Dieudonn\'e, \emph{\'{E}l\'ements de g\'eom\'etrie
  alg\'ebrique}, Inst. Hautes \'Etudes Sci. Publ. Math. (1961--1967), no.~4, 8,
  11, 17, 20, 24, 28, 32.

\end{thebibliography}
\end{document}